\documentclass[
	english
]{scrartcl}

\usepackage[T1]{fontenc}
\usepackage[utf8]{inputenc}

\usepackage{amsmath,amsfonts,amsthm,amssymb}
\usepackage[colorlinks,citecolor=blue]{hyperref}
\usepackage{doi,natbib}
\usepackage{natbib}
\usepackage{cancel}
\usepackage{relsize}

\usepackage{changes}

\usepackage{enumitem}
\setlist[enumerate]{label=(\alph*)}

\usepackage[figure]{hypcap}
\usepackage{graphicx}
\graphicspath{{./img/}}
\usepackage{subcaption} 

\usepackage{preamble}

\usepackage{marginnote}

\usepackage[capitalise,nameinlink]{cleveref}
\allowdisplaybreaks

\numberwithin{equation}{section}

\newcommand\norm[1]{\left\Vert#1\right\Vert}
\newcommand\nnorm[1]{\Vert#1\Vert}

\newcommand\dual[2]{\left\langle #1, #2\right\rangle}

\newcommand\innerprod[2]{\left\langle #1, #2\right\rangle}
\newcommand\ninnerprod[2]{\langle #1, #2\rangle}

\newcommand\N{\mathbb{N}}
\newcommand\R{\mathbb{R}}

\newcommand\SSS{\mathcal S}
\newcommand\PPP{\mathbf P}
\newcommand\QQQ{\mathbf Q}
\newcommand\OOO{\mathbf O}
\newcommand\LLL{\mathbf\Lambda}

\newcommand\tto{\rightrightarrows}

\newcommand{\cl}{\operatorname{cl}}
\newcommand{\spa}{\operatorname{span}}

\newcommand{\sgn}{\operatorname{sgn}}

\newcommand{\dist}{\operatorname{dist}}

\newcommand{\dom}{\operatorname{dom}}
\renewcommand{\Im}{\operatorname{Im}}
\newcommand{\gph}{\operatorname{gph}}
\newcommand{\epi}{\operatorname{epi}}

\newcommand{\dr}{\operatorname{d}}

\newcommand{\trace}{\operatorname{trace}}

\DeclareMathOperator*{\argmax}{\operatorname{argmax}}

\newcommand{\xb}{\bar x}

\newcommand{\yb}{\bar y}

\DeclareMathAlphabet{\mathpzc}{OT1}{pzc}{m}{it}
\newcommand\oo{\mathpzc{o}}

\newtheorem{theorem}{Theorem}[section]
\newtheorem{lemma}[theorem]{Lemma}
\newtheorem{proposition}[theorem]{Proposition}
\newtheorem{assumption}[theorem]{Assumption}

\newtheorem{corollary}[theorem]{Corollary}
\newtheorem{remark}[theorem]{Remark}
\newtheorem{definition}[theorem]{Definition}
\newtheorem{example}[theorem]{Example}

\definecolor{mygreen}{rgb}{0.0,0.7,0.0}
\definecolor{mybrown}{rgb}{0.5,0.5,0.0}

\begin{document}

\title{On the directional asymptotic approach in optimization theory
	}
\author{%
	Mat\'{u}\v{s} Benko%
	\footnote{%
		University of Vienna,
		Applied Mathematics and Optimization,
		1090 Vienna,
		Austria;
		Johann Radon Institute for Computational and Applied Mathematics,
		4040 Linz,
		Austria;
		\email{matus.benko@ricam.oeaw.ac.at},
		\orcid{0000-0003-3307-7939}%
		}
	\and
	Patrick Mehlitz%
	\footnote{%
		University of Duisburg-Essen,
		Faculty of Mathematics,
		45127 Essen,
		Germany;
		\email{patrick.mehlitz@uni-due.de},
		\orcid{0000-0002-9355-850X}%
		}
	}

\publishers{}
\maketitle

\begin{abstract}
	 As a starting point of our research, 
	 we show that, for a fixed order $\gamma\geq 1$, each local minimizer of
	 a rather general nonsmooth optimization problem in Euclidean spaces
	 is either M-stationary in the classical sense (corresponding
	 to stationarity of order $1$), satisfies
	 stationarity conditions in terms of a coderivative construction of
	 order $\gamma$, or is asymptotically stationary with respect
	 to a critical direction as well as order $\gamma$ in a certain sense. 
	 By ruling out the latter case 
	 with a constraint qualification not stronger than directional metric
	 subregularity, we end up with new necessary optimality conditions comprising
	 a mixture of limiting variational tools of orders $1$ and $\gamma$.
	 These abstract findings are carved out for the broad class of geometric
	 constraints and $\gamma:=2$, and visualized by examples 
	 	from complementarity-constrained and nonlinear semidefinite optimization. 
	 As a byproduct of the particular setting
	 $\gamma:=1$, our general approach yields new so-called 
	 directional asymptotic regularity conditions which serve as
	 constraint qualifications guaranteeing
	 M-stationarity of local minimizers.
	 We compare these new regularity conditions
	 with standard constraint qualifications from nonsmooth optimization.
	 Further, we extend directional concepts of pseudo- and quasi-normality
	 to arbitrary set-valued mappings.
	 It is shown that these properties provide sufficient conditions for the validity
	 of directional asymptotic regularity.
	 Finally, a novel coderivative-like variational tool is used to construct
	 sufficient conditions for the presence of directional asymptotic regularity.
	 For geometric constraints, it is illustrated that all appearing objects can be 
	 calculated in terms of initial problem data.
\end{abstract}


\begin{msc}	
	\mscLink{49J52}, \mscLink{49J53}, \mscLink{49K27}, \mscLink{90C22},
	\mscLink{90C30}, \mscLink{90C33}
\end{msc}

\section{Introduction}\label{sec:introduction}

In order to identify local minimizers of optimization problems analytically
or numerically, it is
desirable that such points satisfy applicable necessary optimality conditions.
Typically, under validity of a constraint qualification, first-order
necessary optimality conditions of abstract Karush--Kuhn--Tucker (KKT)-type
hold at local minimizers. Here,
\emph{first-order} refers to the fact that first-order tools of (generalized)
differentiation are used to describe the variation of all involved data functions.
In the case where the celebrated tools of limiting variational analysis are exploited,
one speaks of so-called Mordukhovich (or, briefly, M-) stationarity,
see \cite{Mordukhovich2018}.
In the absence of constraint qualifications, i.e., in a \emph{degenerate} situation,
local minimizers still satisfy a Fritz--John (FJ)-type first-order necessary
optimality condition which allows for a potentially vanishing multiplier
associated with the generalized derivative of the objective function.
Since such a condition allows to discard the objective function, it might be
too weak in practically relevant scenarios. 

In recent years, \emph{asymptotic} (\emph{approximate} or \emph{sequential} are also common)
concepts of stationarity and regularity received much attention
not only in standard nonlinear optimization,
see \cite{AndreaniMartinezSvaiter2010,AndreaniHaeserMartinez2011,AndreaniMartinezRamosSilva2016,AndreaniMartinezRamosSilva2018},
but also in complementarity-, cardinality-, and switching-constrained programming,
see \cite{AndreaniHaeserSecchinSilva2019,KanzowRaharjaSchwartz2021,LiangYe2021,Ramos2019},
conic optimization,
see \cite{AndreaniGomezHaeserMitoRamos2021},
nonsmooth optimization,
see \cite{HelouSantosSimoes2020,Mehlitz2020a,Mehlitz2022},
or even infinite-dimensional optimization,
see \cite{BoergensKanzowMehlitzWachsmuth2019,KanzowSteckWachsmuth2018,KrugerMehlitz2021}.
The interest in asymptotic stationarity conditions is based on the observation that they
hold at local minimizers in the absence of constraint qualifications while being more restrictive
than the corresponding FJ-type conditions, and that 
different types of solution algorithms like multiplier-penalty- and some SQP-methods
naturally compute such points. Asymptotic constraint qualifications provide conditions
which guarantee that an asymptotically stationary point is already stationary in
classical sense.
It has been reported, e.g., in \cite{AndreaniMartinezRamosSilva2016,LiangYe2021,Mehlitz2020a,Ramos2019}
that asymptotic constraint qualifications are comparatively mild. 
Inherently from their construction, asymptotic constraint qualifications simplify
the convergence analysis of some numerical solution algorithms.

The aim of this paper is to apply the \emph{directional} approach
to limiting variational analysis, see e.g.\ \cite{BenkoGfrererOutrata2019}, 
in order to enrich the asymptotic
stationarity and regularity conditions from \cite{KrugerMehlitz2021,Mehlitz2020a}
with the aid of directional information.
Noting that the directional tools of variational analysis were successfully applied 
to find refined M-stationarity-type optimality conditions and mild constraint qualifications 
for diverse problems in optimization theory, see e.g.\
\cite{BaiYeZhang2019,BaiYe2021,BenkoCervinkaHoheisel2019,Gfrerer2013,Gfrerer2014,GfrererKlatte2016,GfrererYeZhou2022}
and the references therein, this seems to be a desirable goal.

\Cref{sec:directional_asymptotic_tools} contains the core of our research.
As a starting point, we show in \cref{sec:stationarity_via_coderivative}
(see, particularly, \cref{thm:higher_order_directional_asymptotic_stationarity}) that local minimizers
of rather general optimization problems in Euclidean spaces, which we formally
introduce in \cref{sec:directional_asymptotic_tools}, 
are either
M-stationary, satisfy a stationarity condition combining the limiting subdifferential
of the objective function and a coderivative-like tool associated with the
constraints of some arbitrary order $\gamma\geq 1$, 
a so-called \emph{pseudo-coderivative}, see \cite{Gfrerer2014a}, or come along with an
asymptotic stationarity condition depending on a critical direction as well as
the order $\gamma$ where the involved sequence of multipliers is diverging. 
Even for $\gamma:=1$, this enhances the findings from
\cite{KrugerMehlitz2021,Mehlitz2020a}. Furthermore, this result opens a new way on
how to come up with applicable necessary optimality conditions for the original problem,
namely, by ruling out the irregular situation of asymptotic stationarity which can
be done in the presence of so-called \emph{metric pseudo-subregularity} of order $\gamma$,
see \cite{Gfrerer2014a} again.
In the case $\gamma:=1$, we end up with M-stationarity, and metric pseudo-subregularity reduces
to metric subregularity, i.e., we obtain results related to \cite{Gfrerer2013}.
For $\gamma>1$, this procedure leads to a mixed-order stationarity condition 
involving the pseudo-coderivative of order $\gamma$, and metric pseudo-subregularity
is weaker than metric subregularity.
If $\gamma:=2$ and so-called geometric constraints, induced by a twice continuously
differentiable mapping $g$ as well as a closed set $D$, are investigated, 
this pseudo-coderivative
can be estimated from above in terms of initial problem data, i.e., 
in terms of (first- and second-order) derivatives associated with $g$ 
as well as tangent and normal cones to $D$, under mild conditions.
These estimates of the pseudo-coderivative of order $2$ are interesting on their own 
and presented in \cref{sec:variational_analysis_constraint_mapping},
which is the essence to all applications of our general findings.
The associated mixed-order necessary optimality conditions and qualification conditions are worked
out in \cref{sec:constraint_mappings}, and in \cref{sec:applications}, 
	they are applied to complementarity-constrained 
	and nonlinear semidefinite optimization problems
in order to illustrate our findings.
Let us note that related necessary optimality conditions for 
optimization problems which comprise 
first- and second-order tools at the same time can be found e.g.\ in
\cite{ArutyunovAvakovIzmailov2008,Avakov1985,Avakov1989,AvakovArutunovIzmailov2007,Gfrerer2007,Gfrerer2014a,IzmailovSolodov2002}.
These results are based on the concept of $2$--regularity and its extensions, see \cite{Avakov1985,TRETYAKOV1984} for its origins.
Indeed, even Gfrerer's metric pseudo-subregularity from \cite{Gfrerer2014a}, utilized in this paper,
can be seen as an extension of $2$--regularity to arbitrary set-valued mappings.
For us, however, these mixed-order conditions are only a by-product - we focus on how they
can be used to find new constraint qualifications guaranteeing M-stationarity of local minimizers.

\Cref{sec:directional_asymptotic_regularity} is dedicated to the investigation of \emph{directional} asymptotic regularity conditions,
which are motivated by the asymptotic stationarity conditions from \cref{thm:higher_order_directional_asymptotic_stationarity} 
(for $\gamma:=1$) and whose validity
directly yields M-stationarity of local minimizers.
Roughly speaking, these conditions demand certain control of unbounded input 
sequences (multipliers) associated with the regular coderivative of the underlying
set-valued mapping in a neighborhood of the reference point.
We enrich and refine the asymptotic regularity conditions from \cite{Mehlitz2020a} in two ways.
First, the directional approach reveals that asymptotic regularity is only necessary in critical directions.
Second, we observe an additional restriction the problematic multipliers satisfy:
while their norm tends to infinity, their direction is tightly controlled.
These insights enable us to
relate our new constraint qualifications with already existing ones
from the literature.
Similarly as standard asymptotic regularity, the directional counterpart is also independent of (directional) metric subregularity.
However, several sufficient conditions for metric subregularity, which are independent of asymptotic regularity,
imply directional asymptotic regularity.
For instance, this is true for the First-Order Sufficient Condition for Metric Subregularity
from \cite{GfrererKlatte2016}, see \cref{sec:directional_asymptotic_regularity_basics}.
Moreover, in \cref{sec:pseudo_quasi_normality}, we extend the (directional) concepts of pseudo- and quasi-normality
from \cite{BaiYeZhang2019,BenkoCervinkaHoheisel2019} to abstract set-valued mappings
and show that these conditions are sufficient
for directional metric subregularity as well as directional asymptotic regularity.
Notably, even standard (nondirectional) versions of pseudo- and quasi-normality
do not imply asymptotic regularity since the latter does not restrict the direction of the problematic multipliers.
Finally, a new directional coderivative-like tool, the \emph{directional super-coderivative},
see \cref{sec:generalized_differentiation},
is used in \cref{sec:asymptotic_regularity_via_super_coderivative} to construct
sufficient conditions for the validity of directional asymptotic regularity.
In the presence of so-called metric pseudo-regularity, 
see \cite{Gfrerer2014a} again, this leads to conditions in terms of the aforementioned
pseudo-coderivatives. Noting that these generalized derivatives can be computed
in terms of initial problem data
for geometric constraint systems,
we can specify our findings in this situation.
As it turns out, the approach 
is closely related to our findings from \cref{sec:constraint_mappings}.
Furthermore, we show that the explicit sufficient conditions for directional asymptotic regularity
provide constraint qualifications for M-stationarity which are not stronger than the
First- and Second-Order Sufficient Condition for Metric Subregularity from \cite{GfrererKlatte2016}.

\section{Notation and preliminaries}\label{sec:notation}

We rely on standard notation taken from
\cite{AubinFrankowska2009,BonnansShapiro2000,RockafellarWets1998,Mordukhovich2018}.

\subsection{Basic notation}

Let $\R$, $\R_+$, and $\R_-$ denote the real, the nonnegative real, and the nonpositive 
real numbers, respectively.
The sign function $\sgn\colon\R\to\{-1,0,1\}$ is defined by
$\sgn(t):=-1$ for all $t<0$, $\sgn(t):=1$ for all $t>0$, and $\sgn(0):=0$.
Throughout the paper, $\mathbb X$ and $\mathbb Y$ denote Euclidean spaces, 
i.e., finite-dimensional
Hilbert spaces. For simplicity, the associated inner product will be represented by
$\innerprod{\cdot}{\cdot}$ since the underlying space will be clear from the context.
The norm induced by the inner product is denoted by $\norm{\cdot}$.
The unit sphere in $\mathbb X$ will be represented by $\mathbb S_{\mathbb X}$.
Furthermore, for $\varepsilon>0$ and $\bar x\in\mathbb X$,
$\mathbb B_\varepsilon(\bar x):=\{x\in\mathbb X\,|\,\norm{x-\bar x}\leq\varepsilon\}$
is the closed $\varepsilon$-ball around $\bar x$.
We are also concerned with so-called (closed) directional neighborhoods of given directions.
These are sets of type
\[
 	\mathbb B_{\varepsilon,\delta}(u)
	:=
	\{w\in\mathbb X\,|\,\norm{\norm{w}u-\norm{u}w}\leq\delta\norm{u}\norm{w},\,\norm{w}\leq\varepsilon\},
\]
where $u\in\mathbb X$ is a reference direction and $\varepsilon,\delta>0$. Clearly, $\mathbb B_{\varepsilon,\delta}(0)=\mathbb B_{\varepsilon}(0)$.
For a nonempty set $Q\subset\mathbb X$, we refer to the closed, convex cone
$Q^\circ:=\{\eta\in\mathbb X\,|\,\forall x\in Q\colon\,\innerprod{\eta}{x}\leq 0\}$
as the polar cone of $Q$.
Furthermore, for some $\bar x\in\mathbb X$, 
$[\bar x]^\perp:=\{\eta\in\mathbb X\,|\,\innerprod{\eta}{\bar x}=0\}$ 
and $\spa(\bar x)$  are the annihilator of $\bar x$
and the smallest subspace of $\mathbb X$ containing $\bar x$, respectively.
By $\dist(\bar x,Q):=\inf_x\{\norm{x-\bar x}\,|\,x\in Q\}$,
we denote the distance of $\bar x$ to $Q$.
For simplicity of notation, we use $\bar x+Q:=Q+\bar x:=\{x+\bar x\in\mathbb X\,|\,x\in Q\}$.
	The closure and the horizon cone of $Q$ are represented by
	$\cl(Q)$ and $Q^\infty$, respectively.
For a given linear operator $A\colon\mathbb X\to\mathbb Y$, 
$A^*\colon\mathbb Y\to\mathbb X$ is used to denote its adjoint
	while $\Im A:=\{Ax\in\mathbb Y\,|\,x\in\mathbb X\}$ is the image of $A$.

Let $g\colon\mathbb X\to\mathbb Y$ be a continuously differentiable mapping.
We use $\nabla g(\bar x)\colon\mathbb X\to\mathbb Y$ to denote the derivative
of $g$ at $\bar x\in\mathbb X$. Note that $\nabla g(\bar x)$ is a linear operator.
Let us emphasize that, in the special case $\mathbb Y:=\R$, $\nabla g(\bar x)$
does not coincide with the standard gradient which would correspond to
$\nabla g(\bar x)^*1$.
For twice continuously differentiable $g$ and a vector $\lambda\in\mathbb Y$,
we set $\langle \lambda,g\rangle(x):=\innerprod{\lambda}{ g(x) }$ 
for each $x\in\mathbb X$ in order
to denote the associated scalarization mapping $\innerprod{\lambda}{g}\colon\mathbb X\to\R$.
By $\nabla\innerprod{\lambda}{g}(\bar x)$ and $\nabla^2\innerprod{\lambda}{g}(\bar x)$
we represent the first- and second-order derivatives of this map at $\bar x\in\mathbb X$ (w.r.t.\ the variable which enters $g$).
Furthermore, for $u,u'\in\mathbb X$, we make use of
\[
	\nabla^2g(\bar x)[u,u']
	:=
	\sum_{i=1}^m \innerprod{u}{\nabla^2\innerprod{e_i^\textup{c}}{g}(\bar x)(u')}\,e_i^\textup{c}
\]
for brevity where 
$m\in\N$ is the dimension of $\mathbb Y$ and
$e_1^\textup{c},\ldots,e_m^\textup{c}\in\mathbb Y$ denote the $m$ canonical unit vectors of $\mathbb Y$.
In the case $\mathbb Y:=\R$, the second-order derivative 
$\nabla ^2g(\bar x)\colon\mathbb X\times\mathbb X\to\R$ is a bilinear mapping, and
for each $u\in\mathbb X$, we identify $\nabla^2g(\bar x)u$ with an element of $\mathbb X$.

\subsection{Fundamentals of variational analysis}

Let us fix a closed set $Q\subset\mathbb X$ and some point $x\in Q$.
We use
\begin{align*}
	\mathcal T_Q(x)
	:=
	\left\{u\in\mathbb X\,\middle|\,
		\begin{aligned}
			&\exists\{u_k\}_{k\in\N}\subset\mathbb X,\,\exists\{t_k\}_{k\in\N}\subset\R_+\colon\\
			&\qquad u_k\to u,\,t_k\downarrow 0,\,x+t_ku_k\in Q\,\forall k\in\N
		\end{aligned}
	\right\}
\end{align*}
to denote the (Bouligand) \emph{tangent cone} to $Q$ at $x$.
Furthermore, we make use of
\begin{align*}
	\widehat{\mathcal N}_Q(x)
	&:=
	\left\{\eta\in\mathbb X\,|\,\forall x'\in Q\colon\,\ninnerprod{\eta}{x'-x}\leq\oo(\nnorm{x'-x})\right\},\\
	\mathcal N_Q(x)
	&:=
	\left\{\eta\in\mathbb X\,\middle|\,
		\begin{aligned}
			&\exists\{x_k\}_{k\in\N}\subset Q,\,\exists\{\eta_k\}_{k\in\N}\subset\mathbb X\colon\\
			&\qquad x_k\to x,\,\eta_k\to\eta,\,\eta_k\in\widehat{\mathcal N}_Q(x_k)\,\forall k\in\N
		\end{aligned}
	\right\},
\end{align*}
the \emph{regular} (or Fr\'{e}chet) and \emph{limiting} (or Mordukhovich) \emph{normal cone} to $Q$ at $x$.
Observe that both of these normal cones coincide with the standard normal cone of convex
analysis as soon as $Q$ is convex.
For $\tilde x\notin Q$, we set $\mathcal T_Q(\tilde x):=\emptyset$ and
$\widehat{\mathcal N}_Q(\tilde x):=\mathcal N_Q(\tilde x):=\emptyset$.
Finally, for some $u\in\mathbb X$, we use
\[
	\mathcal N_Q(x;u)
	:=
	\left\{
		\eta\in\mathbb X\,\middle|\,
			\begin{aligned}
				&\exists\{u_k\}_{k\in\N}\subset\mathbb X,\,\exists\{t_k\}_{k\in\N}\subset\R_+,\,
				\exists\{\eta_k\}_{k\in\N}\subset\mathbb X\colon\\
				&\qquad u_k\to u,\,t_k\downarrow 0,\,\eta_k\to\eta,\,\eta_k\in\widehat{\mathcal N}_Q(x+t_ku_k)\,\forall k\in\N
			\end{aligned}
	\right\}
\]
in order to represent the \emph{directional limiting normal cone} to $Q$ at $x$ in direction $u$.
Note that this set is empty if $u$ does not belong to $\mathcal T_Q(x)$. 
If $Q$ is convex, we have $\mathcal N_Q(x;u)=\mathcal N_Q(x)\cap [u]^\perp$.

The limiting normal cone to a set is well known for its \emph{robustness}, i.e., 
it is outer semicontinuous as a set-valued mapping.
In the course of the paper, we exploit an analogous property of the directional limiting 
normal cone which has been validated in \cite[Proposition~2]{GfrererYeZhou2022}.
\begin{lemma}\label{lem:robustness_directional_limiting_normals}
	Let $Q\subset\mathbb X$ be closed and fix $x\in Q$.
	Then, for each $u\in\mathbb X$, we have
	\[
	\mathcal N_Q(x;u)
	=
	\left\{
		\eta\in\mathbb X\,\middle|\,
			\begin{aligned}
				&\exists\{u_k\}_{k\in\N}\subset\mathbb X,\,
					\exists\{t_k\}_{k\in\N}\subset\R_+,\,
					\exists\{\eta_k\}_{k\in\N}\subset\mathbb X\colon\\
				&\qquad u_k\to u,\,t_k\downarrow 0,\,\eta_k\to\eta,\,
				\eta_k\in \mathcal N_Q(x+t_ku_k)\,\forall k\in\N
			\end{aligned}
	\right\}.
\]
\end{lemma}

In this paper, the concept of \emph{polyhedrality} will be of essential importance.
Let us recall that a set $Q\subset\R^m$ will be called polyhedral if it is
the union of finitely many convex polyhedral sets. Similarly, it is referred to
as locally polyhedral around $x\in Q$ whenever 
$Q\cap\{z\in\R^m\,|\,\forall i\in\{1,\ldots,m\}\colon\,|z_i-x_i|\leq\varepsilon\}$
is polyhedral for some $\varepsilon>0$.
The following lemma provides some basic properties of polyhedral sets.
Statement~\ref{item:exactness_tangential_approximation}
is proven in \cite[Proposition 8.24]{Ioffe2017}.
The equality in statement~\ref{item:normal_cones_to_polyhedral_sets} 
follows from \cite[Proposition 2.11]{BenkoGfrererYeZhangZhou2022}
and the rest is straightforward,
see \cite[Lemma~2.1]{Gfrerer2014} as well.
\begin{lemma}\label{lem:some_properties_of_polyhedral_sets}
	Let $Q\subset\mathbb \R^m$ be a closed set which is locally polyhedral
	around some fixed point $x\in Q$.
	Then the following statements hold.
	\begin{enumerate}
	\item\label{item:exactness_tangential_approximation}
		There exists a neighborhood $U\subset\R^m$ of $x$ 
		such that $(x+\mathcal T_Q(x))\cap U=Q\cap U$.
	\item\label{item:normal_cones_to_polyhedral_sets} 
		For arbitrary $u\in\R^m$, we have
		\begin{equation}\label{eq:normals_to_polyhedral_sets}
			\mathcal N_Q(x;u)
			=
			\mathcal N_{\mathcal T_Q(x)}(u)
			\subset
			\mathcal N_{Q}(x)\cap [u]^\perp.
		\end{equation}
		If $Q$ is, additionally, convex, and $u\in\mathcal T_Q(x)$, 
		then the final inclusion holds as an equality.
	\end{enumerate}
\end{lemma}

It is well known that the regular and limiting normal cone enjoy an exact product rule
which is not true for the tangent cone in general.
However, the following lemma shows that such a product rule also holds for tangents
as soon as polyhedral sets are under consideration.
Its proof is straightforward and, hence, omitted.
\begin{lemma}\label{lem:product_rule_tangents_polyhedral_sets}
\leavevmode
	\begin{enumerate}
		\item For closed sets $P\subset\mathbb X$ and $Q\subset\mathbb Y$ 
			as well as $x\in P$ and $y\in Q$,
			we have $\mathcal T_{P\times Q}(x,y)\subset\mathcal T_P(x)\times\mathcal T_Q(y)$.
		\item For closed sets $P\subset\R^n$ and $Q\subset\R^m$ as well as $x\in P$ and $y\in Q$,
			such that $P$ and $Q$ are locally polyhedral around $x$ and $y$, respectively, 
			we have $\mathcal T_{P\times Q}(x,y)=\mathcal T_P(x)\times\mathcal T_Q(y)$.
	\end{enumerate}
\end{lemma}

Let us mention that a slightly more general version of the above lemma 
can be found in \cite[Proposition~1]{GfrererYe2017}.

For a set-valued mapping $\Phi\colon\mathbb X\tto\mathbb Y$, we use
$\dom \Phi:=\{x\in\mathbb X\,|\,\Phi(x)\neq\emptyset\}$,
$\gph \Phi:=\{(x,y)\in\mathbb X\times\mathbb Y\,|\,y\in\Phi(x)\}$,
$\ker\Phi:=\{x\in\mathbb X\,|\,0\in\Phi(x)\}$, and
$\Im\Phi:=\bigcup_{x\in\mathbb X}\Phi(x)$ in order to
represent the domain, graph, kernel, and image of $\Phi$, respectively.
Furthermore, the so-called inverse mapping $\Phi^{-1}\colon\mathbb Y\tto\mathbb X$
is defined via $\gph\Phi^{-1}:=\{(y,x)\in\mathbb Y\times\mathbb X\,|\,(x,y)\in\gph\Phi\}$.

There exist numerous concepts of local \emph{regularity} or
\emph{Lipschitzness} associated with set-valued mappings.
In this paper, we are mostly concerned with so-called
directional metric pseudo-(sub)regularity which originates from 
\cite[Definition~1]{Gfrerer2014a}.

\begin{definition}\label{def:metric_pseudo_subregularity}
	Fix a set-valued mapping $\Phi\colon\mathbb X\tto\mathbb Y$ which has a closed graph locally
	around $(\bar x,\bar y)\in\gph\Phi$, a pair of directions $(u,v)\in\mathbb X\times\mathbb Y$, and 
	a constant $\gamma\geq 1$.
	\begin{enumerate}
		\item 
		We say that $\Phi$ is \emph{metrically pseudo-regular of order $\gamma$ in
		direction $(u,v)$} at $(\bar x,\bar y)$ if there are constants $\varepsilon>0$, $\delta>0$,
		and $\kappa>0$ such that the estimate
		\begin{equation}\label{eq:estimate_pseudo_regularity}
			\norm{x-\bar x}^{\gamma-1}\dist(x,\Phi^{-1}(y))
			\leq
			\kappa\,\dist(y,\Phi(x))
		\end{equation}
		holds for all $(x,y)\in(\bar x,\bar y)+\mathbb B_{\varepsilon,\delta}(u,v)$ with
		$\dist(y,\Phi(x))\leq\delta\norm{x-\bar x}^\gamma$.
		In the case where this is fulfilled for $(u,v):=(0,0)$, we say that $\Phi$ is 
		\emph{metrically pseudo-regular of order $\gamma$} at $(\bar x,\bar y)$.
		\item We say that $\Phi$ is \emph{metrically pseudo-subregular of order $\gamma$ in direction $u$}
		at $(\bar x,\bar y)$ if there are constants $\varepsilon>0$, $\delta>0$, and $\kappa>0$
		such that \eqref{eq:estimate_pseudo_regularity} holds for $y:=\bar y$ and all
		$x\in\bar x+\mathbb B_{\varepsilon,\delta}(u)$.
		In the case where this is fulfilled for $u:=0$, we say that $\Phi$ is
		\emph{metrically pseudo-subregular of order $\gamma$} at $(\bar x,\bar y)$.
	\end{enumerate}
\end{definition}

Metric pseudo-regularity of order $\gamma\geq 1$ in direction $(u,0)$ at $(\bar x,\bar y)$
is a sufficient condition for metric pseudo-subregularity of order $\gamma$ in direction $u$ at the
same point, see \cite[Lemma~3]{Gfrerer2014a}.
Observe that metric pseudo-subregularity in a specified direction of some order $\gamma\geq 1$
implies metric pseudo-subregularity of arbitrary order larger than $\gamma$ in the same direction.
For $\gamma:=1$, the above definition of (directional) metric pseudo-subregularity 
recovers the one of \emph{(directional) metric subregularity}, 
see \cite[Definition~1.2]{Gfrerer2013}.
On the contrary, for $\gamma:=1$, the above definition of directional metric pseudo-regularity
does not recover the one of \emph{directional metric regularity} which demands
that \eqref{eq:estimate_pseudo_regularity} holds 
for all $(x,y)\in(\bar x,\bar y)+\mathbb B_{\varepsilon,\delta}(u,v)$
such that $\norm{(u,v)}\dist((x,y),\gph\Phi)\leq\delta\norm{(u,v)}\norm{(x,y)-(\bar x,\bar y)}$,
see \cite[Definition~1.1]{Gfrerer2013}.
Particularly, for $(u,v):=(0,0)$, the notion of directional metric regularity
reduces to the classical one of metric regularity,
while directional metric pseudo-regularity does not.
This was shown in \cite[Example~1.1]{Gfrerer2014a},
which is a very natural example, and we will use it to illustrate some novel concepts.

\begin{example}\label{ex:metric_pseudo_regularity}
For every  $\gamma \geq 1$,
the mapping $\Phi \colon \R \tto \R$, 
given by $\Phi(x) := \{ \vert x \vert^\gamma \}$, $x\in\R$,
is metrically pseudo-regular of order $\gamma$ at $(0,0)$. 
The case $\gamma := 1$ provides an example of a mapping
which is metrically pseudo-regular of order $1$ at $(0,0)$ but not metrically regular there.
The violation of metric regularity is clear
as any points $y<0$ approaching $0$ come along with $\Phi^{-1}(y) = \emptyset$,
blowing up the left-hand side of \eqref{eq:estimate_pseudo_regularity}.
These problematic elements $y$ are, however,
ruled out by the condition $\dist(y,\Phi(x))\leq\delta\norm{x-\bar x}^\gamma$ 
in the definition of metric pseudo-regularity,
which reads $\vert y - \vert x \vert \vert \leq \delta \vert x \vert$ in the 
present situation.
\end{example}

Another important case, which we will explore in detail, corresponds to $\gamma:=2$.
In this case, the notions from \cref{def:metric_pseudo_subregularity} 
provide an extension of so-called \emph{2-regularity} from \cite{Avakov1985,TRETYAKOV1984}
to set-valued mappings.
In \cref{sec:comparison_2reg}, we compare our approach with an extension of 2-regularity to
constraint mappings from \cite{ArutyunovAvakovIzmailov2008,ArutyunovIzmailov2020}.

Recall that a single-valued function $g\colon\mathbb X\to\mathbb Y$ is called \emph{calm in direction} $u\in\mathbb X$ 
at $x\in\mathbb X$ whenever there 
are constants $\varepsilon>0$, $\delta>0$, and $L>0$ such that
\[
	\forall x'\in x+\mathbb B_{\varepsilon,\delta}(u)\colon\quad
	\nnorm{g(x')-g(x)}\leq L\nnorm{x'-x}.
\]
If this holds for $u:=0$, we simply say that $g$ is {calm} at $x$.
Clearly, the latter property is weaker than Lipschitzness of $g$ at $x$.

\subsection{Generalized differentiation}\label{sec:generalized_differentiation}

In this section, we recall some notions from generalized differentiation and
introduce some novel derivatives for set-valued mappings.

\subsubsection{Subdifferentials}

Let us start with a lower semicontinuous function $\varphi\colon\mathbb X\to\R\cup\{\infty\}$
and some point $\bar x\in\dom\varphi:=\{x\in\mathbb X\,|\,\varphi(x)<\infty\}$.
	The lower semicontinuous function 
	$\dr\varphi(\bar x)\colon\mathbb X\to\R\cup\{-\infty,\infty\}$ given by
	\[
		\forall u\in\mathbb X\colon\quad
		\dr\varphi(\bar x)(u)
		:=
		\liminf\limits_{t\downarrow 0,\,u'\to u}
		\frac{\varphi(\bar x+tu')-\varphi(\bar x)}{t}
	\]
	is referred to as the \emph{subderivative} of $\varphi$ at $\bar x$.
The \emph{regular} (or Fr\'{e}chet) and \emph{limiting} (or Mordukhovich) \emph{subdifferential} 
of $\varphi$ at $\bar x$ are given by
\begin{align*}
	\widehat\partial \varphi(\bar x)
	&:=
	\left\{
		\eta\in\mathbb X\,\middle|\,
		(\eta,-1)\in\widehat{\mathcal N}_{\epi\varphi}(\bar x,\varphi(\bar x))
	\right\},\\
	\partial \varphi(\bar x)
	&:=
	\left\{
		\eta\in\mathbb X\,\middle|\,
		(\eta,-1)\in\mathcal N_{\epi\varphi}(\bar x,\varphi(\bar x))
	\right\},
\end{align*}
respectively, where $\epi\varphi:=\{(x,\alpha)\in\mathbb X\times\R\,|\,\varphi(x)\leq\alpha\}$
is the epigraph of $\varphi$.
In the case where $\varphi$ is continuously differentiable at $\bar x$, both sets
reduce to the singleton containing only the gradient $\nabla\varphi(\bar x)^*1$.
We note that for any sequences $\{x_k\}_{k\in\N}\subset\dom\varphi$ and $\{x_k^*\}_{k\in\N}\subset\mathbb X$
such that $x_k\to\bar x$, $\varphi(x_k)\to\varphi(\bar x)$, $x_k^*\to x^*$ for some $x^*\in\mathbb X$, and
$x_k^*\in\partial\varphi(x_k)$ for each $k\in\N$, we also have $x^*\in\partial\varphi(\bar x)$,
see \cite[Proposition~8.7]{RockafellarWets1998}.
This property is referred to as robustness of the limiting subdifferential.

In the case where $\varphi$ is locally Lipschitzian around $\bar x$,
and for some direction $u\in\mathbb X$,
\[
	\partial\varphi(\bar x;u)
	:=
	\left\{
		\eta\in\mathbb X\,\middle|\,
			\begin{aligned}
				&\exists\{u_k\}_{k\in\N}\subset\mathbb X,\,
				\exists\{t_k\}_{k\in\N}\subset\R_+,\,
				\exists\{\eta_k\}_{k\in\N}\subset\mathbb X\colon
				\\
				&\qquad
				u_k\to u,\,t_k\downarrow 0,\,\eta_k\to\eta,\,
				\eta_k\in\widehat\partial\varphi(\bar x+t_ku_k)\,\forall k\in\N
			\end{aligned}
	\right\}
\]
is referred to as the \emph{limiting subdifferential} of $\varphi$ at $\bar x$
\emph{in direction $u$}.
We note that $\partial\varphi(\bar x;0)=\partial\varphi(\bar x)$
and $\partial\varphi(\bar x;u)\subset\partial\varphi(\bar x)$
for all $u\in\mathbb X$.
Furthermore, let us mention that, 
in the definition of the directional limiting subdifferential,
we can equivalently replace the requirement 
$\eta_k\in\widehat\partial\varphi(\bar x+t_ku_k)$
by
$\eta_k\in\partial\varphi(\bar x+t_ku_k)$
for each $k\in\N$.
This can be easily checked by means of a classical diagonal sequence argument.
Hence, the directional limiting subdifferential also enjoys a certain 
kind of robustness.

\subsubsection{Graphical derivatives}

Below, we introduce three different graphical derivatives of a set-valued mapping. 
While the standard graphical derivative is well known from the literature, the
concepts of graphical pseudo-derivative and graphical subderivative are,
to the best of our knowledge, new.

\begin{definition}\label{def:graphical_derivative}
	Let $\Phi\colon\mathbb X\tto\mathbb Y$ be a set-valued mapping possessing 
	a closed graph locally around $(\xb,\yb)\in\gph \Phi$.
	\begin{enumerate}
	\item 
		The \emph{graphical derivative} of $\Phi$ at $(\bar x,\bar y)$ is the
		mapping $D\Phi(\bar x,\bar y)\colon\mathbb X\tto\mathbb Y$ given by
		\[
			\gph D\Phi(\bar x,\bar y)=\mathcal T_{\gph\Phi}(\bar x,\bar y).
		\]
		In the case where $\Phi$ is single-valued at $\bar x$, we use 
		$D\Phi(\bar x)\colon\mathbb X\tto\mathbb Y$ for brevity.
        \item\label{item:def_graph_pseudo-derivative}
            Given $\gamma \geq 1$, the \emph{graphical pseudo-derivative of order $\gamma$}
            of $\Phi$ at $(\bar x,\bar y)$
            is the mapping 
		$D_{\gamma}\Phi(\bar x,\bar y)\colon \mathbb X \tto \mathbb Y$ 
		which assigns to $u\in \mathbb X$ the set of all $v\in \mathbb Y$
                such that there are sequences $\{u_k\}_{k\in\N}\subset\mathbb X$, 
		$\{v_k\}_{k\in\N}\subset\mathbb Y$, 
		and $\{t_k\}_{k\in\N}\subset\R_+$
		which satisfy $u_k\to u$, $v_k\to v$, $t_k\downarrow 0$, and 
		$(\bar x+t_ku_k,\bar y + (t_k \norm{u_k})^{\gamma} v_k)\in\gph\Phi$ 
		for all $k\in\N$.
	\item\label{item:def_graph_subderivative}
		The \emph{graphical subderivative} of $\Phi$ at $(\bar x,\bar y)$
		is the mapping 
		$D_\textup{sub}\Phi(\bar x,\bar y)\colon \mathbb S_{\mathbb X} \tto \mathbb S_{\mathbb Y}$ 
		which assigns to $u\in \mathbb S_{\mathbb X}$ 
		the set of all $v\in \mathbb S_{\mathbb Y}$
		such that there are sequences $\{u_k\}_{k\in\N}\subset\mathbb X$, 
		$\{v_k\}_{k\in\N}\subset\mathbb Y$, 
		and $\{t_k\}_{k\in\N},\{\tau_k\}_{k\in\N}\subset\R_+$
		which satisfy $u_k\to u$, $v_k\to v$, $t_k\downarrow 0$, $\tau_k\downarrow 0$, 
		$\tau_k/t_k\to\infty$, and 
		$(\bar x+t_ku_k,\bar y+\tau_kv_k)\in\gph\Phi$ 
		for all $k\in\N$.
	\end{enumerate}
\end{definition}

Let us note that for every set-valued mapping $\Phi\colon\mathbb X\tto\mathbb Y$, 
whose graph is closed locally around $(\bar x,\bar y)\in\gph\Phi$, 
	we have $D_1\Phi(\bar x,\bar y)(u)=D\Phi(\bar x,\bar y)(u)$ 
	for all $u\in\mathbb S_{\mathbb X}$.
	Furthermore, for each $\gamma > 1$,
	one obtains the trivial estimates
	\begin{align*}
		\dom D_\gamma\Phi(\bar x,\bar y)
		\subset
		\ker D\Phi(\bar x,\bar y)
	\end{align*}
	and
\begin{equation}\label{eq:trivial_upper_estimate_graphical_subderivative}
	\forall u\in\mathbb S_{\mathbb X}\colon\quad
	D_\textup{sub}\Phi(\bar x,\bar y)(u)
	\subset 
	D\Phi(\bar x,\bar y)(0)
\end{equation}
right from the definition of these objects. 

In the course of the paper, we are mainly interested in the graphical (sub)derivative 
associated with so-called normal cone mappings. 
In the next lemma, we present some corresponding upper estimates.

\begin{lemma}\label{lem:graphical_derivatives_of_normal_cone_map}
	Let $D\subset\mathbb Y$ be a nonempty, closed, convex set such that
	the (single-valued) projection operator onto $D$, denoted by
	$\Pi_D\colon\mathbb Y\to\mathbb Y$, is directionally differentiable.
	Fix $\bar y\in D$ and $\bar y^*\in\mathcal N_D(\bar y)$.
	Then, for arbitrary $u\in\mathbb Y$, we find
	\begin{align*}
		D\mathcal N_D(\bar y,\bar y^*)(u)
		&\subset
		\{v\in\mathbb Y\,|\,\Pi_D'(\bar y+\bar y^*;u+v)=u\},
	\end{align*}
	and for $u\in\mathbb S_{\mathbb Y}$, we find
	\begin{align*}
		D_\textup{sub}\mathcal N_D(\bar y,\bar y^*)(u)
		&\subset
		\{v\in\mathbb S_{\mathbb Y}\,|\,
			\Pi_D'(\bar y+\bar y^*;v)=0, \innerprod{u}{v}\geq 0\}.
	\end{align*}
	Above, $\Pi_D'(y,v)$ denotes the directional derivative of $\Pi_D$
	at $y\in\mathbb Y$ in direction $v\in\mathbb Y$.
\end{lemma}
\begin{proof}
	By convexity of $D$, we have the well-known equivalence
	\[
		\forall y,y^*\in\mathbb Y\colon\quad
		y^*\in\mathcal N_D(y)
		\quad\Longleftrightarrow\quad
		\Pi_D(y+y^*)=y.
	\]
	In the remainder of the proof, we set $\tilde y:=\bar y+\bar y^*$ for brevity.
	Next, let us fix $u,v\in\mathbb Y$ as well as 
	sequences $\{u_k\}_{k\in\N},\{v_k\}_{k\in\N}\subset\mathbb Y$
	and $\{\tau_k\}_{k\in\N},\{\varepsilon_k\}_{k\in\N}\subset\R_+$
	such that $u_k\to u$, $v_k\to v$, $\tau_k\downarrow 0$, 
	and $\bar y^*+\tau_kv_k\in\mathcal N_D(\bar y+\tau_k\varepsilon_k u_k)$, i.e.,
	$\Pi_D(\tilde y+\tau_k\varepsilon_ku_k+\tau_kv_k)=\bar y+\tau_k\varepsilon_ku_k$,
	for each $k\in\N$. Using $\Pi_D(\tilde y)=\bar y$, we find
	\begin{equation}\label{eq:graphical_derivative_normal_cone_map}
		\forall k\in\N\colon\quad
		\varepsilon_k u_k
		=
		\frac{\Pi_D(\tilde y+\tau_k\varepsilon_ku_k+\tau_kv_k)-\Pi_D(\tilde y)}
		{\tau_k}.
	\end{equation}
	
	In the case where $v\in D\mathcal N_D(\bar y,\bar y^*)(u)$ holds, we can choose $\varepsilon_k=1$
	for each $k\in\N$, and taking the limit $k\to\infty$ in
	\eqref{eq:graphical_derivative_normal_cone_map} 
	while exploiting directional differentiability
	and Lipschitzness of $\Pi_D$ yields $\Pi'_D(\tilde y;u+v)=u$.
	This shows the first estimate.
	
	Now, assume that $v\in D_\textup{sub}\mathcal N_D(\bar y,\bar y^*)(u)$ is valid.
	Then $\varepsilon_k\downarrow 0$ and $u, v \in\mathbb S_{\mathbb Y}$ can be postulated,
	and taking the limit $k\to\infty$ in \eqref{eq:graphical_derivative_normal_cone_map} 
	shows $\Pi'_D(\tilde y;v)=0$.
	By nature of the projection, we have
	\begin{align*}
		\innerprod{
			\tilde y+\tau_k\varepsilon_ku_k+\tau_kv_k
			-
			\Pi_D(\tilde y+\tau_k\varepsilon_ku_k+\tau_kv_k)
		}{
			\Pi_D(\tilde y)
			-
			\Pi_D(\tilde y+\tau_k\varepsilon_ku_k+\tau_kv_k)
		}
		\leq 
		0
	\end{align*}
	for each $k\in\N$.
	Exploiting \eqref{eq:graphical_derivative_normal_cone_map}, this is equivalent to
	\begin{align*}
		\innerprod{
			\tilde y+\tau_kv_k-\Pi_D(\tilde y)
		}{
			\Pi_D(\tilde y)
			-
			\Pi_D(\tilde y+\tau_k\varepsilon_ku_k+\tau_kv_k)
		}
		\leq 
		0
	\end{align*}
	for each $k\in\N$.
	Some rearrangements and the characterization of the projection lead to
	\begin{align*}
		&\tau_k\,
		\innerprod{v_k}{
			\Pi_D(\tilde y)
			-
			\Pi_D(\tilde y+\tau_k\varepsilon_ku_k+\tau_kv_k)
		}
		\\
		&\quad
		\leq
		\innerprod{
			\tilde y-\Pi_D(\tilde y)
		}{
			\Pi_D(\tilde y+\tau_k\varepsilon_ku_k+\tau_kv_k)
			-
			\Pi_D(\tilde y)
		}
		\leq 
		0.
	\end{align*}
	Division by $\tau_k^2\varepsilon_k$ and \eqref{eq:graphical_derivative_normal_cone_map},
	thus, give us $\innerprod{v_k}{ u_k }\geq 0$
	for each $k\in\N$, and taking the limit, we obtain $\innerprod{u}{v}\geq 0$
	which shows the second estimate.
\end{proof}

Let us note that it has been shown in \cite[Theorem~3.1, Corollary~3.1]{WuZhangZhang2014} 
that the estimate on the
graphical derivative of the normal cone mapping $\mathcal N_D$ holds as an equality 
in the situation where
$D$ is the convex cone of positive semidefinite symmetric matrices, 
and that the presented proof extends to arbitrary convex
cones as long as the associated projection operator is directionally differentiable.
This result can also be found in slightly more general form in 
\cite[Theorem~3.3]{MordukhovichOutrataRamirez2015}.
In order to make the estimates from \cref{lem:graphical_derivatives_of_normal_cone_map} 
explicit, one needs to be in
position to characterize the directional derivative of the projection onto the convex set $D$.
This is easily possible if $D$ is polyhedral, see \cite{Haraux1977} and 
\cref{rem:comparison_estimate_pseudocoderivative_polyhedral}, 
but even in
nonpolyhedral situations, e.g., where $D$ is the second-order cone 
or the cone of positive semidefinite symmetric
matrices, closed formulas for this directional derivative are available in the literature, see
\cite[Lemma~2]{OutrataSun2008} and \cite[Theorem~4.7]{SunSun2002}, respectively.

The following technical result will become handy later on.
\begin{lemma}\label{lem:technical_property_D_normal_cone_map}
	Let $D\subset\mathbb Y$ be nonempty and closed, and fix $\bar y\in D$.
	Then the following assertions hold.
	\begin{enumerate}
		\item For each $u\in\mathbb Y$,
			we have $D\mathcal N_D(\bar y,0)(u)\subset\mathcal N_D(\bar y;u)$.
		\item For each $u\in\mathbb S_{\mathbb Y}$, 
			we have $D_\textup{sub}\mathcal N_D(\bar y,0)(u)\subset\mathcal N_D(\bar y;u)$.
	\end{enumerate}
\end{lemma}
\begin{proof}
	We only prove validity of the first assertion.
	The second one can be shown in analogous fashion.
	
		Fix $u\in\mathbb Y$ and $v\in D\mathcal N_D(\bar y,0)(u)$.
		Then we find sequences $\{u_k\}_{k\in\N},\{v_k\}_{k\in\N}\subset\mathbb Y$
		and $\{t_k\}_{k\in\N}\subset\R_+$ with $u_k\to u$, $v_k\to v$, $t_k\downarrow 0$,
		and $t_kv_k\in\mathcal N_D(\bar y+t_ku_k)$ for each $k\in\N$.
		Since, for each $k\in\N$, $\mathcal N_D(\bar y+t_ku_k)$ is a cone,
		we find $v_k\in\mathcal N_D(\bar y+t_ku_k)$,
		and $v\in\mathcal N_D(\bar y;u)$ follows by robustness of the
		directional limiting normal cone, see \cref{lem:robustness_directional_limiting_normals}.
\end{proof}

In the next two results, we investigate the special situation $\mathbb Y:=\R^m$ in detail.
First, in the case where we consider the normal cone mapping associated 
with polyhedral sets, 
there is no difference
between graphical derivative and graphical subderivative as the subsequent lemma shows.
\begin{lemma}\label{lem:normal_cone_map_of_polyhedral_set}
	Let $D\subset\R^m$ be a polyhedral set.
	Then $\gph\mathcal N_D$ is polyhedral as well, and 
	for arbitrary $(\bar y,\bar y^*)\in\gph\mathcal N_D$ and 
	$u, v\in\R^m \setminus \{0\}$, we have
	\[
		v \in D\mathcal N_D(\bar y,\bar y^*)(u)
		\quad\Longleftrightarrow\quad
		v/\norm{v} \in D_\textup{sub}\mathcal N_D(\bar y,\bar y^*)(u/\norm{u}).
	\]
\end{lemma}
\begin{proof}
	It follows from \cite[Theorem~2]{AdamCervinkaPistek2016} that there exist finitely many
	convex polyhedral sets $D_1,\ldots,D_\ell\subset\R^m$ and closed, convex, polyhedral cones
	$K_1,\ldots,K_\ell\subset\R^m$ such that $\gph\mathcal N_D=\bigcup_{i=1}^\ell D_i\times K_i$.
	Particularly, $\gph\mathcal N_D$ is polyhedral.
	
	Next, consider some nonzero $u, v\in\R^m$ with $v/\norm{v} \in D_\textup{sub}\mathcal N_D(\bar y,\bar y^*)(u/\norm{u})$. 
	Then we find $\{\tilde u_k\}_{k\in\N},\{\tilde v_k\}_{k\in\N}\subset\R^m$ 
	and $\{\tilde t_k\}_{k\in\N},\{\tau_k\}_{k\in\N}\subset\R_+$
	such that $u_k:=\tilde u_k \norm{u}\to u$, $v_k:=\tilde v_k \norm{v} \to v$, 
	$t_k:= \tilde t_k/\norm{u} \downarrow 0$, $\tau_k\downarrow 0$, $\tau_k/t_k\to\infty$, and
	$(\bar y+t_k u_k,\bar y^*+ (\tau_k/\norm{v}) v_k)\in\gph\mathcal N_D$ for all $k\in\N$.
	Thus, we can pick $j\in\{1,\ldots,\ell\}$ and a subsequence (without relabeling) such that
	$(\bar y+t_k u_k,\bar y^*+ (\tau_k/\norm{v}) v_k) \in D_j\times K_j$
	and $\tau_k/\norm{v}>t_k$ for all $k\in\N$.
	By convexity of $K_j$, we also have
	$(\bar y+t_k u_k,\bar y^*+ t_k v_k) \in D_j\times K_j$
	which shows $v\in D\mathcal N_D(\bar y,\bar y^*)(u)$.
	The converse implication can be proven in analogous fashion by multiplying the null sequence in the domain space
	with another null sequence.
\end{proof}

The next lemma shows how the graphical derivative of normal cone mappings associated with Cartesian products
of polyhedral sets can be computed.
\begin{lemma}\label{lem:product_rule_graphical_derivative}
	Fix some $\ell\in\N$. For each $i\in\{1,\ldots,\ell\}$, 
	let $D_i\subset\R^{m_i}$ for some $m_i\in\N$ be polyhedral.
	Set $D:=\prod_{i=1}^\ell D_i$, $m:=\sum_{i=1}^\ell m_i$, and $L:=\{1,\ldots,\ell\}$. 
	Then we have 
	\[
		\gph\mathcal N_D
		=
		\{((y_1,\ldots,y_\ell),(y_1^*,\ldots,y_\ell^*))\in\R^m\times\R^m\,|\,
		\forall i\in L\colon\,(y_i,y_i^*)\in\gph\mathcal N_{D_i}\},
	\]
	and for arbitrary 
	$\bar y:=(\bar y_1,\ldots,\bar y_\ell),\bar y^*:=(\bar y_1^*,\ldots,\bar y_\ell^*)\in\R^m$
	satisfying
	$(\bar y,\bar y^*)\in\gph\mathcal N_D$ as well as $u:=(u_1,\ldots,u_\ell)\in\R^m$,
	we find
	\[
		D\mathcal N_D(\bar y,\bar y^*)(u)
		=
		\{v=(v_1,\ldots,v_\ell)\in\R^m\,|\,
		\forall i\in L\colon\,v_i\in D\mathcal N_{D_i}(\bar y_i,\bar y_i^*)(u_i)\}.
	\]
\end{lemma}
\begin{proof}
	The representation of $\gph\mathcal N_D$ is a simple consequence of the product rule 
	for the computation of limiting
	normals, see e.g.\ \cite[Proposition~1.4]{Mordukhovich2018}, 
	and does not rely on the polyhedrality of the underlying
	sets. 
		Thus, $\gph\mathcal N_D$ is, up to a permutation of components, the same as $\prod_{i=1}^\ell\gph\mathcal N_{D_i}$.
		Since, for each $i\in L$, $\gph\mathcal N_{D_i}$ is polyhedral by \cref{lem:normal_cone_map_of_polyhedral_set},
		the same has to hold for $\gph\mathcal N_D$.
		The final formula of the lemma is a simple consequence of
		\cref{lem:product_rule_tangents_polyhedral_sets}
		and 
		\cite[Exercise~6.7]{RockafellarWets1998}.
\end{proof}

\subsubsection{Coderivatives, pseudo-coderivatives, and super-coderivatives}\label{sec:coderivative}

In the subsequently stated definition, we first recall the notion of regular and limiting coderivative of a set-valued mapping
before introducing its so-called directional pseudo-coderivative. The latter will be of essential importance in the course of the paper.
It corresponds to a minor modification
of the notion of directional pseudo-coderivative introduced by Gfrerer in \cite[Definition~2]{Gfrerer2014a},
which we recall as well.
\begin{definition}\label{def:coderivatives}
	Let $\Phi\colon\mathbb X\tto\mathbb Y$ be a set-valued mapping possessing a closed graph locally around $(\xb,\yb)\in\gph \Phi$.
	Furthermore, let $(u,v)\in\mathbb X\times\mathbb Y$ be a pair of directions.
	\begin{enumerate}
	\item The \emph{regular and limiting coderivative} of $\Phi$ at $(\bar x,\bar y)$
	are the set-valued mappings
	$\widehat D^*\Phi(\bar x,\bar y)\colon\mathbb Y\tto\mathbb X$ 
	and
	$D^*\Phi(\bar x,\bar y)\colon\mathbb Y\tto\mathbb X$ given,
	respectively, by
	\begin{align*}
		\forall y^*\in\mathbb Y\colon\quad
		\widehat D^*\Phi(\bar x,\bar y)(y^*)
		&:=
		\left\{x^*\in\mathbb X\,\middle|\,
			(x^*,-y^*)\in\widehat{\mathcal N}_{\gph\Phi}(\bar x,\bar y)\right\},
		\\
		D^*\Phi(\bar x,\bar y)(y^*)
		&:=
		\left\{x^*\in\mathbb X\,\middle|\,
			(x^*,-y^*)\in\mathcal N_{\gph\Phi}(\bar x,\bar y)\right\}.
	\end{align*}
	The set-valued mapping
	$D^*\Phi((\bar x,\bar y);(u,v))\colon\mathbb Y\tto\mathbb X$ given by
	\begin{align*}
		\forall y^*\in\mathbb Y\colon\quad
		D^*\Phi((\bar x,\bar y);(u,v))(y^*)
		&:=
		\left\{x^*\in\mathbb X\,\middle|\,
			(x^*,-y^*)\in\mathcal N_{\gph\Phi}((\bar x,\bar y);(u,v))\right\}
	\end{align*}
	is the \emph{limiting coderivative} of $\Phi$ at $(\bar x,\bar y)$ \emph{in direction} $(u,v)$.
	If $\Phi$ is single-valued at $\bar x$, we use 
	$\widehat D^*\Phi(\bar x),D^*\Phi(\bar x),D^*\Phi(\bar x;(u,v))\colon\mathbb Y\tto\mathbb X$
	for brevity.
	\item\label{item:new_pseudo_coderivative}
	Given $\gamma \geq 1$ and $u\in\mathbb S_{\mathbb X}$, the
 	\emph{pseudo-coderivative of order $\gamma$} of $\Phi$ at $(\bar x,\bar y)$ 
 	\emph{in direction $(u,v)$}
 	is the mapping $D^\ast_{\gamma} \Phi((\xb,\yb); (u,v))\colon\mathbb Y\tto\mathbb X$
 	which assigns to $y^*\in\mathbb Y$ the set of all $x^*\in\mathbb X$ such that 
 	there are sequences $\{u_k\}_{k\in\N},\{x_k^*\}_{k\in\N}\subset\mathbb X$,
 	$\{v_k\}_{k\in\N},\{y_k^*\}_{k\in\N}\subset\mathbb Y$, 
 	and $\{t_k\}_{k\in\N}\subset\R_+$ which satisfy
 	$u_k\to u$, $v_k\to v$, $t_k\downarrow 0$, $x_k^*\to x^*$, $y_k^*\to y^*$, and
 	\begin{equation}\label{eq:def_con_pseudo_cod}
 		\forall k\in\N\colon\quad
 		\left(x^\ast_k,-\frac{y^\ast_k}{(t_k\norm{u_k})^{\gamma-1}}\right)
		\in 
		\widehat{\mathcal N}_{\gph \Phi}(\xb + t_k u_k,\yb + (t_k\norm{u_k})^{\gamma} v_k).
 	\end{equation}
 	In the case $\gamma:=1$, this definition recovers the one of $D^*\Phi((\bar x,\bar y);(u,v))$.
 	\item\label{item:Gfrerer_pseudo_coderivative} 
 	Given $\gamma\geq 1$ and $u\in\mathbb S_{\mathbb X}$, 
 	\emph{Gfrerer's pseudo-coderivative of order $\gamma$} of $\Phi$ 
 	at $(\bar x,\bar y)$ \emph{in direction $(u,v)$}
 	is the mapping 
 	$\widetilde D^\ast_{\gamma} \Phi((\xb,\yb); (u,v))\colon\mathbb Y\tto\mathbb X$
 	which assigns to $y^*\in\mathbb Y$ the set of all $x^*\in\mathbb X$ 
 	such that there are sequences $\{u_k\}_{k\in\N},\{x_k^*\}_{k\in\N}\subset\mathbb X$,
 	$\{v_k\}_{k\in\N},\{y_k^*\}_{k\in\N}\subset\mathbb Y$, 
 	and $\{t_k\}_{k\in\N}\subset\R_+$ which satisfy
 	$u_k\to u$, $v_k\to v$, $t_k\downarrow 0$, $x_k^*\to x^*$, $y_k^*\to y^*$, and
 	\begin{equation}\label{eq:def_con_Gfrerer_pseudo_cod}
 		\forall k\in\N\colon\quad
 		\left(x^\ast_k,-\frac{y^\ast_k}{(t_k\norm{u_k})^{\gamma-1}}\right)
		\in 
		\widehat{\mathcal N}_{\gph \Phi}(\xb + t_k u_k,\yb + t_k v_k).
 	\end{equation}
 	Again, for $\gamma:=1$, we recover the definition of $D^*\Phi((\bar x,\bar y);(u,v))$.
	\end{enumerate}
\end{definition}

Let $\Phi\colon\mathbb X\tto\mathbb Y$ be a set-valued mapping whose graph is closed locally around $(\bar x,\bar y)\in\gph\Phi$
and fix a pair of directions $(u,v)\in\mathbb S_{\mathbb X}\times\mathbb Y$, $(x^*,y^*) \in\mathbb X\times\mathbb Y$, and $\gamma>1$.
Then we obtain the trivial relations
\begin{equation}\label{eq:trivial_upper_estimate_pseudo_coderivative}
	x^* \in D^*_\gamma\Phi((\bar x,\bar y);(u,v))(y^*)
	\ \Longrightarrow \
	\left\{\begin{aligned}
			&0 \in D\Phi(\bar x,\bar y)(u), \ 0 \in D^*\Phi(\xb,\yb)(y^*), \\
			&0 \in D^*\Phi((\xb,\yb);(u,0))(y^*),\\
			&v \in D_\gamma\Phi(\bar x,\bar y)(u),\\
			&x^* \in \widetilde D^*_\gamma\Phi((\bar x,\bar y);(u,0))(y^*).
		\end{aligned}
	\right.
\end{equation}
Note also that the mappings $D^*\Phi((\bar x,\bar y);(u,v))$ and $\widetilde D^*_\gamma\Phi((\bar x,\bar y);(u,v))$
have a nonempty graph if and only if $v \in D\Phi(\bar x,\bar y)(u)$ while
the mapping $D^*_\gamma\Phi((\bar x,\bar y);(u,v))$
has a nonempty graph if and only if $v \in D_\gamma\Phi(\bar x,\bar y)(u)$.%

	Since the (directional) limiting coderivative
	is defined via the (directional) limiting normal cone,
	it possesses a robust behavior as well.
	In the subsequent lemma, we show
	a somewhat robust behavior of the directional
	pseudo-coderivatives under consideration,
	which will be important later on.
	Basically, we prove that one can replace the regular by the
	limiting normal cone in \eqref{eq:def_con_pseudo_cod} and
	\eqref{eq:def_con_Gfrerer_pseudo_cod} without changing the
	resulting pseudo-coderivative.
	The technical proof, which is based on a standard diagonal sequence argument,
	is presented in \cref{sec:appendix} for the purpose of completeness.

	\begin{lemma}\label{lem:pseudo_coderivative_via_limiting_normals}
		\Cref{def:coderivatives}\,\ref{item:new_pseudo_coderivative} and
		\cref{def:coderivatives}\,\ref{item:Gfrerer_pseudo_coderivative} can equivalently be formulated in terms of
		limiting normals.
	\end{lemma}

To illustrate the pseudo-coderivatives from \cref{def:coderivatives},
we revisit \cref{ex:metric_pseudo_regularity}.
\begin{example}\label{ex:pseudo_coderivatives}
For $\gamma > 1$, we consider the mapping $\Phi \colon \R \tto \R$, 
given by $\Phi(x) := \{ \vert x \vert^\gamma \}$, $x\in\R$,
already discussed in \cref{ex:metric_pseudo_regularity}.
Set $(\bar x, \bar y) := (0,0)$ as well as $u:=\pm 1$ and choose $v\in\R$ arbitrarily.
First, $v \in D \Phi(\bar x,\bar y)(u)$ by definition
requires sequences $\{t_k\}_{k\in\N}\subset\R_+$ and
$\{u_k\}_{k\in\N},\{v_k\}_{k\in\N}\subset\R$ satisfying
$t_k \downarrow 0$, $u_k \to u$, $v_k \to v$, 
and $t_k v_k = (t_k \vert u_k \vert)^\gamma$ for all $k\in\N$, showing $v=0$.
Thus, we fix $v:=0$ to find
$D^*\Phi((\bar x,\bar y);(u,0))(y^*) = \{0\}$ for all $y^*\in\R$
as the defining sequences $\{x_k^*\}_{k\in\N},\{y_k^*\}_{k\in\N}\subset\R$
satisfy $x_k^* \to x^*$, $y_k^* \to y^*$, and
$x_k^* = \gamma (t_k \vert u_k \vert)^{\gamma-1} \sgn (u_k) y_k^*$
for all $k\in\N$.
Furthermore,
$\widetilde D^\ast_{\gamma} \Phi((\xb,\yb); (u,0))(y^*) = \{\gamma \sgn (u) y^*\}$
holds for each $y^*\in\R$
as the defining sequences $\{x_k^*\}_{k\in\N},\{y_k^*\}_{k\in\N}\subset\R$
satisfy $x_k^* \to x^*$, $y_k^* \to y^*$, and
$x_k^* = \gamma \sgn (u_k) y_k^*$ for all $k\in\N$.
Using similar arguments as above, one can check that
$v \in D_\gamma\Phi(\bar x,\bar y)(u)$ yields $v=1$,
and for $v:=1$, we get $D^\ast_\gamma\Phi((\bar x,\bar y);(u,1))(y^*)=\{\gamma\sgn(u)y^*\}$
for all $y^*\in\R$.
\end{example}
	
Below, we introduce yet another concept of coderivative which will become important in \cref{sec:asymptotic_regularity_via_super_coderivative}.

\begin{definition}\label{def:super_coderivative}
	Let $\Phi\colon\mathbb X\tto\mathbb Y$ be a set-valued mapping with a closed graph
	and fix $(\bar x,\bar y)\in\gph\Phi$ and 
	$(u,v)\in\mathbb S_{\mathbb X}\times\mathbb S_{\mathbb Y}$.
	The \emph{super-coderivative} of $\Phi$ at $(\bar x,\bar y)$ in direction $(u,v)$ is
	the mapping $D^*_\textup{sup}\Phi((\xb,\yb); (u,v))\colon\mathbb Y \tto \mathbb X$,
	which assigns to $y^* \in\mathbb Y$ the set of all $x^\ast \in \mathbb X$ 
	such that there are sequences
	$\{u_k\}_{k\in\N},\{x_k^*\}_{k\in\N}\subset\mathbb X$, 
	$\{v_k\}_{k\in\N},\{y_k^*\}_{k\in\N}\subset\mathbb Y$,
	and $\{t_k\}_{k\in\N},\{\tau_k\}_{k\in\N}\subset\R_+$ which satisfy
	$u_k\to u$, $v_k\to v$, $x_k^*\to x^*$, $y_k^*\to y^*$, $t_k\downarrow 0$, 
	$\tau_k\downarrow 0$, and $\tau_k/t_k\to 0$
	such that
	\begin{equation}\label{eq:characterization_super_coderivative}
		x_k^* 
		\in 
		\widehat{D}^*\Phi(\xb + t_k u_k,\yb + \tau_k v_k)
			(((t_k \norm{u_k})/(\tau_k \norm{v_k})) y_k^*)
	\end{equation}
	holds for all $k\in\N$.
\end{definition}

We start with some remarks regarding \cref{def:super_coderivative}.
First, observe that we only exploit the super-coderivative w.r.t.\
unit directions $(u,v)\in\mathbb S_{\mathbb X}\times\mathbb S_{\mathbb Y}$
which also means that $\{u_k\}_{k\in\N}\subset\mathbb X$ and
$\{v_k\}_{k\in\N}\subset\mathbb Y$ can be chosen such that $u_k\neq 0$
and $v_k\neq 0$ hold for all $k\in\N$. 
Particularly, condition \eqref{eq:characterization_super_coderivative} is reasonable.

Second, we would like to note that 
$x^*\in D^*_\textup{sup}\Phi((\bar x,\bar y);(u,v))(y^*)$ implies the existence of
sequences $\{u_k\}_{k\in\N}\subset\mathbb X$, 
$\{v_k\}_{k\in\N}\subset\mathbb Y$,
and $\{t_k\}_{k\in\N},\{\tau_k\}_{k\in\N}\subset\R_+$ which satisfy
$u_k\to u$, $v_k\to v$, $t_k\downarrow 0$, $\tau_k\downarrow 0$, and $\tau_k/t_k\to 0$
as well as $(\bar x+t_ku_k,\bar y+\tau_kv_k)\in\gph\Phi$ for all $k\in\N$.
Thus, in the light of \cref{def:graphical_derivative}\,\ref{item:def_graph_subderivative} 
of the graphical subderivative, one might be tempted
to say that the pair $(u,v)$ belongs to the graph of the graphical
\emph{super-derivative} of $\Phi$ at $(\bar x,\bar y)$.
This justifies the terminology in \cref{def:super_coderivative}.

Let us briefly discuss the relation between pseudo-coderivatives and the
novel super-coderivative from \cref{def:super_coderivative}.
Consider $\gamma > 1$ and $x^* \in D^*_\gamma\Phi((\xb,\yb); (u,v))(y^*)$ 
for $(u,v)\in\mathbb S_{\mathbb X}\times\mathbb S_{\mathbb Y}$ and $y^*\in\mathbb Y^*$.
Setting $\tau_k := (t_k \norm{u_k})^{\gamma}$ for each $k\in\N$, 
where $\{t_k\}_{k\in\N}\subset\R_+$ 
and $\{u_k\}_{k\in\N}\subset\mathbb X$ are the sequences from the definition of the 
pseudo-coderivative,
we get $x^* \in D^*_\textup{sup}\Phi((\xb,\yb); (u,v))(y^*)$ 
since $t_k^{\gamma-1}\norm{u_k}^\gamma\to 0$.

In the subsequent lemma, we comment on the converse inclusion which, 
to some extent, 
holds in the presence of a qualification condition in terms of the pseudo-coderivative.
\begin{lemma}\label{lem:super_coderivative_vs_pseudo_coderivative}
	Let $(\bar x,\bar y)\in\gph\Phi$, $(u,v)\in\mathbb S_{\mathbb X}\times\mathbb S_{\mathbb Y}$,
	$y^* \in \mathbb Y$, and $\gamma > 1$ be fixed.
	Furthermore, assume that $\ker D^*_\gamma\Phi((\bar x,\bar y);(u,0))\subset\{0\}$ holds.
	Then there exists $\alpha > 0$ such that
	\begin{align*}
		D^*_\textup{sup}\Phi((\bar x,\bar y);(u,v))(y^*)
		&\subset
		\widetilde{D}^*_\gamma\Phi((\bar x,\bar y);(u,0))(0)
		\cup
		D^*_\gamma\Phi((\bar x,\bar y);(u,\alpha v))(y^*/\alpha)
		\\
		&\qquad
		\cup
		\Im D^*_\gamma\Phi((\bar x,\bar y);(u,0))
		\\
		&\subset
		\Im \widetilde D^*_\gamma\Phi((\bar x,\bar y);(u,0)).
	\end{align*}
\end{lemma}
\begin{proof}
	Let $x^*\in D^*_\textup{sup}\Phi((\bar x,\bar y);(u,v))(y^*)$ be arbitrarily chosen.
	Then we find sequences $\{u_k\}_{k\in\N},\{x_k^*\}_{k\in\N}\subset\mathbb X$, 
	$\{v_k\}_{k\in\N},\{y_k^*\}_{k\in\N}\subset\mathbb Y$,
	and $\{t_k\}_{k\in\N},\{\tau_k\}_{k\in\N}\subset\R_+$ which satisfy
	$u_k\to u$, $v_k\to v$, $x_k^*\to x^*$, $y_k^*\to y^*$, $t_k\downarrow 0$, 
	$\tau_k\downarrow 0$, and $\tau_k/t_k\to 0$
	as well as \eqref{eq:characterization_super_coderivative} for all $k\in\N$.
	This also gives us
	\begin{equation}\label{eq:insert_order_gamma_into_super_coderivative}
		x_k^*
		\in
		\widehat D^*\Phi
			\left(
				\bar x+t_ku_k,
				\bar y+(t_k\norm{u_k})^\gamma\frac{\tau_kv_k}{(t_k\norm{u_k})^\gamma}
			\right)
			\left(
				(t_k\norm{u_k})^{1-\gamma}
				\frac{(t_k\norm{u_k})^\gamma}{\tau_k\norm{v_k}}y_k^*
			\right)
	\end{equation}
	for all $k\in\N$.
	Set $\tilde y_k^*:=(t_k\norm{u_k})^\gamma/(\tau_k\norm{v_k})y_k^*$ for each $k\in\N$. 
	In the case where $\{\tilde y_k^*\}_{k\in\N}$ is not bounded, 
	we have $(\tau_k\norm{v_k})/(t_k\norm{u_k})^\gamma\to 0$
	along a subsequence (without relabeling), and taking the limit in
	\[
		x_k^*/\nnorm{\tilde y_k^*}
		\in
		\widehat D^*\Phi
		\left(
			\bar x+t_ku_k,
			\bar y+(t_k\norm{u_k})^\gamma\frac{\tau_kv_k}{(t_k\norm{u_k})^\gamma}
		\right)
		\left(
			(t_k\norm{u_k})^{1-\gamma} 
			\tilde y_k^*/\nnorm{\tilde y_k^*}
		\right)
	\]
	yields 
	that $\ker D^*_\gamma\Phi((\bar x,\bar y);(u,0))$ contains a nonzero element, 
	which is a contradiction.
	Hence, $\{\tilde y_k^*\}_{k\in\N}$ is bounded.
	
	For each $k\in\N$, we set $\alpha_k:=\tau_k\norm{v_k}/(t_k\norm{u_k})^\gamma$.
	First, suppose that $\{\alpha_k\}_{k\in\N}$ is not bounded.
	Then, along a subsequence (without relabeling), we may assume $\alpha_k\to\infty$.
	By boundedness of $\{y_k^*\}_{k\in\N}$, $\tilde y_k^*\to 0$ follows.
	Rewriting \eqref{eq:insert_order_gamma_into_super_coderivative} yields
	\[
		x_k^*\in\widehat D^*\Phi\left(\bar x+t_ku_k,\bar y+t_k\frac{\tau_kv_k}{t_k}\right)
		\left((t_k\norm{u_k})^{1-\gamma}\tilde y_k^*\right)
	\]
	for each $k\in\N$, and taking the limit $k\to\infty$ 
	while respecting $\tau_k/t_k\to 0$, thus, gives 
	$x^*\in \widetilde D^*_\gamma\Phi((\bar x,\bar y);(u,0))(0)$.
	In the case where $\{\alpha_k\}_{k\in\N}$ converges to some $\alpha>0$ (along a subsequence
	without relabeling), we can simply take the limit $k\to\infty$ in
	\eqref{eq:insert_order_gamma_into_super_coderivative} in order to find
	$x^*\in D^*_\gamma\Phi((\bar x,\bar y);(u,\alpha v))(y^*/\alpha)$.
	Finally, let us consider the case $\alpha_k\to 0$ (along a subsequence without
	relabeling). Then, by boundedness of $\{\tilde y_k^*\}_{k\in\N}$,
	taking the limit $k\to\infty$ in \eqref{eq:insert_order_gamma_into_super_coderivative}
	gives $x^*\in \Im D^*_\gamma\Phi((\bar x,\bar y);(u,0))$.
	Thus, we have shown the first inclusion.
	
	The second inclusion follows by the upper estimate 
	\eqref{eq:trivial_upper_estimate_pseudo_coderivative} for the pseudo-coderivative.
\end{proof}

\subsubsection{Sufficient conditions for pseudo-(sub)regularity}\label{sec:sufficient_conditions_pseudo_regularity}

Graphical derivative and (directional) limiting coderivative 
are powerful tools for studying regularity properties of set-valued mappings,
such as (strong) metric regularity and subregularity, as well as their 
inverse counterparts of Lipschitzness,
such as Aubin property and (isolated) calmness.
Indeed, given a closed-graph set-valued mapping $\Phi\colon\mathbb X\tto\mathbb Y$, 
metric regularity and strong metric subregularity at some point $(\xb,\yb) \in \gph \Phi$
are characterized, respectively, by
\begin{subequations}\label{eq:char_MS_via_generalized_derivatives}
	\begin{align}
		\label{eq:Mordukhovich_criterion}
				\ker D^*\Phi(\xb,\yb) &= \{0\},\\
		\label{eq:Levy_Rockafellar_criterion}
				\ker D\Phi(\xb,\yb) &= \{0\},
	\end{align}
\end{subequations}
see e.g.\ \cite{Levy96,Mordukhovich2018,RockafellarWets1998} for the definition of
these Lipschitzian properties as well as the above results.
	Let us mention that \eqref{eq:Mordukhovich_criterion} is referred to as
	\emph{Mordukhovich criterion} in the literature,
	while \eqref{eq:Levy_Rockafellar_criterion} is called
	\emph{Levy--Rockafellar criterion}.

	For fixed $u\in\mathbb S_{\mathbb X}$, we will refer to
	\begin{equation}\label{eq:FOSCMS_u}
		\ker D^*\Phi((\bar x,\bar y);(u,0)) 
		\subset
		\{0\},
	\end{equation}
	which implies that $\Phi$ is metrically subregular at $(\bar x,\bar y)$ in direction $u$,
	see e.g.\ \cite[Theorem~5]{Gfrerer2013}, as FOSCMS$(u)$.
		Note that it is formulated as an inclusion as
		the left-hand side in \eqref{eq:FOSCMS_u} is empty whenever
		$u\notin\ker D\Phi(\bar x,\bar y)$.
		Indeed, in this case, $\Phi$ is trivially metrically subregular at $(\bar x,\bar y)$
		in direction $u$.
	Furthermore, whenever \eqref{eq:FOSCMS_u} holds for all $u\in \ker D\Phi(\xb,\yb)\cap\mathbb S_{\mathbb X}$,
	which we will refer to as FOSCMS, then $\Phi$ is already metrically subregular at $(\bar x,\bar y)$,
	see \cite[Lemma~2.7]{Gfrerer2014}.
	Above, FOSCMS abbreviates
	\emph{First-Order Sufficient Condition for Metric Subregularity},
	and this terminology has been coined in \cite{Gfrerer2013}.
	Clearly, each of the conditions from \eqref{eq:char_MS_via_generalized_derivatives}
	is sufficient for FOSCMS.
The relations \eqref{eq:trivial_upper_estimate_pseudo_coderivative} suggest that
the pseudo-coderivative can be useful particularly in situations where the above regularity properties,
which are related to (first-order) coderivatives, fail.

Note that the aforementioned notions of regularity and Lipschitzness 
express certain linear rate of change of the mapping.
Similarly, there is an underlying linearity in the definition 
of graphical derivative and coderivatives.
Take the graphical derivative for instance. 
Since the same sequence $\{t_k\}_{k\in\N}$ appears in the domain 
as well as in the range space, if $v \in D\Phi(\xb,\yb)(u)$ 
implies that $u\in\mathbb X$ and $v\in\mathbb Y$ are both nonzero,
it suggests a proportional (linear) rate of change.
Thus, in order to characterize pseudo-(sub)regularity of order $\gamma>1$ of $\Phi$, 
it is not very surprising that 
we need to exploit derivative-like objects based on sub- or superlinear structure.
Exemplary, this has been successfully visualized in \cite[Corollary~2]{Gfrerer2014a} 
by means of Gfrerer's directional pseudo-coderivative
of order $\gamma>1$ from \cref{def:coderivatives}\,\ref{item:Gfrerer_pseudo_coderivative}.
Here, we show 
that the fundamental result from \cite[Theorem~1(2)]{Gfrerer2014a} yields also
an analogous sufficient condition for metric pseudo-subregularity 
via the pseudo-coderivative from \cref{def:coderivatives}\,\ref{item:new_pseudo_coderivative}.

\begin{lemma}\label{lem:sufficient_condition_pseudo_subregularity_abstract}
	Let $\Phi\colon\mathbb X\tto\mathbb Y$ be a set-valued mapping 
	having a closed graph locally around
	$(\bar x,\bar y)\in\gph\Phi$, 
	fix a direction $u\in\mathbb S_{\mathbb X}$, and some $\gamma\geq 1$.
	Assume that
	\begin{equation}\label{eq:FOSCMS_gamma}
		\ker D^\ast_{\gamma}  \Phi((\xb,\yb ); (u,0)) \subset \{0\}
	\end{equation}
	holds.
	Then $\Phi$ is metrically pseudo-subregular of order $\gamma$ 
	at $(\bar x,\bar y)$ in direction $u$.
\end{lemma}
\begin{proof}
	Suppose that $\Phi$ is not metrically pseudo-subregular 
	of order $\gamma$ at $(\bar x,\bar y)$
	in direction $u$. 
	Due to \cite[Theorem~1(2)]{Gfrerer2014a}, we find sequences $\{t_k\}_{k\in\N}\subset\R_+$,
	$\{u_k\}_{k\in\N},\{x_k^*\}_{k\in\N}\subset\mathbb X$, 
	and $\{v_k\}_{k\in\N},\{y_k^*\}_{k\in\N}\subset\mathbb Y$
	satisfying (among other things) $t_k\downarrow 0$, $u_k\to u$, $t_k^{1-\gamma}v_k\to 0$, as well as
	$x_k^*\to 0$, such that
	$\norm{y_k^*}=1$ and 
	\[
		(x_k^*,-y_k^*/(t_k\norm{u_k})^{\gamma-1})
		\in
		\widehat{\mathcal N}_{\gph\Phi}((\bar x,\bar y)+t_k(u_k,v_k))
	\]
	for each $k\in\N$. Let us set $\tilde v_k:=t_k^{1-\gamma}\norm{u_k}^{-\gamma} v_k$ 
	for each $k\in\N$.
	Then we have 
	\[
		(x_k^*,-y_k^*/(t_k\norm{u_k})^{\gamma-1})
		\in
		\widehat{\mathcal N}_{\gph\Phi}(\bar x+t_ku_k,\bar y+(t_k\norm{u_k})^\gamma \tilde v_k)
	\]
	for each $k\in\N$ and $\tilde v_k\to 0$ from $t_k^{1-\gamma}v_k\to 0$.
	Observing that $\{y_k^*\}_{k\in\N}$ possesses 
	a nonvanishing accumulation point $y^*\in\mathbb Y$,
	taking the limit along a suitable subsequence 
	yields $0\in  D^*_\gamma\Phi((\bar x,\bar y);(u,0))(y^*)$
	which contradicts the assumptions of the lemma.
\end{proof}

Let us remark that due to \eqref{eq:trivial_upper_estimate_pseudo_coderivative}, condition
\begin{equation}\label{eq:FO_characterization_of_dir_metric_pseudo_reg}
	\ker \widetilde{D}^*_\gamma\Phi((\bar x,\bar y);(u,0))
	\subset
	\{0\}
\end{equation}
is stronger than \eqref{eq:FOSCMS_gamma} and, thus, also sufficient 
for metric pseudo-subregularity
of $\Phi$ of order $\gamma\geq 1$ at $(\bar x,\bar y)$ in direction $u$.
By means of \cite[Corollary~2]{Gfrerer2014a},
\eqref{eq:FO_characterization_of_dir_metric_pseudo_reg}
is actually equivalent to $\Phi$ being
metrically pseudo-regular at $(\bar x,\bar y)$ in
direction $(u,0)$.
Note that in the case $\gamma:=1$, both conditions \eqref{eq:FOSCMS_gamma} 
and \eqref{eq:FO_characterization_of_dir_metric_pseudo_reg}
recover 
FOSCMS$(u)$.
	In \cref{ex:pseudo_coderivatives}, \eqref{eq:FOSCMS_gamma} 
	and \eqref{eq:FO_characterization_of_dir_metric_pseudo_reg} hold simultaneously.
	The following example illustrates that
	\eqref{eq:FOSCMS_gamma} can be strictly milder than 
	\eqref{eq:FO_characterization_of_dir_metric_pseudo_reg}.

\begin{example}
For $\gamma > 1$, we consider the mapping $\Phi \colon \R \tto \R$ given by
\[
	\gph\Phi
	:=
	\{(x,y)\,|\,|x|^\gamma\leq y \leq 2|x|^\gamma\}
	\cap
	\left(\bigcup\nolimits_{k\in\N}\R\times\{1/2^k\}\right).
\]
Essentially, $\gph\Phi$ is a closed staircase enclosed by the graphs of 
the functions $x\mapsto|x|^\gamma$ and
$x\mapsto 2|x|^\gamma$.
Set $(\bar x, \bar y) := (0,0)$ and $u:=1$.
First, it is easy to see that \eqref{eq:FOSCMS_gamma} is satisfied, because one can show
$\ker D^*_\gamma\Phi((\bar x,\bar y);(u,0)) = \emptyset$.
Indeed, the sequences 
$\{t_k\}_{k\in\N}\subset\R_+$ and $\{u_k\}_{k\in\N},\{v_k\}_{k\in\N}\subset\R$ 
from the definition of the pseudo-coderivative satisfy,
among others, $t_k \downarrow 0$, $u_k \to u$, $v_k \to v$, 
and $(t_k |u_k|)^\gamma \leq (t_k |u_k|)^\gamma v_k \leq 2 (t_k |u_k|)^\gamma$ for each $k\in\N$.
Thus, $D^*_\gamma\Phi((\bar x,\bar y);(u,v))$ can have a nonempty graph 
only for $v\in[1,2]$.
Next, let us argue that \eqref{eq:FO_characterization_of_dir_metric_pseudo_reg} 
fails due to
$1 \in \ker \widetilde{D}^*_\gamma\Phi((\bar x,\bar y);(u,0))$.
We consider the sequences $\{t_k\}_{k\in\N}\subset\R_+$ and
$\{u_k\}_{k\in\N},\{v_k\}_{k\in\N}\subset\R$ given by
\[
 	\forall k\in\N\colon\quad
 	t_k := \left(\frac{3}{2^{k+2}}\right)^{1/\gamma}, 
 	\qquad
 	u_k := 1,
 	\qquad
	v_k := \frac1{2^{k}} \left(\frac{3}{2^{k+2}}\right)^{-1/\gamma}.
\]
We obviously have $t_k\downarrow 0$, $u_k\to 1$, as well as $v_k\to 0$, and 
one can easily check that $(t_ku_k,t_kv_k)\in\gph\Phi$ holds for all $k\in\N$.
By construction, there exist vertical normals to $\gph\Phi$ at 
$(t_k u_k, t_k v_k)$ for each $k\in\N$,
so we can choose $x_k^* := 0$ and $y_k^* := 1$ satisfying 
\eqref{eq:def_con_Gfrerer_pseudo_cod}.
Taking the limit $k\to\infty$ shows
$1 \in \ker \widetilde{D}^*_\gamma\Phi((\bar x,\bar y);(u,0))$.
\end{example}

\begin{remark}\label{Rem:Strong_pseudo_subreg}
	Let $\Phi\colon\mathbb X\tto\mathbb Y$ be a set-valued mapping 
	having locally closed graph around
	$(\bar x,\bar y)\in\gph\Phi$, 
	and fix some $\gamma\geq 1$.
	Note that if we replace the set $\Phi^{-1}(\bar y)$
	by just the singleton $\{\bar x\}$ in \cref{def:metric_pseudo_subregularity} 
	of metric pseudo-subregularity,
	the estimate \eqref{eq:estimate_pseudo_regularity}
	simplifies to
	\[
		\norm{x-\bar x}^{\gamma}
		\leq
		\kappa\,\dist(\bar y,\Phi(x)).
	\]
	Asking this to hold for all $x\in\mathbb B_\varepsilon(\bar x)$ and some $\varepsilon>0$ 
	seems like a natural way to define \emph{strong metric pseudo-subregularity}
	of order $\gamma$ of $\Phi$ at $(\bar x,\bar y)$.
	It is an easy exercise to verify that this condition is satisfied
	if and only if $\ker D_{\gamma}\Phi(\bar x,\bar y)=\{0\}$.
	This characterization is clearly an extension of the
	Levy--Rockafellar criterion \eqref{eq:Levy_Rockafellar_criterion},
	and it provides a justification for the graphical pseudo-derivative.
\end{remark}

Finally, by definition of the pseudo-coderivatives, we easily find the inclusions
\begin{align*}
	\ker D^*_{\gamma+\varepsilon}\Phi((\bar x,\bar y);(u,0))
	&\subset
	\ker D^*_\gamma\Phi((\bar x,\bar y);(u,0)),\\
	\ker \widetilde D^*_{\gamma+\varepsilon}\Phi((\bar x,\bar y);(u,0))
	&\subset
	\ker \widetilde D^*_\gamma\Phi((\bar x,\bar y);(u,0))
\end{align*}
for each $\gamma\geq 1$ and $\varepsilon>0$.
Hence, as $\gamma$ increases, the qualification conditions
\eqref{eq:FOSCMS_gamma} and \eqref{eq:FO_characterization_of_dir_metric_pseudo_reg} 
become weaker.

\section{Pseudo-(sub)regularity of order 2 for constraint mappings}\label{sec:variational_analysis_constraint_mapping}

In this section, we address the pseudo-coderivative calculus for so-called constraint mappings 
$\Phi\colon\mathbb X\tto\mathbb Y$ which are given by $\Phi(x):=g(x)-D$ for all $x\in\mathbb X$,
where $g\colon\mathbb X\to\mathbb Y$ is a single-valued continuous function and $D\subset\mathbb Y$ 
is a closed set, 
and apply our findings from \cref{sec:sufficient_conditions_pseudo_regularity}
in order to derive sufficient conditions for directional metric pseudo-(sub)regularity of order 2.
Let us emphasize that this representation of $\Phi$ will be a standing assumption
in the overall section.
The constraint mapping $\Phi$ plays an important role for the analysis 
of so-called geometric constraint systems of type $g(x)\in D$.

\subsection{Directional pseudo-coderivatives and sufficient conditions}
\label{sec:pseudo_coderivatives_constraint_systems_theory}

The first lemma of this subsection addresses 
upper estimates of the regular, limiting, and
directional limiting coderivative of constraint
mappings.
These results are in principle quite standard,
with the exception of the lower estimates in~\ref{item:constraint_maps_regular_coderivative}
and~\ref{item:constraint_maps_directional_limiting_coderivative},
which can be shown using \cite[Theorem~3.1]{BenkoMehlitz2020} 
and \cite[Lemma 6.1]{BenkoGfrererOutrata2019}, respectively.
However, since we proceed in a fairly mild setting
where $g$ is assumed to be merely continuous,
we cannot simply rely on change-or-coordinates
formulas, see e.g.\ \cite[Exercise~6.7]{RockafellarWets1998},
even for the proof of the standard parts
in~\ref{item:constraint_maps_regular_coderivative} 
and~\ref{item:constraint_maps_limiting_coderivative}.
Thus, we prove everything using the results from our recent paper \cite{BenkoMehlitz2020}.

\begin{lemma}\label{lem:coderivatives_constraint_maps}
	Fix $(x,y)\in\gph\Phi$. Then the following statements hold.
	\begin{enumerate}
		\item\label{item:constraint_maps_regular_coderivative} 
		For each $y^*\in\mathbb Y$, we have
        \[
		\widehat{D}^*\Phi(x,y)(y^*)
		\subset
		\begin{cases}
			\widehat{D}^*g(x) (y^*)	&	y^*\in \widehat{\mathcal N}_D(g(x)-y),\\
			\emptyset				&	\text{otherwise,}
		\end{cases}
		\]
        and the opposite inclusion holds if $g$ is calm at $x$.
		\item\label{item:constraint_maps_limiting_coderivative}
		 For each $y^*\in\mathbb Y$, we have
		\[
		D^*\Phi(x,y)(y^*)
		\subset
		\begin{cases}
			D^*g(x) (y^*)	&	y^*\in \mathcal N_D(g(x)-y),\\
			\emptyset	&	\text{otherwise},
		\end{cases}
		\]
		and the opposite inclusion holds whenever $g$ is continuously
		differentiable at $x$.
		\item\label{item:constraint_maps_directional_limiting_coderivative} 
		 For each pair of directions $(u,v)\in\mathbb X\times\mathbb Y$
			and each $y^*\in\mathbb Y$, we have
			\[
				D^*\Phi((x,y);(u,v))(y^*)
				\subset
				\begin{cases}
					\bigcup\limits_{w\in Dg(x)(u)}
					D^*g(x;(u,w))(y^*)	& y^*\in\mathcal N_D(g(x)-y;w-v),\\
					\emptyset				&\text{otherwise}
				\end{cases}
			\]
			provided $g$ is calm at $x$,
			and the opposite inclusion holds whenever $g$ is continuously
			differentiable at $x$.
	\end{enumerate}
\end{lemma}
\begin{proof}
	\begin{enumerate}
		\item For the proof, we observe that
			$\gph\Phi=\gph g+(\{0\}\times (-D))$
			is valid.
			Now, we exploit the sum rule from \cite{BenkoMehlitz2020}.
			Therefore, let us introduce the surrogate mapping 
			$M\colon\mathbb X\times\mathbb Y\tto(\mathbb X\times\mathbb Y)\times(\mathbb X\times\mathbb Y)$ given by
			\begin{equation}\label{eq:intermediate_map_sum_rule}
				\begin{aligned}
				M(x,y)
				:=&
				\left\{
					((x_1,y_1),(x_2,y_2))\in\gph g\times(\{0\}\times(-D))\,\middle|\,
						\begin{aligned}
							x&=x_1+x_2\\
							y&=y_1+y_2
						\end{aligned}
				\right\} 
				\\
				=&
				\begin{cases}
					\{((x,g(x)),(0,y-g(x)))\}	& g(x)-y\in D,\\
					\emptyset			& \text{otherwise}
				\end{cases}
				\end{aligned}
			\end{equation}
			for all $(x,y)\in\mathbb X\times\mathbb Y$,
			and observe that $\gph\Phi=\dom M$ holds while $M$ is single-valued and continuous on $\gph\Phi$.
			Now, we find
			\[
				\widehat{\mathcal N}_{\gph\Phi}(x,y)
				\subset
				\widehat{D}^*M((x,y),((x,g(x)),(0,y-g(x))))((0,0),(0,0))
			\]
			for all $(x,y)\in\gph\Phi$
			from \cite[Theorem~3.1]{BenkoMehlitz2020}, and the converse
			inclusion holds if $g$ is calm at $x$ since this ensures that $M$ is
			so-called isolatedly calm at the point of interest, 
			see \cite[Corollary~4.4, Section~5.1.1]{BenkoMehlitz2020}.
			Now, computing the regular normal cone to $\gph M$ via 
			\cite[Lemmas~2.1, 2.2]{BenkoMehlitz2020} 
			and applying the definition of the regular coderivative yields the claim.
		\item
			The proof of the inclusion $\subset$ is similar as the one of the first statement.
			Again, we exploit the mapping $M$ given in \eqref{eq:intermediate_map_sum_rule} and
			apply \cite[Theorem~3.1]{BenkoMehlitz2020} while observing that $M$ is 
			so-called inner semicompact w.r.t.\ its domain 
			at each point $(x,y)\in\gph\Phi$ by continuity of $g$.
			In the presence of continuous differentiability, the converse inclusion
			$\supset$ follows easily by applying the change-of-coordinates formula 
			provided in \cite[Exercise~6.7]{RockafellarWets1998}. 
		\item 
			This assertion can be shown in similar way as the second one,
			see \cite[Lemma~2.1]{BenkoMehlitz2020} as well.
	\end{enumerate}
\end{proof}

Let us note that the upper estimate in~\ref{item:constraint_maps_regular_coderivative} was also shown in 
\cite[Lemma~3.2]{BaiYeZhang2019}, but it actually follows directly from
\cite[Exercise~6.44]{RockafellarWets1998} upon realizing 
$\gph\Phi=\gph g+(\{0\}\times (-D))$.
In the case where $g$ is not calm at the reference point, one can still obtain an upper
estimate for the directional limiting coderivative from \cite[Theorem~3.1]{BenkoMehlitz2020}
which is slightly more technical since it comprises another union over
$w\in Dg(x)(0)\cap\mathbb S_{\mathbb Y}$.

Next, we estimate the directional pseudo-coderivatives of order $2$ 
of constraint mappings in terms of initial problem data.

\begin{theorem}\label{The : NCgen}
Let $g$ be twice continuously differentiable.
Given $(\xb,0) \in \gph \Phi$ and a direction $u \in \mathbb S_{\mathbb X}$, let
\[
	x^* \in \widetilde{D}^\ast_{2} \Phi((\xb,0);(u,v))(y^*)
\]
for some $v, y^*\in\mathbb Y$.
Then there exists $z^*\in \mathbb Y$ such that
\begin{subequations}\label{eq:upper_estimate_pseudo_coderivative_order_two}
	\begin{align}
		\label{eq:2ordEstimNC}
 		x^* &= \nabla^2\langle y^*,g\rangle(\bar x)(u) + \nabla g(\xb)^* z^*,
 		\\
 		\label{eq:multiplier_from_dir_lim_normal_cone_and_kernel}
		y^* &\in \mathcal N_{D}(g(\xb);\nabla g(\xb) u - v) \cap \ker \nabla g(\xb)^*.
	\end{align}
\end{subequations}
Further specifications of $z^*$ satisfying
\eqref{eq:upper_estimate_pseudo_coderivative_order_two}
are available under additional assumptions.
\begin{enumerate}
 \item\label{item:general_estimate_+CQ}
Each of the following two conditions
\begin{subequations}\label{eq:some:CQ}
	\begin{align}
		\label{eq:some:CQ_1}
 		& D\mathcal N_{D}(g(\xb),y^*)(0) \cap \ker \nabla g(\xb)^*=\{0\},
 		\\
 		\label{eq:some:CQ_2}
		\nabla g(\xb)u \neq v,
		\quad
		& D_{\textup{sub}}\mathcal N_{D}(g(\xb),y^*)
		\left(\frac{\nabla g(\xb)u - v}{\norm{\nabla g(\xb)u - v}} \right) 
		\cap \ker \nabla g(\xb)^*
		=\emptyset
	\end{align}
\end{subequations}
implies that we can find $z^* \in D\mathcal N_{D}(g(\xb),y^*)(\nabla g(\xb)u - v)$ satisfying \eqref{eq:upper_estimate_pseudo_coderivative_order_two}.
\item\label{item:polyhedral_estimate}
If $\mathbb Y:=\R^m$ and $D$ is locally polyhedral around $g(\bar x)$, 
then $\mathcal N_{D}(g(\xb);\nabla g(\xb) u - v)= \mathcal N_{\mathcal T_D(g(\xb))}(\nabla g(\xb) u - v)$,
and there are two elements $z^*_1,z^*_2\in\R^m$ satisfying \eqref{eq:upper_estimate_pseudo_coderivative_order_two} 
(for $z^*:=z_i^*$ with $i=1,2$, respectively) with
$z_1^*\in\mathcal N_{\mathcal T_D(g(\xb))}(\nabla g(\xb) u - v)$
and
$z_2^*\in\mathcal T_{\mathcal N_{\mathcal T_D(g(\xb))}(\nabla g(\xb) u - v)}(y^*)$.
\end{enumerate}
\end{theorem}
\begin{proof}
 Since $x^* \in \widetilde{D}^\ast_{2} \Phi((\xb,0);(u,v))(y^*)$, we find 
 $\{t_k\}_{k\in\N}\subset\R_+$, $\{u_k\}_{k\in\N},\{x_k^*\}_{k\in\N}\subset\mathbb X$, and $\{v_k\}_{k\in\N},\{y_k^*\}_{k\in\N}\subset\mathbb Y$
 with $t_k\downarrow 0$, $u_k\to u$, $v_k\to v$, $x_k^*\to x^*$, $y_k^*\to y^*$, as well as
 \[(x_k^*,-y_k^*/\tau_k) \in \widehat{\mathcal N}_{\gph \Phi}(\xb + t_ku_k,t_k v_k)\]
 for all $k\in\N$ where we used $\tau_k := t_k \norm{u_k}$ for brevity of notation.
 \cref{lem:coderivatives_constraint_maps} yields 
 $x_k^*=\nabla g(\bar x+t_ku_k)^* y_k^*/\tau_k$ and $y_k^*\in\tau_k\widehat{\mathcal N}_D(g(\bar x+t_ku_k)-t_k v_k)$
 for each $k\in\N$. Taking the limit in $\tau_kx_k^*=\nabla g(\bar x+t_ku_k)^* y_k^*$, we find $y^*\in\ker\nabla g(\bar x)^*$.
 Combining this with a Taylor expansion and 
 denoting $\tilde w_k:=g(\bar x+t_ku_k)-t_k v_k$ gives us
\begin{subequations}
	\begin{align}
	\label{eqn2 : Domain}
		& x_k^* - \nabla^2\langle y_k^*,g\rangle(\bar x)(u) + \oo(1)
		=
		\nabla g(\xb)^* \frac{y_k^*}{\tau_k}
		=
		\nabla g(\xb)^* \frac{y_k^* - y^*}{\tau_k},\\
	\label{eqn1 : Image}
		& y_k^*
		\in
		 \widehat{\mathcal N}_D(\tilde w_k)
		=
		 \widehat{\mathcal N}_D\left(
		 	g(\xb) + t_k \left(\nabla g(\xb) u - v + \oo(1)\right) 
		 	\right)
	\end{align}
\end{subequations}
for each $k\in\N$.
We readily obtain $y^* \in \mathcal N_{D}(g(\xb);\nabla g(\xb) u - v)$, i.e.,
\eqref{eq:multiplier_from_dir_lim_normal_cone_and_kernel}, as well as
\[
	x^* - \nabla^2\langle y^*,g\rangle(\bar x)(u) \in \Im \nabla g(\xb)^*,
\]
i.e., \eqref{eq:2ordEstimNC},
due to the closedness of 
$\Im \nabla g(\xb)^*$.

In the general case~\ref{item:general_estimate_+CQ}, we will use \eqref{eqn2 : Domain} only with
the right-hand side $\nabla g(\bar x)^*(y_k^* - y^*)/\tau_k$, 
but in the polyhedral case~\ref{item:polyhedral_estimate}, 
it is also reasonable to take a closer look at the expression $\nabla g(\bar x)^*y_k^*/\tau_k$.

Let us now prove~\ref{item:general_estimate_+CQ}.
Using the notation from above, let us first assume that $\{z_k^*\}_{k\in\N}$, given by
$z_k^*:=(y_k^* - y^*)/\tau_k$ for each $k\in\N$, remains bounded.
Then we may pass to a subsequence (without relabeling)
so that it converges to some $z^*\in\mathbb Y$. We get
\[
	y^* + \tau_k z_k^* = y_k^* 
	\in  
	\widehat{\mathcal N}_D\left(g(\xb) + \tau_k (\nabla g(\xb)u - v + \oo(1))\right)
\]
and $z^* \in D\mathcal N_{D}(g(\xb),y^*)(\nabla g(\xb)u - v)$ follows.
Clearly, taking the limit in \eqref{eqn2 : Domain} yields \eqref{eq:2ordEstimNC} as well.

On the other hand, if $\{z_k^*\}_{k\in\N}$ does not remain bounded,
we pass to a subsequence (without relabeling) such that
$\tau_k/\norm{y_k^* - y^*} \to 0$ and $\hat{z}_k^* \to \hat{z}^*$
for some $\hat z^*\in\mathbb S_{\mathbb Y}$
where we used $\hat z_k^*:= (y_k^* - y^*)/\norm{y_k^* - y^*}$ for each $k\in\N$. 
Multiplying \eqref{eqn2 : Domain} by $\tau_k/\norm{y_k^* - y^*}$
and taking the limit yields $\nabla g(\xb)^* \hat{z}^* = 0$. 
Taking into account $(\tilde w_k - g(\xb))/\tau_k \to \nabla g(\xb)u - v$,
we get
\begin{equation}\label{eq:some_convergence_of_surrogate_sequences}
	\frac{\norm{\tilde w_k - g(\xb)}}{\nnorm{y_k^* - y^*}}
	=
	\frac{\norm{\tilde w_k - g(\xb)}}{\tau_k}\frac{\tau_k}{\nnorm{y_k^* - y^*}}
	\to
	0.
\end{equation}
Let us assume that $\nabla g(\bar x)u\neq v$.
Then, for sufficiently large $k\in\N$, we have $\tilde w_k\neq g(\bar x)$, so we can 
set $\hat q_k:=(\tilde w_k-g(\bar x))/\norm{\tilde w_k-g(\bar x)}$ for any
such $k\in\N$ and find $\hat q\in\mathbb S_{\mathbb Y}$ such that $\hat q_k\to\hat q$
(along a subsequence without relabeling).
Moreover, we have
\[
	y^* + \nnorm{y_k^* - y^*} \hat{z}_k^* 
	= 
	y_k^* \in  \widehat{\mathcal N}_D\left(g(\xb) + \norm{\tilde w_k - g(\xb)} \hat{q}_k\right)
\]
from \eqref{eqn1 : Image}, so that \eqref{eq:some_convergence_of_surrogate_sequences}
yields $\hat{z}^* \in D_{\textup{sub}}\mathcal N_{D}(g(\xb),y^*)(\hat{q})$.
This contradicts \eqref{eq:some:CQ_2}.
In the case where $\nabla g(\bar x)u=v$ holds, \eqref{eq:some:CQ_2} is not applicable.
However, we still have
\[
	y^*+\nnorm{y_k^*-y^*}\hat z_k^*
	=
	y_k^* 
	\in 
	\widehat{\mathcal N}_D\left(g(\bar x)+\nnorm{y_k^*-y^*}\,\frac{\tilde w_k-g(\bar x)}{\nnorm{y_k^*-y^*}}\right),
\]
so that taking the limit $k\to\infty$ 
while respecting \eqref{eq:some_convergence_of_surrogate_sequences} 
yields $\hat z^*\in D\mathcal N_D(g(\bar x),y^*)(0)$
which contradicts \eqref{eq:some:CQ_1}.

In the polyhedral case~\ref{item:polyhedral_estimate},
we will show that one can always replace
the potentially unbounded sequences from \eqref{eqn2 : Domain} by bounded ones.
To start, we prove that
$y_k^* \in \mathcal N_{\mathcal T_D(g(\xb))}\left(\nabla g(\xb) u - v\right)$
for all sufficiently large $k\in\N$.
\Cref{lem:some_properties_of_polyhedral_sets}\,\ref{item:exactness_tangential_approximation} yields
the existence of a neighborhood $V\subset\R^m$ of $0$ such that
\begin{equation}\label{eq:exactness_of_tangents_D}
	\mathcal T_D(g(\xb)) \cap V = \big(D - g(\xb) \big) \cap V,
\end{equation}
as well as the fact that $\mathcal T_D(g(\xb))$ is polyhedral.
Thus, from \eqref{eqn1 : Image} we conclude
\begin{align*}
	y_k^*
	& \in 
	\widehat{\mathcal N}_D\left(g(\xb) + t_k \left(\nabla g(\xb) u - v + \oo(1)\right)\right)\\
	& = 
	\widehat{\mathcal N}_{g(\xb) + \mathcal T_D(g(\xb))}\left(g(\xb) + t_k \left(\nabla g(\xb) u - v + \oo(1)\right)\right)\\
	& = 
	\widehat{\mathcal N}_{\mathcal T_D(g(\xb))}\left(\nabla g(\xb) u - v + \oo(1)\right)
	\ \subset \
	\mathcal N_{\mathcal T_D(g(\xb))}\left(\nabla g(\xb) u - v\right)
\end{align*}
for all sufficiently large $k\in\N$.

Next, let us set $K:=\mathcal N_{\mathcal T_D(g(\bar x))}(\nabla g(\bar x)u-v)$ for brevity of notation,
and note that $K$ is a polyhedral cone.
From above we know that $y_k^* \in K$ holds for all sufficiently large $k\in\N$.
Then we also get $y^*, y_k^*/\tau_k \in K$
and, by \Cref{lem:some_properties_of_polyhedral_sets}\,\ref{item:exactness_tangential_approximation},
$(y_k^* - y^*)/\tau_k \in \mathcal T_K(y^*)$,
where $\mathcal T_K(y^*)$ is also a polyhedral cone.
Thus, referring to \eqref{eqn2 : Domain}, we may invoke Hoffman's lemma,
see \cite[Lemma~3C.4]{DontchevRockafellar2014},
to find some bounded sequences 
$\{z_{1,k}^*\}_{k\in\N} \subset K$ and
$\{z_{2,k}^*\}_{k\in\N} \subset \mathcal T_K(y^*)$
satisfying
 \[
 	\nabla g(\xb)^* z_{i,k}^* = x_k^* - \nabla^2\langle y^*_k,g\rangle(\bar x)(u) + \oo(1)
 \]
for $i=1,2$.
Thus, accumulation points $z_i^*\in\R^m$ of $\{z_{i,k}^*\}_{k\in\N}$
for $i=1,2$ satisfy \eqref{eq:2ordEstimNC} and
$z_1^* \in K$ and $z_2^* \in \mathcal T_K(y^*)$.
\end{proof}

Below, we comment on the findings of \cref{The : NCgen}.
To start, we illustrate that the additional information on the multiplier $z^*$
provided in statements~\ref{item:general_estimate_+CQ} and~\ref{item:polyhedral_estimate}
is the same whenever $D$ is a convex polyhedral set in $\mathbb Y:=\R^m$.
\begin{remark}\label{rem:comparison_estimate_pseudocoderivative_polyhedral}
	We use the notation from \cref{The : NCgen}.
	Suppose that $D$ is a \emph{convex}
	polyhedral set in $\mathbb Y:=\R^m$.
	First, we claim that
	\begin{align*}
	\mathcal N_{\mathcal T_D(g(\xb))}(\nabla g(\xb)u - v)
	& \subset
	\mathcal N_{\mathcal T_D(g(\xb))}(\nabla g(\xb)u - v) + \spa(y^*)
	\\
	&=
	\mathcal T_{\mathcal N_{\mathcal T_D(g(\xb))}(\nabla g(\xb) u - v)}(y^*)
	=
	D\mathcal N_{D}(g(\xb),y^*)(\nabla g(\xb)u - v).
	\end{align*}
	The first two relations are straightforward and so let us prove the last one.
	Based on the so-called reduction lemma, see \cite[Lemma~2E.4]{DontchevRockafellar2014}, 
	and \cite[Proposition~2A.3]{DontchevRockafellar2014},
	for each pair $(\bar z, \bar z^*) \in \gph \mathcal N_D$, we get
	\begin{align*}
		&\big(\gph \mathcal N_D - (\bar z,\bar z^*)\big) \cap \mathcal{O} 
		\\
		&\qquad
		=
		 \{(w,w^*) \,|\, w \in \mathcal K_D(\bar z,\bar z^*),\, w^* \in \mathcal K_D(\bar z,\bar z^*)^{\circ},\, \innerprod{w}{ w^* }=0\} \cap \mathcal{O},
	\end{align*}
	where $\mathcal{O}\subset\R^m\times\R^m$ is a neighborhood of $(0,0)$ and 
	$\mathcal K_D(\bar z,\bar z^*) := \mathcal T_D(\bar z) \cap [\bar z^*]^\perp$ represents the \emph{critical cone} to $D$ at $(\bar z,\bar z^*)$.
	By \cref{lem:some_properties_of_polyhedral_sets}\,\ref{item:exactness_tangential_approximation}, this simply means
	\[
		\mathcal T_{\gph \mathcal N_D}(\bar z,\bar z^*) 
		=
		\{(w,w^*) \,|\, w \in \mathcal K_D(\bar z,\bar z^*),\, w^* \in \mathcal K_D(\bar z,\bar z^*)^{\circ},\, \innerprod{w}{ w^* }=0\}.
	\]
	Thus, $z^* \in D\mathcal N_{D}(g(\xb),y^*)(\nabla g(\xb)u - v)$ 
	means $\nabla g(\xb)u - v\in \mathcal T_D(g(\xb)) \cap [y^*]^\perp$, which gives us
	\begin{align*}
		y^*
		\in \mathcal N_D(g(\bar x))\cap[\nabla g(\bar x)u-v]^\perp
		= 
		\mathcal N_{\mathcal T_D(g(\xb))}(\nabla g(\xb) u - v),
	\end{align*}
	and
	\begin{align*}
		z^*
		& \in 
		\mathcal N_{\mathcal K_D(g(\xb),y^*)}(\nabla g(\xb)u - v)
		=
		\big(\mathcal T_D(g(\xb)) \cap [y^*]^\perp\big)^{\circ} \cap [\nabla g(\xb)u - v]^\perp\\
		& =
		\big(\mathcal N_D(g(\xb))+\spa(y^*)\bigr)\cap [\nabla g(\bar x)u-v]^\perp\\
		& =
		\mathcal N_D(g(\bar x))\cap[\nabla g(\bar x)u-v]^\perp + \spa(y^*)\\
		& = 
		\mathcal N_{\mathcal T_D(g(\xb))}(\nabla g(\xb)u - v) + \spa(y^*)
		= 
		\mathcal T_{\mathcal N_{\mathcal T_D(g(\xb))}(\nabla g(\xb) u - v)}(y^*)
	\end{align*}
	by the basic properties of convex polyhedral cones and \cref{lem:some_properties_of_polyhedral_sets}\,\ref{item:normal_cones_to_polyhedral_sets}.
	
	Hence, in the convex polyhedral case, 
	the information on $y^*$ and $z^*$ from statements~\ref{item:general_estimate_+CQ} and~\ref{item:polyhedral_estimate} (case $z_2^*$)
	of \cref{The : NCgen} is the same,
	while the information from statement~\ref{item:polyhedral_estimate} (case $z_1^*$)
	is seemingly sharper.
	Let us now demonstrate that it is actually also equivalent to the others.
	
	Note that \eqref{eq:multiplier_from_dir_lim_normal_cone_and_kernel} can be equivalently written as
	$y^*	\in \mathcal N_D(g(\bar x))\cap[v]^\perp \cap \ker \nabla g(\xb)^*$ 
	due to \cref{lem:some_properties_of_polyhedral_sets}\,\ref{item:normal_cones_to_polyhedral_sets}
	and
	$[\nabla g(\bar x)s-v]^\perp\cap\ker\nabla g(\bar x)^*=[v]^\perp\cap\ker\nabla g(\bar x)^*$
	for all $s\in\mathbb X$.
	This also means that, for any such $y^*$, the sets
	\begin{align*}
		A_1(y^*,v)	&:= \left\{s \in \mathbb X \,\middle|\, \nabla g(\xb)s - v \in \mathcal T_D(g(\xb))\right\},
		\\
		A_2(y^*,v) 	&:= \left\{s \in \mathbb X \,\middle|\, \nabla g(\xb)s - v \in \mathcal K_D(g(\xb),y^*)\right\}
	\end{align*}
	coincide, and viewing $x^*$, $y^*$, $u$, and $v$ as parameters, the linear programs
	\begin{equation}\label{eq:LPi}\tag{LP$(i)$}
		\min\limits_s \{\nabla^2\langle y^*,g\rangle(\bar x)[u,s] - \innerprod{x^*}{s} \,|\, s \in A_i(y^*,v)\}
	\end{equation}
	are the same for $i=1,2$.
	On the other hand, \eqref{eq:2ordEstimNC} with $z^* \in \mathcal N_{\mathcal T_D(g(\xb))}(\nabla g(\xb)u - v)$
	and $z^* \in \mathcal N_{\mathcal K_D(g(\xb),y^*)}(\nabla g(\xb)u - v)$, respectively,
	precisely characterizes the fact that $u$ is a minimizer of 
	\hyperref[eq:LPi]{\textup{(LP$(1)$)}} and \hyperref[eq:LPi]{\textup{(LP$(2)$)}}.
	Hence, this information on $z^*$ is the same.
\end{remark}

Some additional comments on \cref{The : NCgen} are stated subsequently.
\begin{remark}\label{rem:upper_estimate_pseudo_coderivative_constraint_maps}
	We use the notation from \cref{The : NCgen}.
	\begin{enumerate}
	\item 
	Note that, in the case $\nabla g(\bar x)u\neq v$, assumption \eqref{eq:some:CQ_2}, 
	which is stated in terms of the graphical subderivative,
	is milder than \eqref{eq:some:CQ_1} in terms or the standard graphical derivative,
	and it preserves the connection to the direction $\nabla g(\xb)u - v$.
	Let us also note that the case $\nabla g(\bar x)u=v$ is, anyhow, special since this
	would annihilate the directional information in \eqref{eq:multiplier_from_dir_lim_normal_cone_and_kernel}
	completely.
	\item\label{item:graphical_subderivative_assumption_polyhedral}
	If $\mathbb Y:=\R^m$ and $D$ is locally polyhedral around $g(\bar x)$,
	conditions \eqref{eq:some:CQ} reduce to
	\begin{equation*}
		D\mathcal N_{D}(g(\xb),y^*)(\nabla g(\xb)u - v) 
		\cap \ker \nabla g(\xb)^*
		\subset \{0\}
	\end{equation*}
	thanks to \cref{lem:normal_cone_map_of_polyhedral_set}.
	\end{enumerate}
\end{remark}

In the polyhedral case, we can derive yet sharper information on $z^*$ 
if we start with the new pseudo-coderivative instead of the one utilized by Gfrerer.
This is also the main reason for introducing the new definition.
Throughout the paper, we will rely on the following result.
Particularly, it plays an important role in \cref{Pro:Milder_than_SOSCMS} and \cref{cor:FOSCMS_SOSCMS_imply_sAR},
which we were not able to get using the estimates from \cref{The : NCgen}.

\begin{theorem}\label{The : NCgen_2}
Let $g$ be twice continuously differentiable.
Given $(\xb,0) \in \gph \Phi$,
assume that $\mathbb Y:=\R^m$ and $D$ is locally polyhedral around $g(\bar x)$.
For a direction $u \in \mathbb S_{\mathbb X}$, let
\[
	x^* \in D^\ast_{2} \Phi((\xb,0);(u,v))(y^*) 
\]
for some $v, y^*\in\R^m$.
Then there exists $s\in\mathbb X$ satisfying
$y^* \in \mathcal N_{\mathbf T(u)}(w_s(u,v))\cap \ker \nabla g(\bar x)^*$ where
\begin{equation}\label{eq:Tu_and_ws}
 \mathbf T(u):=\mathcal T_{\mathcal T_D(g(\xb))}(\nabla g(\xb) u),\qquad
 w_s(u,v):= \nabla g(\xb) s + 1/2 \nabla^2 g(\xb)[u,u] - v,
\end{equation}
together with two elements
$z_1^*\in\mathcal N_{\mathbf T(u)}(w_s(u,v))$
and
$z_2^*\in\mathcal T_{\mathcal N_{\mathbf T(u)}(w_s(u,v))}(y^*)$
satisfying $x^* = \nabla^2\langle y^*,g\rangle(\bar x)(u) + \nabla g(\xb)^* z_i^*$ for $i=1,2$.
Moreover, $v \in D_2 \Phi(\xb,0)(u)$ is equivalent to
the existence of $s\in\mathbb X$ with $w_s(u,v) \in \mathbf T(u)$.%
\end{theorem}
\begin{proof}
Similar arguments as in the proof of \cref{The : NCgen} yield 
\eqref{eqn2 : Domain}
together with $y_k^* \in \widehat{\mathcal N}_D(w_k)$ for each $k\in\N$ where
\begin{align*}
	w_k
	& := 
	g(\bar x+t_ku_k)-t_k^2 v_k
	=
	g(\xb) + t_k \nabla g(\xb) u + t_k^2 z_k,\\
	z_k
	& := 
	\frac{\big( w_k - g(\xb) \big) / t_k - \nabla g(\xb) u}{t_k}
	=
	\nabla g(\xb) \frac{u_k - u}{t_k} + \frac12 \nabla^2 g(\xb)[u,u] -  v + \oo(1).
\end{align*}
As in the final part of the proof of \cref{The : NCgen},
all we need to show is $y^*_k \in \mathcal N_{\mathbf T(u)}(w_s(u,v))$
for all sufficiently large $k\in\N$ and some appropriately chosen $s\in\mathbb X$.

Noting that $D$ is polyhedral while $\mathbf T(u)$ is a polyhedral cone,
we can apply \cref{lem:some_properties_of_polyhedral_sets}\,\ref{item:exactness_tangential_approximation}
to find neighborhoods $V,W\subset\R^m$ of $0$ such that \eqref{eq:exactness_of_tangents_D} and
\[
 	\mathbf T(u) \cap W 
 	= 
 	\mathcal T_{\mathcal T_D(g(\xb))}(\nabla g(\xb) u) \cap W 
 	=
 	\big(\mathcal T_D(g(\xb)) - \nabla g(\xb) u \big) \cap W.
\]
 Consequently,  we have  $w_k - g(\xb) \in \mathcal T_D(g(\xb))$ and, hence, also
 $\big( w_k - g(\xb) \big) / t_k \in \mathcal T_D(g(\xb))$ for sufficiently large $k\in\N$.
 Similarly, we conclude that $z_k \in \mathbf T(u)$.
 Taking into account that for each cone $K$, $q \in K$, and $\alpha > 0$, one has $\mathcal T_K(q) = \mathcal T_K(\alpha q) $, 
 we find
 \begin{align*}
 	\mathcal T_D(w_k)
 	&=
 	\mathcal T_{g(\bar x)+\mathcal T_D(g(\bar x))}(w_k)
 	=
 	\mathcal T_{\mathcal T_D(g(\bar x))}((w_k-g(\bar x))/t_k)
 	\\
 	&=
 	\mathcal T_{\mathcal T_D(g(\bar x))}\left(\nabla g(\xb) u + t_k z_k \right)
 	=
 	\mathcal T_{\mathcal T_D(g(\bar x))-\nabla g(\bar x)u}\left(z_k\right)
 	=
 	\mathcal T_{\mathbf T(u)}\left(z_k\right)
 \end{align*}
 for all sufficiently large $k\in\N$, and we obtain $y_k^* \in \widehat{\mathcal N}_D(w_k) = \widehat{\mathcal N}_{\mathbf T(u)}(z_k)$.

 Since $\mathbf T(u)$ is polyhedral, so is $\gph \mathcal N_{\mathbf T(u)}$,
 see \cref{lem:normal_cone_map_of_polyhedral_set}, 
 and it can be written as the union of finitely many convex polyhedral sets, 
 say $C_1,\ldots,C_\ell\subset\R^m\times\R^m$.
 Thus, we have
 \[
 	(z_k,y_k^*) \in \gph \widehat{\mathcal N}_{\mathbf T(u)} \subset \gph \mathcal N_{\mathbf T(u)} = \bigcup_{j=1}^\ell C_j
 \]
 for sufficiently large $k\in\N$.
 We may pick an index $\bar j\in\{1,\ldots,\ell\}$ such that $(z_k,y_k^*) \in C_{\bar j}$ holds for infinitely many $k\in\N$
 and suppose that $C_{\bar j}$ can be represented as $C_{\bar j} = \{(z,y) \,|\, A z + B y \leq c\}$
 for some matrices $A$, $B$, as well as $c$ of appropriate dimensions.
 Hence, by passing to a subsequence (without relabeling), we get
 \[
 	A \nabla g(\xb)\frac{u_k - u}{t_k} \leq c - A \left(\frac12 \nabla^2 g(\xb)[u,u] - v + \oo(1)\right)- B y_k^*.
 \]
 For each $k\in\N$, a generalized version of Hoffman's lemma, 
 see \cite[Theorem~3]{Ioffe1979},
 now yields the existence of $s_k\in\mathbb X$ with
 \begin{align*}
 	A \nabla g(\xb) s_k &\leq c - A  ((1/2) \nabla^2 g(\xb)[u,u] - v+ \oo(1))- B y_k^*,\\
 	\norm{s_k} &\leq \beta\bigl\Vert c - A ((1/2) \nabla^2 g(\xb)[u,u] - v + \oo(1)) - B y_k^*\bigr\Vert
 \end{align*}
 for some constant $\beta > 0$ not depending on $k$.
 Thus, $\{s_k\}_{k\in\N}$ is bounded and satisfies
 \[
 	\forall k\in\N\colon\quad
 	\bigl(\nabla g(\xb)s_k + 1/2 \nabla^2 g(\xb)[u,u] - v + \oo(1),y_k^*\bigr) \in C_{\bar j} \subset \gph \mathcal N_{\mathbf T(u)}.
 \]
 We may assume that $\{s_k\}_{k\in\N}$ converges to some $s\in\mathbb X$.
 Exploiting \eqref{eq:Tu_and_ws}, we infer
 \[
 	y^*_k \in \mathcal  N_{\mathbf T(u)}(w_s(u,v) + \oo(1)) \subset \mathcal N_{\mathbf T(u)}(w_s(u,v))
 \]
 for all sufficiently large $k\in\N$ from polyhedrality of $\mathbf T(u)$ 
 and the definition of the limiting normal cone.
 
To show the second statement, note that
 $v \in D_2 \Phi(\xb,0)(u)$ is equivalent to
 $0 \in D^\ast_{2} \Phi((\xb,0);(u,v))(0)$,
 so that any of these two conditions readily yields the existence of $s\in\mathbb X$ with $w_s(u,v) \in \mathbf T(u)$.
 Conversely, suppose that there exists
 $s\in\mathbb X$ with $w_s(u,v) \in \mathbf T(u)$.
 Let $\{t_k\}_{k \in \N} \subset \R_+$ be an arbitrary sequence with $t_k\downarrow 0$, 
 and define the sequences $\{u_k\}_{k \in \N} \subset \mathbb X$ and $\{v_k\}_{k \in \N},\{\hat w_k\}_{k\in\N} \subset \mathbb Y$ by
 $u_k := u + t_k s$ and
 \[
	v_k := \big(g(\bar x + t_k u_k) - \hat w_k\big)/(t_k\norm{u_k})^2,
	\qquad
	\hat w_k := g(\bar x) + t_k \nabla g(\bar x) u + t_k^2 w_s(u,v)
 \]
 for all $k \in \N$.
 First, a second-order Taylor expansion together with $\norm{u_k}\to 1$ yields $v_k \to v$.
 Next, using similar arguments as before, polyhedrality of $\mathcal T_D(g(\bar x))$
 and, locally around $g(\bar x)$, $D$, together with $w_s(u,v) \in \mathbf T(u)$, yields
 $g(\bar x + t_k u_k) - (t_k\norm{u_k})^2 v_k = \hat w_k \in D$, i.e.,
 $(\bar x+t_ku_k,(t_k\norm{u_k})^2v_k)\in\gph\Phi$, for sufficiently large $k\in\N$.
 Taking the limit $k\to\infty$ gives $v\in D_2\Phi(\bar x,0)(u)$, and this completes the proof.
\end{proof}

\begin{remark}\label{rem:refined_polyhedral_situation_normal_cone_relation}
	Let us mention that if $\mathbb Y:=\R^m$ and $D$ is locally polyhedral around $g(\bar x)$,
	we get the relations 
	\begin{align*}
            \mathcal N_{\mathbf T(u)}(w_s(u,v))
            & = 
            \mathcal N_{\mathcal T_D(g(\bar x))}(\nabla g(\bar x)u;w_s(u,v))
            \\
            & \subset 
            \mathcal N_{\mathcal T_D(g(\bar x))}(\nabla g(\bar x)u)\cap[w_s(u,v)]^\perp\\
            & \subset 
            \mathcal N_{\mathcal T_D(g(\bar x))}(\nabla g(\bar x)u)
             = 
            \mathcal N_D(g(\bar x);\nabla g(\bar x)u) 
	=
	\mathcal N_{\mathbf T(u)}(0).
     \end{align*}
     from 
     \cref{lem:some_properties_of_polyhedral_sets}\,\ref{item:normal_cones_to_polyhedral_sets}.   
     This also yields 
     $\mathcal T_{\mathcal N_{\mathbf T(u)}(w_s(u,v))}(y^*) 
     \subset\mathcal T_{\mathcal N_{\mathcal T_D(g(\bar x))}(\nabla g(\bar x)u)}(y^*)$.
\end{remark}

Again, in the convex polyhedral case, the two options provided by \cref{The : NCgen_2} coincide.
This can be shown using the same arguments as in \cref{rem:comparison_estimate_pseudocoderivative_polyhedral} 
but with the sets
\begin{equation}\label{eq:A_i_sets}
		\begin{aligned}
		\widetilde A_1(y^*,u,v) 
		&:= 
		\left\{\tilde s \in \mathbb X \,\middle|\, w_{\tilde s}(u,v) \in \mathbf T(u)\right\},
		\\
		\widetilde A_2(y^*,u,v) 
		&:= 
		\left\{\tilde s \in \mathbb X \,|\, w_{\tilde s}(u,v) \in \mathcal K_{\mathcal T_D(g(\bar x))}(\nabla g(\bar x)u,y^*)\right\}
		\end{aligned}
\end{equation}
which coincide because the required existence of $s\in\mathbb X$ with
$y^* \in \mathcal N_{\mathbf T(u)}(w_s(u,v))\cap\ker\nabla g(\bar x)^*=(\mathbf T(u))^\circ\cap [w_s(u,v)]^\perp\cap\ker\nabla g(\bar x)^*$ yields
$y^*\in[1/2\nabla^2 g(\bar x)[u,u]-v]^\perp$ and, thus,
$\innerprod{y^*}{w_{\tilde s}(u,v)} = 0$ for all $\tilde s \in \mathbb X$.
This means that our conditions from \cref{The : NCgen_2}
precisely state that the associated linear programs \eqref{eq:LPi}, $i=1,2$, with
$A_i(y^*,v)$ replaced by $\widetilde A_i(y^*,u,v)$, have a solution.

From \cref{The : NCgen,The : NCgen_2}
we obtain the following explicit sufficient
conditions for metric pseudo-(sub)regularity
of constraint mappings.

\begin{corollary}\label{cor:sufficient_condition_pseudo_subregularity}
 Let $g$ be twice continuously differentiable.
 Consider $(\xb,0) \in \gph \Phi$ and a direction $u \in \mathbb S_{\mathbb X}$.
 The characterization \eqref{eq:FO_characterization_of_dir_metric_pseudo_reg}
 of metric pseudo-regularity of order $2$ of $\Phi$ in direction $(u,0)$ at $(\bar x,0)$ holds
 under conditions~\ref{item:trivial_SC_for_MPR},~\ref{item:SC_for_MPR}, and~\ref{item:SC_for_MPR_polyhedral},
 while the sufficient condition \eqref{eq:FOSCMS_gamma}
 for metric pseudo-subregularity of order $2$ of $\Phi$ in direction $u$ at $(\bar x,0)$
 is valid also under~\ref{item:SC_for_MPSR_polyhedral}.
 \begin{enumerate}
 	 \item\label{item:trivial_SC_for_MPR} 
 	 	One has 
  		\begin{equation*}
   			\left. \begin{aligned}
        		&\nabla g(\xb)^* y^* = 0, \, 
        		\nabla^2\langle y^*,g\rangle(\bar x)(u) + \nabla g(\xb)^* z^* = 0, \\
        		&y^* \in \mathcal N_{D}(g(\xb);\nabla g(\xb) u)
       			\end{aligned}
				\right\} 
				\quad \Longrightarrow \quad 
				y^* = 0.
  		\end{equation*}
  	\item\label{item:SC_for_MPR} 
  		One has
  		\begin{equation}\label{eq:CQ_pseudo_subregularity_II}
  			\left. \begin{aligned}
        		&\nabla g(\xb)^* y^* = 0, \, 
        		\nabla^2\langle y^*,g\rangle(\bar x)(u) + \nabla g(\xb)^* z^* = 0, \\
        		&y^* \in \mathcal N_{D}(g(\xb);\nabla g(\xb) u), \,
        		z^*\in D\mathcal N_D(g(\bar x),y^*)(\nabla g(\bar x)u)
       		\end{aligned}
       		\right\} 
       		\quad \Longrightarrow \quad 
       		y^* = 0.
  		\end{equation}
  		Furthermore, we either have
  		\begin{equation}\label{eq:CQ_pseudo_subregularity_Ia}
  			\left. \begin{aligned}
        		&\nabla g(\xb)^* y^* = 0, \, 
        		\nabla g(\xb)^* \hat z^* = 0, \\
        		&y^* \in \mathcal N_{D}(g(\xb);\nabla g(\xb) u), \,
        		 \hat z^* \in D\mathcal N_D(g(\bar x),y^*)(0) 
       			\end{aligned}
				\right\} 
				\quad \Longrightarrow \quad 
				\hat z^* = 0
  		\end{equation}
		or $\nabla g(\bar x)u\neq 0$ and
		\begin{equation}\label{eq:CQ_pseudo_subregularity_Ib}
  			\left. \begin{aligned}
        		&\nabla g(\xb)^* y^* = 0, \, 
        		\nabla g(\xb)^* \hat z^* = 0, \\
        		&y^* \in \mathcal N_{D}(g(\xb);\nabla g(\xb) u)
       			\end{aligned}
				\right\} 
				\ \Longrightarrow \ 
				\hat z^* \notin D_\textup{sub}\mathcal N_D(g(\bar x),y^*)
				\left(\frac{\nabla g(\bar x)u}{\norm{\nabla g(\bar x)u}}\right).
  		\end{equation}
  		\item\label{item:SC_for_MPR_polyhedral}
  			It holds $\mathbb Y:=\R^m$, $D$ is locally polyhedral around $g(\xb)$,
  			and 
			\begin{equation}\label{eq:some_sufficient_condition_for_pseudo_regularity}
			\left.
			\begin{aligned}
			&\nabla g(\bar x)^*y^*=0,\,
			\nabla^2\langle y^*,g\rangle(\bar x)(u) + \nabla g(\bar x)^* z^*=0,\\
			&y^*\in\mathcal N_{\mathcal T_D(g(\bar x))}(\nabla g(\bar x)u),\\
			&z^*\in\mathcal N_{\mathcal T_D(g(\bar x))}(\nabla g(\bar x)u)\,
			\big(
        			\textrm{or } \, z^* \in \mathcal T_{\mathcal N_{\mathcal T_D(g(\bar x))}(\nabla g(\bar x)u)}(y^*)
        			\big) 
			\end{aligned}
			\right\}
			\quad\Longrightarrow\quad
			y^*=0.
			\end{equation}
 		\item\label{item:SC_for_MPSR_polyhedral} 
 			It holds $\mathbb Y:=\R^m$, $D$ is locally polyhedral around $g(\xb)$, 
 			and for each $s\in\mathbb X$ one has
 			\begin{equation}\label{eq:CQ_pseudo_subregularity_polyhedral_II}
   				\left. \begin{aligned}
        			&\nabla g(\xb)^* y^* = 0, \, 
        			\nabla^2\langle y^*,g\rangle(\bar x)(u) + \nabla g(\xb)^* z^* = 0, \\
        			&y^*\in \mathcal N_{\mathbf T(u)}(w_s(u,0)),\\
        			&z^* \in \mathcal N_{\mathbf T(u)}(w_s(u,0))\,\big(
        			\textrm{or } \, z^* \in \mathcal T_{\mathcal N_{\mathbf T(u)}(w_s(u,0))}(y^*)
        			\big)  
       				\end{aligned}
				\right\} 
				\quad \Longrightarrow \quad y^* = 0
  			\end{equation}
  			where $\mathbf T(u)$ and $w_s(u,0)$ are given as in \eqref{eq:Tu_and_ws}.
 		\end{enumerate}
\end{corollary}

	Due to \cref{rem:refined_polyhedral_situation_normal_cone_relation},
	\eqref{eq:some_sufficient_condition_for_pseudo_regularity} 
	indeed implies validity of \eqref{eq:CQ_pseudo_subregularity_polyhedral_II}
	for arbitrarily chosen $s\in\mathbb X$.

	\begin{remark}\label{rem:polyhedral_case_in_general_framework}
		Let us note that if $\mathbb Y:=\R^m$ and $D$ is locally polyhedral around $g(\bar x)$,
		then \eqref{eq:CQ_pseudo_subregularity_Ia} and \eqref{eq:CQ_pseudo_subregularity_Ib}
		appearing in \cref{cor:sufficient_condition_pseudo_subregularity}\,\ref{item:SC_for_MPR}
		reduce to 
		\begin{equation}\label{eq:CQ_pseudo_subregularity_I_polyhedral_new}
			\left. \begin{aligned}
        		&\nabla g(\xb)^* y^* = 0, \, 
        		\nabla g(\xb)^* \hat z^* = 0, \\
        		&y^* \in \mathcal N_{D}(g(\xb);\nabla g(\xb) u), \,
        		 \hat z^* \in D\mathcal N_D(g(\bar x),y^*)(\nabla g(\xb)u) 
       			\end{aligned}
				\right\} 
				\quad \Longrightarrow \quad 
				\hat z^* = 0
		\end{equation}
		thanks to \cref{rem:upper_estimate_pseudo_coderivative_constraint_maps}\,\ref{item:graphical_subderivative_assumption_polyhedral}.
	\end{remark}

\subsection{The convex polyhedral case: a comparison with related results}\label{sec:comparison_2reg}

	Throughout the subsection,
	we assume that $D$ is a convex polyhedral set
	in $\mathbb Y:=\R^m$, and
aim to compare our findings, at least partially, 
with available results from the literature.
To start, we recall the definition of directional 2-regularity
taken from \cite[Definition~1]{ArutyunovAvakovIzmailov2008}.
\begin{definition}\label{def:2-regularity}
	Set $\mathbb Y:=\R^m$, let $D$ be convex and polyhedral, and fix $(\bar x,0)\in\gph\Phi$
	as well as $u\in\mathbb X$.
	Then the \emph{2-regularity condition} is said to hold at $\bar x$ 
	\emph{in direction $u$} if the following is valid:
	\begin{equation}\label{eq:2regularity_primal}
			\Im \nabla g(\bar x)
			+
			\nabla^2g(\bar x)[u,\nabla g(\bar x)^{-1}\mathcal T_D(g(\bar x))]
			-
			\mathcal T_D(g(\bar x))
			=
			\R^m.
		\end{equation}
\end{definition}

Let us mention that the original definition of directional 2-regularity from
\cite[Definition~1]{ArutyunovAvakovIzmailov2008} is different from the one
stated in \cref{def:2-regularity}. However, both conditions are equivalent
by \cite[Proposition~1]{ArutyunovAvakovIzmailov2008}.
Furthermore, it should be noted that, in the setting of \cref{def:2-regularity}, 
the 2-regularity condition in direction $u:=0$ reduces to 
Robinson's constraint qualification,
see \cite[Proposition~2.97]{BonnansShapiro2000}.
Observe that, since $\Im\nabla g(\bar x)$, $\nabla g(\bar x)^{-1}\mathcal T_D(g(\bar x))$,
and $\mathcal T_D(g(\bar x))$ are cones, 2-regularity in a nonzero direction $u$
is equivalent to 2-regularity in direction $\alpha u$ for arbitrary $\alpha>0$.
Hence, it is reasonable to consider merely directions from $\mathbb S_\mathbb{X}$
in \cref{def:2-regularity}.
In \cref{prop:4_conditions_convex_polyhedral_case} below,
we derive a dual characterization of 2-regularity in direction $u$, which states that the conditions
\begin{equation}\label{eq:2-reg_assumptions}\tag{C$(u, y^*)$}
	\nabla g(\bar x)^* y^*=0,\,
	\nabla^2\langle y^*,g\rangle(\bar x)(u)+\nabla g(\bar x)^*z^*=0,\,
	y^*,z^*\in\mathcal N_D(g(\bar x))
\end{equation}
can be satisfied only for $y^* = 0$.
Note that \eqref{eq:2-reg_assumptions} can be stated in a $z^*$-free manner by means of
\begin{equation*}
	\nabla g(\bar x)^* y^*=0,\,
	0 \in \nabla^2\langle y^*,g\rangle(\bar x)(u)+\nabla g(\bar x)^* \mathcal N_D(g(\bar x)),\,
	y^* \in\mathcal N_D(g(\bar x))
\end{equation*}
which is why we did not include $z^*$ in the abbreviation \eqref{eq:2-reg_assumptions}.

Second, we will compare our findings with the ones from \cite{Gfrerer2014a}.
Again, we just consider the situation where $D$ is a convex polyhedral set.
In \cite[Theorem 2\,2.]{Gfrerer2014a},
pseudo-subregularity of the feasibility mapping
$\Phi$ of order $2$ at $(\bar x,0)\in\gph\Phi$ in some direction $u\in\mathbb S_{\mathbb X}$ 
which satisfies $\nabla g(\bar x)u\in\mathcal T_D(g(\bar x))$
(for other directions, the concept is trivial)
was shown to be present under the following condition:
\begin{equation}\label{eq:Gfrerers_sufficient_condition_implicit}
		\text{\eqref{eq:2-reg_assumptions}}, \
		y^* \in \argmax\limits_{
							\hat y^* \in \mathcal N_D(g(\xb))
							\cap
    						\ker \nabla g(\xb)^*
    						} 
    						1/2\nabla^2\dual{ \hat y^* }{g}(\bar x)[u,u]
		\quad\Longrightarrow\quad
		y^*=0.
\end{equation}

We will now derive alternative representations of
\eqref{eq:some_sufficient_condition_for_pseudo_regularity} and \eqref{eq:2regularity_primal},
which are sufficient for directional pseudo-regularity of $\Phi$ of order 2,
as well as 
\eqref{eq:CQ_pseudo_subregularity_polyhedral_II} and 
\eqref{eq:Gfrerers_sufficient_condition_implicit},
being sufficient for directional pseudo-subregularity of $\Phi$ of order 2,
which allow for a comparison of all these conditions.

	To start, let us present a technical lemma, collecting
	some consequences of having $s \in \mathbb X$ with $w_s(u,0) \in \mathbf T(u)$,
	see \eqref{eq:Tu_and_ws} for the definition of
	$w_{s}(u,0)$ and $\mathbf T(u)$.

\begin{lemma}\label{lem:w_s_versus_w_tilde_s}
	Set $\mathbb Y:=\R^m$, let $D$ be convex and polyhedral, and fix $(\bar x,0)\in\gph\Phi$
	as well as $u\in\mathbb S_\mathbb X$ such that $\nabla g(\bar x)u\in\mathcal T_D(g(\bar x))$.
	The existence of $s \in \mathbb X$ with $w_s(u,0) \in \mathbf T(u)$
	is equivalent to 
	the existence of $\tilde s \in \mathbb X$ with $w_{\tilde s}(u,0) \in \mathcal T_D(g(\xb))$,
	and these conditions imply $\nabla^2 \dual{y^*}{g}(\bar x)[u,u]\leq 0$
	for arbitrary $y^* \in \mathcal N_D(g(\bar x)) \cap \ker \nabla g(\bar x)^*$.
	If, additionally, \eqref{eq:2-reg_assumptions} holds, then
	we even have $\nabla^2 \dual{y^*}{g}(\bar x)[u,u]=0$.
\end{lemma}
\begin{proof}
	Let us start to prove the first assertion.
	Note that $\mathbf T(u) = \mathcal T_D(g(\bar x))+\spa(\nabla g(\bar x)u)$
	holds due to polyhedrality of $D$ yielding polyhedrality of $\mathcal T_D(g(\bar x))$.
	Hence, if $s\in\mathbb X$ satisfies $w_s(u,0) \in \mathbf T(u)$,
	then $\tilde s := s + \alpha u$ for some $\alpha \in \R$ satisfies 
	$w_{\tilde s}(u,0) = w_{s}(u,0) + \alpha \nabla g(\xb) u \in \mathcal T_D(g(\xb))$.
	The converse relation is trivial due to $\mathcal T_D(g(\bar x))\subset\mathbf T(u)$.
	
	The second assertion is a consequence of the definition of $w_s(u,0)$.
	
	To show the final assertion,
	note that \eqref{eq:2-reg_assumptions} gives
	\[
		\nabla^2 \langle y^*,g\rangle(\bar x)[u,u] = - \dual{z^*}{\nabla g(\bar x)u} \geq 0
	\]
	as $z^*\in\mathcal N_D(g(\bar x))$ and $\nabla g(\bar x)u\in\mathcal T_D(g(\bar x))$.
\end{proof}

Now, we are in position to state the central result of this subsection.

\begin{proposition}\label{prop:4_conditions_convex_polyhedral_case}
	Set $\mathbb Y:=\R^m$, let $D$ be convex and polyhedral, and fix $(\bar x,0)\in\gph\Phi$
	as well as $u\in\mathbb S_\mathbb X$ such that $\nabla g(\bar x)u\in\mathcal T_D(g(\bar x))$. 
	Then the following statements hold.
	\begin{enumerate}
		\item\label{item:2reg_char}
		 The 2-regularity condition \eqref{eq:2regularity_primal} 
		 is equivalent to the implication
		\begin{equation}\label{eq:2regularity_dual}
			\eqref{eq:2-reg_assumptions}
			\quad\Longrightarrow\quad
			y^*=0.
		\end{equation}
		\item\label{item:suff_cond_pseudo_reg} 
		Condition \eqref{eq:some_sufficient_condition_for_pseudo_regularity} is equivalent to
		\begin{equation}\label{eq:some_sufficient_condition_for_pseudo_regularity_explicit}
			\eqref{eq:2-reg_assumptions},\,
			\dual{z^*}{\nabla g(\bar x)u} = 0
			\quad\Longrightarrow\quad
			y^*=0.
		\end{equation}
		\item\label{item:suff_cond_pseudo_subreg_Gfrerer} 
		Gfrerer's condition \eqref{eq:Gfrerers_sufficient_condition_implicit}
		and condition \eqref{eq:CQ_pseudo_subregularity_polyhedral_II}
		are both equivalent to
		\begin{equation}\label{eq:Gfrerers_sufficient_condition_fully_explicit}
			\eqref{eq:2-reg_assumptions},\,
			w_s(u,0)\in\mathbf{T}(u)
			\quad\Longrightarrow\quad
			y^*=0.
		\end{equation}
	\end{enumerate}
\end{proposition}
\begin{proof}
	Let us start to prove~\ref{item:2reg_char}.
	If the 2-regularity condition holds at $\bar x$ in direction $u$, 
	then computing the polar cone on both sides of
	\eqref{eq:2regularity_primal} while respecting \cite[Exercises~3.4(d) and 3.5]{BertsekasNedicOzdaglar2003} 
	gives
	\[
		\{0\}
		=
		\mathcal N_D(g(\bar x))\cap\ker \nabla g(\bar x)^*
		\cap
		\{y^*\in\R^m\,|\,-\nabla^2\langle y^*,g\rangle(\bar x)(u)\in(\nabla g(\bar x)^{-1}\mathcal T_D(g(\bar x)))^\circ\}.
	\]
	Relying on \cite[Exercise~3.5]{BertsekasNedicOzdaglar2003} again 
	while taking convexity and polyhedrality of $D$ (and, thus, of $\mathcal N_D(g(\bar x))$) into account, we find
	\begin{equation}\label{eq:2regularity_dual_II}	
		\{0\}=
		\mathcal N_D(g(\bar x))\cap\ker \nabla g(\bar x)^*
		\cap
		\{y^*\in\R^m\,|\,-\nabla^2\langle y^*,g\rangle(\bar x)(u)\in\nabla g(\bar x)^*\mathcal N_D(g(\bar x))\}.
	\end{equation}
	Hence, \eqref{eq:2regularity_dual} holds.
	Conversely, if \eqref{eq:2regularity_dual} is valid, then \eqref{eq:2regularity_dual_II} holds as well.
	Computing the polar cone on both sides, we can exploit \cite[Exercises~3.4(d) and 3.5]{BertsekasNedicOzdaglar2003}
	once again in order to obtain
	\[
		\cl\left(
			\Im \nabla g(\bar x)+\nabla^2g(\bar x)[u,\nabla g(\bar x)^{-1}\mathcal T_D(g(\bar x))]-\mathcal T_D(g(\bar x))
		\right)
		=
		\R^m.
	\]
	Finally, one has to observe that the set within the closure operator is a convex polyhedral cone and, thus, closed
	in order to find validity of the 2-regularity condition at $\bar x$ in direction $u$.

	Statement~\ref{item:suff_cond_pseudo_reg} follows immediately from \cref{lem:some_properties_of_polyhedral_sets}\,\ref{item:normal_cones_to_polyhedral_sets}.
	
	Finally, let us turn to the proof of statement~\ref{item:suff_cond_pseudo_subreg_Gfrerer}.
	In order to show the equivalence between conditions
	\eqref{eq:Gfrerers_sufficient_condition_implicit} 
	and \eqref{eq:Gfrerers_sufficient_condition_fully_explicit},
	it suffices to prove that \eqref{eq:Gfrerers_sufficient_condition_implicit} is equivalent to
	\begin{equation}\label{eq:Gfrerers_sufficient_condition_explicit}
	\eqref{eq:2-reg_assumptions},\,\nabla^2\dual{y^*}{g}(\bar x)[u,u]=0,\,
	w_{\tilde s}(u,0)\in\mathcal T_D(g(\bar x))
	\quad
	\Longrightarrow
	\quad
	y^*=0,
	\end{equation}
	since the latter is equivalent to
	\eqref{eq:Gfrerers_sufficient_condition_fully_explicit}
	by \cref{lem:w_s_versus_w_tilde_s}.
	The maximization problem appearing in \eqref{eq:Gfrerers_sufficient_condition_implicit}
	is a linear program whose feasible set is a nonempty, convex polyhedral cone.
	Furthermore, $y^*\in\mathcal N_D(g(\xb))\cap\ker \nabla g(\xb)^*$ 
	is a maximizer if and only if
\begin{equation*}
	\begin{aligned}
	1/2\nabla^2g(\bar x)[u,u]
	&\in
	\mathcal N_{\mathcal N_D(g(\bar x))\cap\ker\nabla g(\bar x)^*}(y^*)
	=
	\left(\mathcal N_D(g(\bar x))\cap\ker\nabla g(\bar x)^*\right)^\circ\cap[y^*]^\perp
	\\
	&=
	\left(\mathcal T_D(g(\bar x))+\Im \nabla g(\bar x)\right)\cap[y^*]^\perp.
	\end{aligned}
\end{equation*}
Here, we made use of \cite[Exercise~3.4(d)]{BertsekasNedicOzdaglar2003} to compute
the polar cone of the appearing intersection, and the latter is a polyhedral cone
and, thus, closed.
This inclusion, in turn, is equivalent to the existence of $\tilde s\in\mathbb X$ such that
\begin{equation*}
	\nabla^2\dual{y^*}{g}(\bar x)[u,u]=0,\,
	w_{\tilde s}(u,0)\in\mathcal T_D(g(\bar x)),
\end{equation*}
showing the claimed equivalence between
\eqref{eq:Gfrerers_sufficient_condition_implicit}
and \eqref{eq:Gfrerers_sufficient_condition_explicit}
as $y^*\in\mathcal N_D(g(\bar x))\cap\ker\nabla g(\bar x)^*$ is
already included in \eqref{eq:2-reg_assumptions}.

Clearly, \eqref{eq:Gfrerers_sufficient_condition_fully_explicit} 
implies \eqref{eq:CQ_pseudo_subregularity_polyhedral_II}
by \cref{lem:some_properties_of_polyhedral_sets}\,\ref{item:normal_cones_to_polyhedral_sets}
and \cref{rem:refined_polyhedral_situation_normal_cone_relation},
so we only need to verify the converse implication.
Thus, let us prove the premise of \eqref{eq:CQ_pseudo_subregularity_polyhedral_II}
assuming that \eqref{eq:2-reg_assumptions} holds 
while there exists some $s \in \mathbb X$ with $w_s(u,0)\in\mathbf{T}(u)$.
Particularly, from these two we infer $\dual{z^*}{\nabla g(\bar x)u}=0$
with the help of \cref{lem:w_s_versus_w_tilde_s},
so the premise of \eqref{eq:some_sufficient_condition_for_pseudo_regularity_explicit} is valid.
Taking into account \cref{rem:comparison_estimate_pseudocoderivative_polyhedral}, 
this means that $u$ is a solution of the linear program
\hyperref[eq:LP_q]{($\widehat{\text{LP}}(0)$)}
where we used
\begin{equation}\label{eq:LP_q}\tag{$\widehat{\text{LP}}(q)$}
	\min\limits_s \{\nabla^2\langle y^*,g\rangle(\bar x)[u,s] \,|\, 
		\nabla g(\xb)s \in \mathcal T_D(g(\xb)) - q\}
\end{equation}
for some parameter $q\in\R^m$.

For arbitrary $q\in\R^m$, 
we claim that whether \eqref{eq:LP_q} has a solution
depends only on its feasibility since, for feasible problems, 
the issue of boundedness is independent of $q$.
This follows from \cite[Lemma~4]{BenkoGfrerer2018}, 
stating that, whenever \eqref{eq:LP_q} is feasible,
then it possesses a solution if and only if there does not exist $s\in\mathbb X$ 
satisfying $\nabla^2\langle y^*,g\rangle(\bar x)[u,s] < 0$ and
\begin{equation*}
	\nabla g(\xb)s 
	\in 
	(\mathcal T_D(g(\xb)) - q)^\infty 
	= 
	(\mathcal T_D(g(\xb)))^\infty
	=
	\mathcal T_D(g(\bar x)),
\end{equation*}
and these conditions are, indeed, independent of $q$.
Above, we have used \cite[Exercises~3.12 and~6.34(c)]{RockafellarWets1998}.
Since \hyperref[eq:LP_q]{($\widehat{\text{LP}}(0)$)} has a solution, 
\eqref{eq:LP_q} has a solution for each $q\in\R^m$ for which it is feasible.
Particularly, \cref{lem:w_s_versus_w_tilde_s} thus yields that 
\hyperref[eq:LP_q]{($\widehat{\text{LP}}(\bar q)$)}
has a solution $\bar s\in\mathbb X$ for $\bar q:=\nabla^2 g (\bar x)[u,u]$.

	Finally, we claim that $\bar s$ is also a solution of the (feasible) linear program
	\begin{align*}
		\min\limits_{\tilde s}\{\nabla^2\dual{y^*}{g}(\bar x)[u,\tilde s]\,|\,
			w_{\tilde s}(u,0)\in\mathbf T(u)\},
	\end{align*}
	whose feasible set equals the set $\widetilde A_1(y^*,u,0)$ from \eqref{eq:A_i_sets}.
	As explained just below \eqref{eq:A_i_sets}, 
	this will confirm the premise of \eqref{eq:CQ_pseudo_subregularity_polyhedral_II}
	and thus conclude the proof.
	Suppose that $\bar s$ is not a solution of this problem, i.e., there exists $\hat s$ with
	$\nabla^2\langle y^*,g\rangle(\bar x)[u,\hat s - \bar s] < 0$ and
	\begin{equation*}
		w_{\hat s}(u,0) \in \mathbf{T}(u) =  \mathcal T_D(g(\bar x))+\spa(\nabla g(\bar x)u).
	\end{equation*}
	Then $\hat s + \alpha u$ is a feasible point 
	of \hyperref[eq:LP_q]{($\widehat{\text{LP}}(\bar q)$)} for some $\alpha \in \R$, while
	\begin{equation*}
		\nabla^2\langle y^*,g\rangle(\bar x)[u,(\hat s + \alpha u) - \bar s]
		=
		\nabla^2\langle y^*,g\rangle(\bar x)[u,\hat s - \bar s]
		<
		0
	\end{equation*}
	follows from $\nabla^2\langle y^*,g\rangle(\bar x)[u,u] = 0$ 
	which holds by \cref{lem:w_s_versus_w_tilde_s}.
	The latter, however,
	means that $\bar s$ is not optimal for \hyperref[eq:LP_q]{($\widehat{\text{LP}}(\bar q)$)}
	 - a contradiction.
\end{proof}

	Let us mention that the first assertion of \cref{prop:4_conditions_convex_polyhedral_case}
	generalizes \cite[Proposition~2]{GfrererOutrata2016c}.

As a corollary of \cref{prop:4_conditions_convex_polyhedral_case}, 
we now can easily interrelate the different sufficient conditions
for pseudo-(sub)regularity.

\begin{corollary}\label{cor:4_conditions_convex_polyhedral_case}
	Set $\mathbb Y:=\R^m$, let $D$ be convex and polyhedral, and fix $(\bar x,0)\in\gph\Phi$
	as well as $u\in\mathbb S_{\mathbb X}$ 
	such that $\nabla g(\bar x)u\in\mathcal T_D(g(\bar x))$. 
	Then the following implications hold:
	\[
		\eqref{eq:2regularity_primal}
		\quad\Longrightarrow\quad
		\eqref{eq:some_sufficient_condition_for_pseudo_regularity}
		\quad\Longrightarrow\quad
		\eqref{eq:Gfrerers_sufficient_condition_implicit}
		\quad\Longleftrightarrow\quad
		\eqref{eq:CQ_pseudo_subregularity_polyhedral_II}.
	\]
	Particularly, \eqref{eq:2regularity_primal} implies that
	$\Phi$ is metrically pseudo-regular of order 2 at $(\bar x,0)$ in direction $(u,0)$.
	Moreover, if there exists $s \in \mathbb X$ with $w_s(u,0)\in\mathbf{T}(u)$, 
	all four conditions are equivalent.
\end{corollary}
\begin{proof}
	The first implication and the equivalence are immediately clear 
	by \cref{prop:4_conditions_convex_polyhedral_case}.
	In order to show the second implication, we first make use of
	\cref{prop:4_conditions_convex_polyhedral_case} in order to see
	that it suffices to verify that
	\eqref{eq:some_sufficient_condition_for_pseudo_regularity_explicit}
	implies \eqref{eq:Gfrerers_sufficient_condition_fully_explicit}.
	This, however, is clear since \eqref{eq:2-reg_assumptions} and
	the existence of $s \in \mathbb X$ such that $w_s(u,0)\in\mathbf T(u)$
	imply $\nabla^2\dual{y^*}{g}(\bar x)[u,u]=0$, see \cref{lem:w_s_versus_w_tilde_s},
	and $\dual{z^*}{\nabla g(\bar x)u}=0$ follows by \eqref{eq:2-reg_assumptions}.

		The fact that \eqref{eq:2regularity_primal} is sufficient for directional
		pseudo-regularity of $\Phi$ now follows 
		from \cref{cor:sufficient_condition_pseudo_subregularity}.
		The final statement is obvious from \cref{prop:4_conditions_convex_polyhedral_case}.
\end{proof}

The following example shows that 
our sufficient condition \eqref{eq:some_sufficient_condition_for_pseudo_regularity} 
for directional pseudo-regularity is strictly milder than directional 2-regularity
from \eqref{eq:2regularity_primal}.

\begin{example}\label{ex:2regularity_stronger}
	Let $g\colon\R \to \R^2$ and $D_i \subset \R^2$, $i=1,2$, be given by
	$g(x) := (x,-x^2)$, $x\in\R$,  and 
	\[
		D_1:=\R \times \R_+,
		\qquad
		D_2:=\R_-\times \R_+.
	\]
	Observe that $D_i$ is a convex polyhedral set for $i=1,2$.
	We consider the constraint mappings $\Phi_i\colon\R\tto\R^2$ given by $\Phi_i(x):=g(x)-D_i$,
	$x\in\R$, for $i=1,2$ and fix $\bar x:=0$ and $u:=-1$.
	Note that $(\bar x,0)\in\gph\Phi_i$ for $i=1,2$.
	 
	Let us start with the investigation of the mapping $\Phi_1$.
	Due to $\mathcal N_{D_1}(g(\xb))=\{0\}\times\R_-$ and
	\[
		\nabla g(\bar x)^*y^*=y_1^*,
		\qquad
		\nabla^2\langle y^*,g\rangle(\bar x)(u) + \nabla g(\bar x)^* z^* = 2y_2^* + z_1^*,
	\]
	one can easily check that
	\eqref{eq:2regularity_dual} and
	\eqref{eq:some_sufficient_condition_for_pseudo_regularity_explicit}
	are both satisfied.
	Consequently, due to \cref{prop:4_conditions_convex_polyhedral_case},
	\eqref{eq:some_sufficient_condition_for_pseudo_regularity} and
	\eqref{eq:2regularity_primal} hold in parallel.
	
	Let us now consider the mapping $\Phi_2$.
	Clearly, \eqref{eq:some_sufficient_condition_for_pseudo_regularity_explicit}
	remains valid 
	since the appearing variable $z^*$ has to be chosen from the set
	$\mathcal N_{D_2}(g(\xb)) \cap [\nabla g(\bar x)u]^\perp=\{0\} \times \R_-$.
	Hence, due to \cref{prop:4_conditions_convex_polyhedral_case},
	\eqref{eq:some_sufficient_condition_for_pseudo_regularity} holds
	(and, thus, pseudo-regularity of order $2$ of $\Phi_2$ in direction $u$ at $(\bar x,0)$).
	However, we have $\mathcal N_{D_2}(g(\xb))=\R_+\times\R_-$,
	so that choosing $y^*:=(0,-1)$ and $z^*:=(2,0)$ yields a violation of 
	\eqref{eq:2regularity_dual} in this situation.
	Consulting \cref{prop:4_conditions_convex_polyhedral_case} once again,
	\eqref{eq:2regularity_primal} is violated as well.
	Let us also note that, for each $s\in\R$, 
	we have 
	\[
		w_s(u,0)
		=
		\nabla g(\bar x)s+1/2\nabla^2g(\bar x)[u,u]
		=
		(s,-1)
		\notin
		\R_-\times\R_+
		=
		\mathcal T_{D_2}(g(\bar x)),
	\]
	see \cref{cor:4_conditions_convex_polyhedral_case}.
	Hence, \eqref{eq:some_sufficient_condition_for_pseudo_regularity}
	is strictly milder than \eqref{eq:2regularity_primal}.
\end{example}

Let us take a closer look at the particular situation where $D:=\{0\}$.
\begin{remark}\label{rem:D=0}
	Set $\mathbb Y:=\R^m$, $D:=\{0\}$, 
	and fix $(\bar x,0)\in\gph\Phi$ as well as $u\in\mathbb S_{\mathbb X}$
	such that $u\in\ker\nabla g(\bar x)$.
	Let us consider the sufficient conditions for directional metric 
	pseudo-(sub)regularity discussed
	in \cref{prop:4_conditions_convex_polyhedral_case}.
	The constraint qualification \eqref{eq:2regularity_dual} obviously reduces to 
	\begin{equation}\label{eq:CQ_for_D=0}
		\nabla g(\bar x)^*y^*=0,\,
		\nabla^2\langle y^*,g\rangle(\bar x)(u)+\nabla g(\bar x)^*z^*=0
		\quad
		\Longrightarrow
		\quad
		y^*=0,
	\end{equation}
	and the latter is equivalent to the 2-regularity condition 
	\eqref{eq:2regularity_primal} at $\bar x$ in direction $u$
	by \cref{prop:4_conditions_convex_polyhedral_case}.
	One can easily check that 
	\eqref{eq:some_sufficient_condition_for_pseudo_regularity}
	also reduces to \eqref{eq:CQ_for_D=0}.
	Furthermore, due to \cref{prop:4_conditions_convex_polyhedral_case},
	\eqref{eq:CQ_pseudo_subregularity_polyhedral_II} and
	\eqref{eq:Gfrerers_sufficient_condition_implicit}
	reduce to
	\begin{equation*}
		\nabla g(\bar x)^*y^*=0,\,
		\nabla^2\langle y^*,g\rangle(\bar x)(u)+\nabla g(\bar x)^*z^*=0,\,
		w_s(u,0)=0
		\quad
		\Longrightarrow
		\quad
		y^*=0,
	\end{equation*}
	and the latter is strictly milder than \eqref{eq:CQ_for_D=0}
	as we will illustrate in \cref{ex:some_trivial_example} below.
	
	To close the remark, let us mention that whenever \eqref{eq:CQ_for_D=0}
	has to hold for all $u\in\mathbb S_{\mathbb X}$ (this implies metric pseudo-subregularity
	of order 2 of $\Phi$ at $(\bar x,0)$ for all unit directions),
	then either $\nabla g(\bar x)$ is surjective or the zero operator,
	see \cite[Remark~2.1]{FischerIzmailovJelitte2023},
	i.e., this situation is rather special.
	We believe, however, that this is mainly because $D$ is trivial
	and partially due to the precise definition of 2-regularity.
	Let us point the interested reader to \cite[Example~2]{Gfrerer2014a},
	which suggests that metric pseudo-subregularity of order $2$ in all unit directions
	might be a reasonable assumption.
\end{remark}

The following example,
which has been motivated by \cref{rem:D=0}, 
indicates that \eqref{eq:Gfrerers_sufficient_condition_implicit}
is strictly milder than \eqref{eq:some_sufficient_condition_for_pseudo_regularity}.
\begin{example}\label{ex:some_trivial_example}
	Let $g\colon\R^2\to\R^3$ and $D\subset\R^3$
	be given by $g(x):=(x_1^2,x_2^2,x_1x_2)$, $x\in\R^2$,
	and $D:=\{0\}$. We consider the point $\bar x:=0$.
	As $\nabla g(\bar x)$ vanishes while we have $\mathcal T_D(g(\bar x))=\{0\}$, 
	each direction $u\in\mathbb S_{\R^2}$
	satisfies $\nabla g(\bar x)u\in\mathcal T_D(g(\bar x))$,
	and we pick any such $u$.
	Due to \cref{rem:D=0}, \eqref{eq:some_sufficient_condition_for_pseudo_regularity}
	and \eqref{eq:2regularity_primal} reduce to
	\[
		y_1^*(2u_1,0)+y_2^*(0,2u_2)+y_3^*(u_2,u_1)=(0,0)
		\quad\Longrightarrow\quad
		y^*=0,
	\]
	and since three vectors in $\R^2$ are always linearly dependent,
	this condition is trivially violated.
	On the other hand, 
	\eqref{eq:CQ_pseudo_subregularity_polyhedral_II} and 
	\eqref{eq:Gfrerers_sufficient_condition_implicit} can be stated as
	\[
		\left.
		\begin{aligned}
		&y_1^*(2u_1,0)+y_2^*(0,2u_2)+y_3^*(u_2,u_1)=(0,0),\\
		&(u_1^2,u_2^2,u_1u_2)=(0,0,0)
		\end{aligned}
		\right\}
		\quad\Longrightarrow\quad
		y^*=0,
	\]
	and this condition holds as the premise regarding $u$ cannot be satisfied
	by any $u\in\mathbb S_{\R^2}$.
\end{example}

	We close this subsection with some more general remarks about (directional)
	2-regularity and Gfrerer's sufficient condition for metric pseudo-(sub)regularity
	from \cite[Theorem~2]{Gfrerer2014a}.

In this subsection, for simplicity, we restricted ourselves to the convex polyhedral case, 
but neither our approach nor the other results are limited to this case.
The original definition of directional 2-regularity
in \cite{ArutyunovAvakovIzmailov2008} is stated for merely convex sets $D$
(no polyhedrality is assumed in the latter paper),
but involves the radial cone to $D$ which is not necessarily closed for curved sets $D$.
Interestingly, \cite[Example~2]{Gfrerer2014a}, already mentioned in \cref{rem:D=0}, 
provides a mapping which is metrically pseudo-regular
of order $2$ in every direction $(u,0)$ with $u \neq 0$,
particularly metrically pseudo-subregular
of order $2$ in every unit direction,
but the 2-regularity condition is violated for every direction;
the chosen set $D$ in this example is the Euclidean unit ball
in $\R^2$ which is not polyhedral.

Let us mention that
\cite[Theorem~2]{Gfrerer2014a} is stated
in the general polyhedral case (no convexity is assumed), and it
yields the existence of several elements $s \in \mathbb X$
corresponding to the active components of the set $D$.
Looking into the proof of \cref{The : NCgen_2}, 
it seems like we could get a similar result with only minor adjustments,
but since we do not need such a result here, we did not develop this approach 
for the purpose of brevity.

Let us also note that the conditions
from statements~\ref{item:trivial_SC_for_MPR} and~\ref{item:SC_for_MPR}
of \cref{cor:sufficient_condition_pseudo_subregularity} are not
covered by \cite{ArutyunovAvakovIzmailov2008}
(since $D$ does not need to be convex for our findings)
or by \cite[Theorem~2]{Gfrerer2014a} (since $D$ does not need
to be polyhedral).

Finally, let us point out that the concept of 2-regularity is useful for the design and the
convergence analysis of Newton-type methods, aiming to solve smooth and nonsmooth equations, 
see e.g.\ \cite{FischerIzmailovJelitte2021,IzmailovKurennoySolodov2018} 
and the references therein.

\section{Directional asymptotic stationarity in nonsmooth optimization}\label{sec:directional_asymptotic_tools}

This section is devoted to directional asymptotic stationarity conditions
and related results.
It contains the foundation of our research, 
\cref{thm:higher_order_directional_asymptotic_stationarity},
which also motivates our considerations in \cref{sec:directional_asymptotic_regularity}.

For a locally Lipschitz continuous function $\varphi\colon\mathbb X\to\R$,
a set-valued mapping $\Phi\colon\mathbb X\tto\mathbb Y$ with a closed graph,
and $\bar y\in\Im\Phi$, we investigate the rather abstract optimization problem
\begin{equation}\label{eq:nonsmooth_problem}\tag{P}
	\min\limits_x\{\varphi(x)\,|\,\bar y\in\Phi(x)\}.
\end{equation}
Throughout the section, the feasible set of \eqref{eq:nonsmooth_problem} 
will be denoted by $\mathcal F\subset\mathbb X$. 
Clearly, we have $\mathcal F\neq\emptyset$ from $\bar y\in\Im\Phi$.
Note that the model \eqref{eq:nonsmooth_problem} covers numerous 
classes of optimization problems from the literature including
standard nonlinear problems, problems with geometric 
(particularly, disjunctive or conic) constraints, problems with
(quasi) variational inequality constraints, and bilevel optimization problems. 
Furthermore, we would like to mention that choosing $\bar y:=0$ would not be restrictive
since one could simply consider $\widetilde\Phi\colon\mathbb X\tto\mathbb Y$ given by
$\widetilde{\Phi}(x):=\Phi(x)-\bar y$, $x\in\mathbb X$, in the case where $\bar y$ does not
vanish.
Optimality conditions and constraint qualifications for problems of this type can be
found, e.g., in \cite{Gfrerer2013,Mehlitz2020a,Mordukhovich2006,YeYe1997}.
A standard notion of stationarity, which applies to \eqref{eq:nonsmooth_problem} and
is based on the tools of limiting variational analysis, is the
one of M-stationarity.

\begin{definition}\label{def:M_stationarity}
	A feasible point $\bar x\in\mathcal F$ of \eqref{eq:nonsmooth_problem} 
	is called \emph{M-stationary} whenever there is a multiplier $\lambda\in\mathbb Y$ such that
	\[
		0\in\partial\varphi(\bar x)+D^*\Phi(\bar x,\bar y)(\lambda).
	\]
\end{definition}

Later in \cref{lem:directional_M_stationarity_via_metric_subregularity}, we will show that
directional metric subregularity of $\Phi$ serves as a constraint qualification for M-stationarity.
In the following lemma, whose proof is analogous to the one of
\cite[Lemma~3.1]{BaiYe2021}, we point out that directional metric subregularity
of $\Phi$ implies that penalizing the constraint in \eqref{eq:nonsmooth_problem}
with the aid of the distance function yields a \emph{directionally} exact penalty function.
\begin{lemma}\label{lem:directional_exact_penalization}
	Let $\bar x\in\mathcal F$ be a local minimizer of \eqref{eq:nonsmooth_problem},
	and assume that $\Phi$ is metrically subregular at $(\bar x,\bar y)$
	in direction $u\in\mathbb S_{\mathbb X}$. 
	Then there are constants $\varepsilon>0$, $\delta>0$, and $C>0$ such
	that $\bar x$ is a global minimizer of
	\begin{equation}\label{eq:directionally_penalized_problem}
		\min\limits_x\{
			\varphi(x)+C\dist(\bar y,\Phi(x))
			\,|\,
			x\in\bar x+\mathbb B_{\varepsilon,\delta}(u)
		\}.
	\end{equation}
\end{lemma} 

Let us note that this result refines well-known findings about classical 
exact penalization in the presence of metric subregularity, see e.g.\ \cite{Burke1991,Clarke1983,KlatteKummer2002}.

\subsection{Approaching mixed-order stationarity conditions}\label{sec:Stat_cond_order_gamma}

	To start, let us introduce a quite general notion of critical directions associated 
	with \eqref{eq:nonsmooth_problem}.
\begin{definition}\label{def:critical_direction}
	For some feasible point $\bar x\in\mathcal F$ and a pair $(\gamma_0,\gamma)\in\R\times\R$ 
	such that $\gamma_0\geq 1$ as well as $\gamma \geq 1$, 
	a direction $u\in \mathbb X$ is called 
	\emph{critical of order $(\gamma_0,\gamma)$} for \eqref{eq:nonsmooth_problem} at $\bar x$ 
	whenever there are sequences $\{u_k\}_{k\in\N}\subset\mathbb X$, $\{\alpha_k\}_{k\in\N} \subset \R$,
	$\{v_k\}_{k\in\N}\subset\mathbb Y$, and $\{t_k\}_{k\in\N}\subset\R_+$ 
	satisfying $u_k\to u$, $t_k\downarrow 0$, $\alpha_k \to 0$, $v_k\to 0$, and, for all $k\in\N$,
	\begin{equation}\label{eq:Crit_Dir}
		(\bar x+t_ku_k,\varphi(\bar x) + (t_k \norm{u_k})^{\gamma_0} \alpha_k) \in \epi \varphi,
		\qquad
		(\bar x+t_ku_k,\bar y+(t_k \norm{u_k})^{\gamma} v_k)\in\gph\Phi.
	\end{equation}
	If $(\gamma_0,\gamma):=(1,1)$, we simply call $u$ a \emph{critical direction} for \eqref{eq:nonsmooth_problem} at $\bar x$.
\end{definition}

Clearly, $u:=0$ is critical of every order. 
Moreover, the set of all critical directions of any fixed order is a cone.
The most standard case $(\gamma_0,\gamma):=(1,1)$ corresponds to \cite[Definition~5]{Gfrerer2013}.
If $\varphi$ is directionally differentiable at $\bar x$, it is easily seen that
$u\in \mathbb X$ is critical for \eqref{eq:nonsmooth_problem} at $\bar x$ if and only if
$\varphi'(\bar x;u)\leq 0$ and $u\in\ker D\Phi(\bar x,\bar y)$,
see \cite[Proposition~3.5]{Shapiro1990} as well.
Let us note that whenever $\bar x\in\mathcal F$ is a feasible point of \eqref{eq:nonsmooth_problem}
such that no critical direction for \eqref{eq:nonsmooth_problem} at $\bar x$ exists, then $\bar x$
is a strict local minimizer of \eqref{eq:nonsmooth_problem}.
Conversely, there may exist strict local minimizers of \eqref{eq:nonsmooth_problem} such that
a critical direction for \eqref{eq:nonsmooth_problem} at this point exists.

While in this paper, we will not go beyond the case $\gamma_0 := 1$ (the case $\gamma_0:=2$ is briefly mentioned in \cref{Lem:Crit_dir_constraints}),
the situation $\gamma > 1$ (particularly $\gamma := 2$) will be very important.
For $\gamma_0 := 1$ and arbitrary $\gamma\geq 1$, 
a critical direction $u \in \mathbb X$ still satisfies
$\dr \varphi (\bar x)(u) \leq 0$ and $u\in\ker D_{\gamma}\Phi(\bar x,\bar y)$, and
the converse is valid whenever $\varphi$ is continuously differentiable at $\bar x$.
In the next lemma, we show that if
$\Phi$ is metrically pseudo-subregular of order $\gamma$ at $(\bar x,\bar y)$,
then $u$ is actually critical of order $(1,\gamma')$ for each $\gamma' \geq 1$.

\begin{lemma}\label{lem: crit_dir+pseudo_subreg}
	Fix a feasible point $\bar x\in\mathcal F$ of \eqref{eq:nonsmooth_problem}, $\gamma \geq 1$,
	and a critical direction $u \in \mathbb{X}$ of order $(1,\gamma)$
	for \eqref{eq:nonsmooth_problem} at $\bar x$.
	If $\Phi$ is metrically pseudo-subregular of order $\gamma$
	in direction $u$ at $(\bar x,\bar y)$, then $u$ is critical
	of order $(1,\gamma')$ for \eqref{eq:nonsmooth_problem} at $\bar x$ for each $\gamma' \geq 1$.
\end{lemma}
\begin{proof}
	Let $\{u_k\}_{k\in\N}\subset\mathbb X$, $\{\alpha_k\}_{k\in\N} \subset \R$,
	$\{v_k\}_{k\in\N}\subset\mathbb Y$, and $\{t_k\}_{k\in\N}\subset\R_+$
	be sequences satisfying $u_k\to u$, $t_k\downarrow 0$, $\alpha_k \to 0$, $v_k\to 0$, 
	as well as \eqref{eq:Crit_Dir} for all $k\in\N$.
	By metric pseudo-subregularity of order $\gamma$ of $\Phi$ at $(\bar x,\bar y)$, 
	there is a constant $\kappa>0$ such that, for sufficiently large $k\in\N$, 
	we get the existence of $\tilde x_k \in \Phi^{-1}(\yb)$ with
	\begin{equation*}
		\norm{\tilde x_k - (\bar x+t_ku_k)}  
		\leq
		\kappa \frac{\dist(\bar y,\Phi(\bar x+t_ku_k))}{(t_k \norm{u_k})^{\gamma - 1}} 
		\leq
		\kappa \frac{(t_k \norm{u_k})^{\gamma} \norm{v_k}}{(t_k \norm{u_k})^{\gamma - 1}} 
		=
		\kappa t_k \norm{u_k} \norm{v_k}.
	\end{equation*}
	Particularly, we find
	$\nnorm{(\tilde x_k - \bar x)/t_k - u_k}\to 0$
	from $u_k\to u$ and $v_k \to 0$.
	Moreover, Lipschitzianity of $\varphi$ yields
	\begin{align*}
        \varphi(\tilde x_k)
		& \leq
		\varphi(\tilde x_k) - \varphi(\bar x + t_k u_k)
		+ \varphi(\bar x) + t_k \norm{u_k} \alpha_k
		\\
		&\leq
		L\kappa t_k \norm{u_k} \norm{v_k}
		+ \varphi(\bar x) + t_k \norm{u_k} \alpha_k
		=
		\varphi(\bar x) + t_k \norm{u_k} (L\kappa \norm{v_k} + \alpha_k)
	\end{align*}
	for some constant $L > 0$ and sufficiently large $k \in \N$.
	Thus, setting $\tilde u_k := (\tilde x_k - \bar x)/t_k$,
	$\tilde \alpha_k := L\kappa \norm{v_k} + \alpha_k$,
	and $\tilde t_k := t_k$ for large enough $k \in \N$
	yields $\tilde u_k\to u$, $\tilde t_k\downarrow 0$, $\tilde\alpha_k\to 0$, 
	as well as 
	\[
		(\bar x+\tilde t_k\tilde u_k,\varphi(\bar x) + \tilde t_k \nnorm{\tilde u_k} \tilde\alpha_k) \in \epi \varphi,
		\qquad
		(\bar x+\tilde t_k\tilde u_k,\bar y)\in\gph\Phi
	\]
	for large enough $k\in\N$,
	and so $u$ is critical of order $(1,\gamma')$ for \eqref{eq:nonsmooth_problem} at $\bar x$
	for each $\gamma'\geq 1$.
\end{proof}

The following result, inspired by and based on \cite[Proposition 2]{Gfrerer2014a},
provides an important interpretation of the notion from \cref{def:critical_direction} 
in terms of the so-called \emph{epigraphical} mapping
$M_0\colon\mathbb X \tto \R$ associated with $\varphi$ and given by $M_0(x) := \varphi(x) + \R_+$, $x\in\mathbb X$.
The proof follows simply from the fact that $\gph M_0 = \epi \varphi$
together with \cref{Rem:Strong_pseudo_subreg}.

\begin{proposition}\label{Pro : critical_direction_interpretation}
 	Given a feasible point $\bar x\in\mathcal F$ and a pair $(\gamma_0,\gamma)\in\R\times\R$ 
	such that $\gamma_0\geq 1$ as well as $\gamma \geq 1$,
	a direction $u\in\mathbb S_{\mathbb X}$ is critical of order $(\gamma_0,\gamma)$
	for \eqref{eq:nonsmooth_problem} at $\bar x$
	if and only if there exist sequences $\{u_k\}_{k\in\N}\subset\mathbb X$ and $\{t_k\}_{k\in\N}\subset\R_+$ such that $u_k \to u$, $t_k\downarrow 0$, and
	\begin{equation}\label{crit_dir_equiv}
		\frac{\dist(\varphi(\bar x),M_0(\bar x + t_k u_k))}{(t_k \norm{u_k})^{\gamma_0}} \to 0,
		\qquad
		\frac{\dist(\bar y,\Phi(\bar x + t_k u_k))}{(t_k \norm{u_k})^\gamma} \to 0.
	\end{equation}
	Moreover, if $\gamma_0 = \gamma$, this is further equivalent to the condition
	\begin{equation*}
		u\in\ker D_{\gamma}M(\bar x, (\varphi(\bar x), \bar y))
	\end{equation*}
	for the mapping $M\colon\mathbb X \tto \R \times \mathbb Y$ 
	given by $M(x):=M_0(x)\times\Phi(x)$, $x\in\mathbb X$.
\end{proposition}

Interestingly, Gfrerer used the conditions \eqref{crit_dir_equiv} as a basis
of his optimality conditions in \cite[Proposition 2]{Gfrerer2014a},
but he did not notice, or at least did not mention,
that these conditions actually provide
a natural extension of his own notion of a critical direction from \cite[Definition~5]{Gfrerer2013}.
This observation enables us to formulate an extension of the common pattern
``for every critical direction there is a multiplier satisfying an FJ-type optimality condition'' in \cref{cor: F-J} below.

\begin{remark}
Gfrerer recognized the importance of considering Cartesian product mappings
$M\colon\mathbb X\to\mathbb Y_0\times\mathbb Y_1\times\ldots\times\mathbb Y_s$, given by
\[
	\forall x\in\mathbb X\colon\quad
	M(x):=M_0(x)\times M_1(x)\times\ldots\times M_s(x)
\]
for the component maps $M_i\colon\mathbb X\to\mathbb Y_i$, $i=0,1,\ldots,s$,
and Euclidean spaces $\mathbb Y_0,\mathbb Y_1,\ldots,\mathbb Y_s$,
and to allow different orders $\gamma_i\geq 1$ of pseudo-(sub)regularity of these component mappings, 
see \cite[Definition~1]{Gfrerer2014a}.
In the same manner, he defined his pseudo-coderivative \cite[Definition~2]{Gfrerer2014a}.
This was essential for his approach to optimality conditions.
For brevity of presentation, we avoid these definitions and bypass
explicitly using these notions by applying
\cite[Proposition~2]{Gfrerer2014a}
in combination with the sufficient conditions
for pseudo-subregularity from 
\cite[Theorem~1(2)]{Gfrerer2014a} to prove \cref{cor: F-J}.%
\end{remark}

\begin{corollary}\label{cor: F-J}
 Let $\bar x\in\mathcal F$ be a local minimizer of \eqref{eq:nonsmooth_problem} and
 let $u\in\mathbb S_{\mathbb X}$ be a critical direction of order $(1,\gamma)$ 
 for \eqref{eq:nonsmooth_problem} at $\bar x$ with $\gamma \geq 1$.
 Then there exist multipliers $(0,0) \neq (\alpha^*, \lambda) \in \R_+ \times \mathbb Y$ satisfying
 \begin{equation*}
	0
	\in
	\alpha^* \partial \varphi(\bar x;u)
	+
	D^\ast_{\gamma}\Phi((\bar x, \bar y);(u,0))(\lambda).
 \end{equation*}
 If the sufficient condition \eqref{eq:FOSCMS_gamma} for metric pseudo-subregularity of order $\gamma$ of $\Phi$
 in direction $u$ at $(\bar x, \bar y)$ holds,
 then the above condition holds with $\alpha^* := 1$.
\end{corollary}
\begin{proof}
	Applying \cref{Pro : critical_direction_interpretation} and then 
	\cite[Proposition 2 and Theorem~1(2)]{Gfrerer2014a}
	yields an element $z^* \in \mathbb X$ and sequences $\{t_k\}_{k\in\N}\subset\R_+$,
	$\{u_k\}_{k\in\N},\{x_k^*\}_{k\in\N}\subset\mathbb X$,
	$\{\alpha_k\}_{k\in\N},\{\alpha_k^*\}_{k\in\N}\subset\R$,
	and $\{v_k\}_{k\in\N},\{\lambda_k\}_{k\in\N}\subset\mathbb Y$
	satisfying (among other things) $t_k\downarrow 0$, $u_k\to u$, $\alpha_k\to 0$, $t_k^{1-\gamma}v_k\to 0$, as well as
	$x_k^*\to 0$, such that, for each $k\in\N$, $\norm{(\alpha_k^*,\lambda_k)}=1$ and 
	\begin{equation}\label{eq:what_we_get_from_Gfrerer}
		(x_k^*,-\alpha_k^*, -\lambda_k/(t_k\norm{u_k})^{\gamma-1})
		\in
		\widehat{\mathcal N}_{\gph M^{z^*}}(\bar x +t_ku_k,\varphi(\bar x) + t_k \alpha_k,\bar y + t_k^{\gamma} \tilde v_k),
	\end{equation}
	where we used $\tilde v_k:=t_k^{1-\gamma}v_k$ as well as the mapping
	$M^{z^*}\colon\mathbb X \tto \R \times \mathbb Y$ defined by
	$M^{z^*}(x):=M^{z^*}_0(x)\times\Phi(x)$, $x\in\mathbb X$,
	with the perturbed epigraphical mapping $M^{z^*}_0\colon\mathbb X\tto\R$ given by
	\[
		\forall x\in\mathbb X\colon\quad
		M^{z^*}_0(x) := \varphi(x) +  \vert \innerprod{z^*}{x - \bar x} \vert^{3} + \R_+.
	\]
	Note that we have $\gph M^{z^*} = \big( \gph M^{z^*}_0 \times \mathbb Y \big) \cap \mathcal P(\gph \Phi \times \R)$,
	where the permutation mapping $\mathcal P\colon \mathbb X \times \mathbb Y \times \R \to \mathbb X \times \R \times \mathbb Y$
	just swaps the last two components.
	After replacing the regular by the larger limiting normal cone in \eqref{eq:what_we_get_from_Gfrerer}
	and noting that $x\mapsto\varphi(x)+|\innerprod{z^*}{x-\bar x}|^3$ is locally Lipschitzian,
	we can apply the intersection rule for limiting normals from
	\cite[Theorem~6.42]{RockafellarWets1998}.
	The latter yields, for each $k\in\N$,
	$x_{k,1}^*,x_{k,2}^*\in\mathbb X$ with $x_k^* = x_{k,1}^* + x_{k,2}^*$ and
	\begin{align*}
		(x_{k,1}^*,-\alpha_k^*) & \in  \mathcal N_{\gph M^{z^*}_0 }(\bar x +t_ku_k,\varphi(\bar x) + t_k \alpha_k),\\
		(x_{k,2}^*,-\lambda_k/(t_k\norm{u_k})^{\gamma-1}) & \in  \mathcal N_{\gph \Phi}(\bar x +t_ku_k,\bar y + t_k^{\gamma} \tilde v_k).
	\end{align*}
	Now, local Lipschitzness of $x\mapsto\varphi(x)+|\innerprod{z^*}{x-\bar x}|^3$ together with boundedness of $\{\alpha_k^*\}_{k\in\N}$
	implies boundedness of $\{x_{k,1}^*\}_{k\in\N}$.
	This, in turn, gives boundedness of $\{x_{k,2}^*\}_{k\in\N}$.
	Since $\{\lambda_k\}_{k\in\N}$ is also bounded,
	taking the limit along a suitable subsequence 
	yields some $x^*\in\mathbb X$, $\alpha^*\in\R$, and $\lambda\in\mathbb Y$ satisfying $(\alpha^*,\lambda)\neq(0,0)$ as well as
	\[
		x^* \in D^* M^{z^*}_0 ((\bar x, \varphi(\bar x));(u,0))(\alpha^*),
		\quad
		-x^* \in D^\ast_{\gamma}\Phi((\bar x, \bar y);(u,0))(\lambda).
	\]
	Here, we used the robustness of the directional limiting coderivative,
	see \cref{lem:robustness_directional_limiting_normals}, 
	as well as \cref{lem:pseudo_coderivative_via_limiting_normals}.
	Taking into account that $x\mapsto\vert \innerprod{z^*}{x - \bar x} \vert^{3}$ is smooth with its gradient vanishing at $\bar x$
	and using \cref{lem:coderivatives_constraint_maps}\,\ref{item:constraint_maps_directional_limiting_coderivative}
	as well as \cite[Proposition~5.1]{BenkoGfrererOutrata2019}, we get $\alpha^* \geq 0$ and
	$D^*  M^{z^*}_0 ((\bar x, \varphi(\bar x));(u,0))(\alpha^*) \subset \alpha^* \partial \varphi(\bar x;u)$.
	This proves the first statement.

	Finally, \eqref{eq:FOSCMS_gamma} clearly implies $\alpha^* > 0$, and by rescaling, we can set $\alpha^* := 1$.
\end{proof}

\subsection{Mixed-order and asymptotic stationarity conditions}\label{sec:stationarity_via_coderivative}

The following result provides asymptotic necessary optimality conditions for \eqref{eq:nonsmooth_problem}
which hold in the absence of constraint qualifications.
The derived conditions depend on a certain \emph{order} $\gamma\geq 1$.
Furthermore, our result specifies how
the asymptotic case \ref{item:irregular_case_pseudo_MSt}
can be ruled out by
metric pseudo-subregularity of $\Phi$ of order $\gamma$
at the reference point.%

\begin{theorem}\label{thm:higher_order_directional_asymptotic_stationarity}
	Let $\bar x\in\mathcal F$ be a local minimizer of \eqref{eq:nonsmooth_problem} 
	and consider $\gamma \geq 1$.
	Then one of the following conditions holds.
	\begin{enumerate}
	  \item\label{item:trivial_MSt}
	  	The point $\bar x$ is M-stationary for \eqref{eq:nonsmooth_problem}.
	  \item\label{item:pseudo_MSt_zero_mult}%
		There exists a critical direction $u\in\mathbb S_{\mathbb X}$ 
		for \eqref{eq:nonsmooth_problem} at $\bar x$ such that
			\begin{equation}\label{eq:Pseudo-M-Stat_Gferer}
				0\in \partial\varphi(\bar x;u)+\widetilde{D}^*_\gamma\Phi((\bar x,\bar y);(u,0))(0).
			\end{equation}
	  \item\label{item:pseudo_dir_MSt}%
		There exist a critical direction $u\in\mathbb S_{\mathbb X}$ 
		for \eqref{eq:nonsmooth_problem} at $\bar x$,
		a nonvanishing multiplier $\lambda \in \mathbb Y$,
		and $\alpha \geq 0$
		such that, for $v:= \alpha \lambda$, we have
		\begin{equation}\label{eq:Pseudo-M-stat}	
		  0 \in \partial\varphi(\bar x;u) + D^\ast_{\gamma} \Phi((\xb,\yb);(u,v))(\lambda).
		\end{equation}
	  \item\label{item:irregular_case_pseudo_MSt}
		There exist a critical direction $u\in\mathbb S_{\mathbb X}$
		of order $(\gamma_0,\gamma)$ for \eqref{eq:nonsmooth_problem} at $\bar x$ for each $\gamma_0\geq 0$, 
		some $y^*\in\mathbb Y$, and sequences
		$\{x_k\}_{k\in\N},\{\eta_k\}_{k\in\N}\subset\mathbb X$ as well as
	$\{y_k\}_{k\in\N}\subset\mathbb Y$ such that
	$x_k\notin\Phi^{-1}(\bar y)$ and $y_k\neq\bar y$ for all $k\in\N$,
	satisfying the convergence properties 
	\begin{equation*}
		\begin{aligned}
			x_k&\to\bar x,&	\qquad	y_k&\to\bar y,& \qquad	\eta_k&\to 0,&
			\qquad
			\frac{x_k-\bar x}{\nnorm{x_k-\bar x}}&\to u,&
			\\
			v_k^\gamma&\to 0,&
			\qquad
			\nnorm{\lambda_k^\gamma} &\to \infty,&\qquad
			\nnorm{v_k^\gamma}\lambda_k^\gamma&\to y^*,&
			&&
		\end{aligned}
	\end{equation*}
	where we used
	\begin{equation}\label{eq:definition_of_surrogate_sequences}
		\forall k\in\N\colon\quad
		v^\gamma_k := \frac{y_k-\bar y}{\norm{x_k-\bar x}^\gamma},\qquad
		\lambda_k^\gamma:= k\nnorm{x_k-\bar x}^{\gamma-1}(y_k-\bar y),
	\end{equation}
	as well as
	\begin{equation}\label{eq:asymptotic_stationarity_regular_tools_gamma}
		\forall k\in\N\colon\quad
		\eta_k\in\partial\varphi(x_k)+D^*\Phi(x_k,y_k)\left(\frac{\lambda_k^\gamma}{\norm{x_k - \bar x}^{\gamma-1}}\right).
	\end{equation}
	\end{enumerate}
	Moreover, if $\Phi$ is metrically pseudo-subregular of order $\gamma$ at $(\bar x, \bar y)$
	in each direction $u\in\ker D\Phi(\bar x,\bar y)\cap \mathbb S_{\mathbb X}$, 
	$\bar x$ satisfies one of the alternatives~\ref{item:trivial_MSt},~\ref{item:pseudo_MSt_zero_mult}, or~\ref{item:pseudo_dir_MSt}.
\end{theorem}
\begin{proof}
	Let $\varepsilon>0$ be chosen such that $\varphi(x)\geq\varphi(\bar x)$ holds for all $x\in\mathcal F\cap\mathbb B_\varepsilon(\bar x)$
	and, for each $k\in\N$, consider the optimization problem
	\begin{equation}\label{eq:penalized_nonsmooth_problem}\tag{P$(k)$}
		\min\limits_{x,y}\{\varphi(x)+\tfrac k2\norm{y-\bar y}^2+\tfrac12\norm{x-\bar x}^2\,|\,(x,y)\in\gph\Phi,\,x\in\mathbb B_\varepsilon(\bar x)\}.
	\end{equation}
	For each $k\in\N$, the objective function of \eqref{eq:penalized_nonsmooth_problem} is bounded from below, continuous on the closed
	feasible set of this problem, and coercive in the variable $y$,
	so \eqref{eq:penalized_nonsmooth_problem} possesses a global minimizer $(x_k,y_k)\in\mathbb X\times\mathbb Y$. 
	By feasibility of $(\bar x,\bar y)$ for
	\eqref{eq:penalized_nonsmooth_problem}, we find
	\begin{equation}\label{eq:estimate_from_penalization}
		\forall k\in\N\colon\quad 
		\varphi(x_k)+\tfrac{k}{2}\norm{y_k-\bar y}^2+\tfrac12\norm{x_k-\bar x}^2\leq \varphi(\bar x).
	\end{equation}

	By boundedness of $\{x_k\}_{k\in\N}\subset\mathbb B_\varepsilon(\bar x)$,
	we may assume $x_k\to\tilde x$ for some $\tilde x\in\mathbb B_\varepsilon(\bar x)$.
	Observing that $\{\varphi(x_k)\}_{k\in\N}$ is bounded by continuity of $\varphi$, $y_k\to\bar y$ easily follows from
	\eqref{eq:estimate_from_penalization}. Furthermore, the closedness of $\gph\Phi$ guarantees $(\tilde x,\bar y)\in\gph\Phi$, i.e.,
	$\tilde x\in\mathcal F\cap\mathbb B_\varepsilon(\bar x)$ leading to $\varphi(\bar x)\leq\varphi(\tilde x)$.
	From \eqref{eq:estimate_from_penalization}, we find
	\begin{align*}
		\varphi(\bar x)
		\leq
		\varphi(\tilde x)
		\leq
		\varphi(\tilde x)+\tfrac12\norm{\tilde x-\bar x}^2
		=
		\lim\limits_{k\to\infty}\left(\varphi(x_k)+\tfrac12\norm{x_k-\bar x}^2\right)
		\leq
		\varphi(\bar x),
	\end{align*}
	and $\tilde x=\bar x$ follows. Thus, we have $x_k\to\bar x$.
	
	Let us assume that there is some $k_0\in\N$ such that $x_{k_0}$ is feasible to
	\eqref{eq:nonsmooth_problem}. 
	By \eqref{eq:estimate_from_penalization}, we find 
	\begin{align*}
		&\varphi(\bar x) +\tfrac{k_0}{2}\norm{y_{k_0}-\bar y}^2+\tfrac12\norm{x_{k_0}-\bar x}^2
		\\
		&\qquad
		\leq 
		\varphi(x_{k_0})+\tfrac{k_0}{2}\norm{y_{k_0}-\bar y}^2+\tfrac12\norm{x_{k_0}-\bar x}^2
		\leq
		\varphi(\bar x),
	\end{align*}
	i.e., $x_{k_0}=\bar x$ and $y_{k_0}=\bar y$. 
	Applying \cite[Theorem~6.1]{Mordukhovich2018}, 
	the subdifferential sum rule \cite[Theorem~2.19]{Mordukhovich2018}, and the
	definition of the limiting coderivative to find stationarity conditions
	of \eqref{eq:penalized_nonsmooth_problem} at $(\bar x,\bar y)$ yields
	$0\in\partial\varphi(\bar x)+D^*\Phi(\bar x,\bar y)(0)$,
	which is covered by~\ref{item:trivial_MSt}.
	
	Thus, we may assume that $x_k\notin\mathcal F=\Phi^{-1}(\bar y)$ holds for all $k\in\N$. 
	Particularly, $x_k\neq\bar x$ and $y_k\neq\bar y$ is valid for all $k\in\N$ in this situation.	
	Assume without loss of generality that $\{x_k\}_{k\in\N}$ belongs to the interior
	of $\mathbb B_\varepsilon(\bar x)$.

	We can apply Fermat's rule, see \cite[Proposition~1.30~(i)]{Mordukhovich2018}, 
	the semi-Lipschitzian sum rule for limiting subgradients from \cite[Corollary~2.20]{Mordukhovich2018}, 
	and the
	definition of the limiting coderivative in order to find
	\begin{equation}\label{eq:asymptotic_necessary_condition}
		\forall k\in\N\colon\quad
		\bar x-x_k\in\partial\varphi(x_k)+D^*\Phi(x_k,y_k)(k(y_k-\bar y)).
	\end{equation}
	Setting $\eta_k:=\bar x-x_k$ for each $k\in\N$, we find $\eta_k\to 0$.
	Since $\{(x_k - \bar x)/\nnorm{x_k - \bar x}\}_{k\in\N}\subset\mathbb S_{\mathbb X}$, 
	we may assume $(x_k - \bar x)/\nnorm{x_k - \bar x} \to u$ 
	for some $u\in\mathbb S_{\mathbb X}$.

	Next, we claim that 
	$\{y_k^*\}_{k\in\N}\subset\mathbb Y$, given by
	$y_k^* := k(y_k-\bar y)\norm{y_k-\bar y}/\norm{x_k - \bar x}$ for each $k\in\N$, is bounded.
	Rearranging \eqref{eq:estimate_from_penalization}, leaving a nonnegative term away,
	and division by $\nnorm{x_k-\bar x}$ give us
	\begin{equation}\label{eq:another_estimate_from_penalization}
		\forall k\in\N\colon
		\quad
		\frac{\varphi(x_k)-\varphi(\bar x)}{\nnorm{x_k-\bar x}}
		+
		\frac k2\frac{\nnorm{y_k-\bar y}^2}{\nnorm{x_k-\bar x}}
		\leq 
		0
		.
	\end{equation}
	Lipschitzianity of $\varphi$ yields boundedness of the first fraction, 
	so $\{k\nnorm{y_k-\bar y}^2/\nnorm{x_k-\bar x}\}_{k\in\N}$ must be bounded
	and, consequently, $\{y_k^*\}_{k\in\N}$ as well.
	Thus, we may assume $y_k^*\to y^*$ for some $y^*\in\mathbb Y$.
	
	Suppose that $\{(y_k-\bar y)/\nnorm{x_k-\bar x}\}_{k\in\N}$ does not converge to
	zero. This, along a subsequence (without relabeling), yields boundedness
	of the sequence $\{k(y_k-\bar y)\}_{k\in\N}$, 
	and taking the limit $k\to\infty$ in \eqref{eq:asymptotic_necessary_condition}
	along yet another subsequence while respecting robustness of the limiting
	subdifferential and the limiting coderivative yields~\ref{item:trivial_MSt}.
	
	Thus, we may assume $(y_k-\bar y)/\norm{x_k-\bar x}\to 0$ 
	for the reminder of the proof.
	Observe that we have
	\begin{equation*}
		(\bar x+ \norm{x_k - \bar x} (x_k - \bar x)/\norm{x_k - \bar x},\bar y+\norm{x_k - \bar x} (y_k - \bar y)/\norm{x_k - \bar x})
		\in
		\gph\Phi
	\end{equation*}
	for all $k\in\N$.
	Additionally, \eqref{eq:another_estimate_from_penalization} yields 
	\begin{equation*}
		\varphi(\bar x+ \norm{x_k - \bar x} (x_k - \bar x)/\norm{x_k - \bar x})-\varphi(\bar x) \leq 0,
	\end{equation*}
	so $u$ is a critical direction
	of order $(\gamma_0, 1)$ for each $\gamma_0 \geq 1$
	for \eqref{eq:nonsmooth_problem} at $\bar x$.
	
	In the remainder of the proof, we are going to exploit the sequences
	$\{v_k^\gamma\}_{k\in\N},\{\lambda_k^\gamma\}_{k\in\N}\subset\mathbb Y$ 
	given as in \eqref{eq:definition_of_surrogate_sequences}. 
	Observe that $y_k^* = \nnorm{v^\gamma_k} \lambda_k^\gamma $, 
	i.e., $\lambda_k^\gamma = y_k^*\norm{x_k - \bar x}^\gamma/\norm{y_k - \bar y}$ is valid for 
	each $k\in\N$.
	Note that the optimality condition \eqref{eq:asymptotic_necessary_condition} can be rewritten as
	\begin{equation}\label{eq:asymptotic_necessary_condition_of_order_gamma}
		\forall k\in\N\colon\quad
		\eta_k
		\in
		\partial\varphi(x_k)+D^*\Phi(x_k,\yb 
		+
		 \norm{x_k - \bar x}^\gamma v^\gamma_k)\left(\frac{\lambda_k^\gamma}{\norm{x_k - \bar x}^{\gamma-1}}\right).
	\end{equation}	
	Now, we need to distinguish three options.

	Let us assume that $\lambda_k^\gamma\to 0$.
	Using $t_k:=\norm{x_k-\bar x}$, 
	we can reformulate \eqref{eq:asymptotic_necessary_condition_of_order_gamma} as
	\[
		\forall k\in\N\colon\quad
		\eta_k
		\in
		\partial\varphi(x_k)
		+
		D^*\Phi\left(\bar x+t_k\frac{x_k-\bar x}{\norm{x_k-\bar x}},\bar y+t_k\frac{y_k-\bar y}{\norm{x_k-\bar x}}\right)
		\left(\frac{\lambda_k^\gamma}{t_k^{\gamma-1}}\right).
	\]
	Taking the limit $k\to\infty$ while respecting robustness of the directional limiting subdifferential
	as well as \cref{lem:pseudo_coderivative_via_limiting_normals} yields~\ref{item:pseudo_MSt_zero_mult}
	since $(y_k-\bar y)/\norm{x_k-\bar x}\to 0$ 
	and $u$ has already been shown to be critical for \eqref{eq:nonsmooth_problem} at $\bar x$.
	
	If $\{\lambda_k^\gamma\}_{k\in\N}$ remains bounded but, along a subsequence
	(without relabeling), stays away from zero, 
	we also get boundedness of $\{v_k^{\gamma}\}_{k\in\N}$
	from boundedness of $\{y_k^*\}_{k\in\N}$,
	and taking the limit along a convergent subsequence (without relabeling) in 
	\eqref{eq:asymptotic_necessary_condition_of_order_gamma} 
	while respecting robustness of the directional limiting subdifferential and \cref{lem:pseudo_coderivative_via_limiting_normals}
	yields precisely \eqref{eq:Pseudo-M-stat}, 
	where $\lambda, v \in\mathbb Y$ 
	with $\lambda\neq 0$ satisfy
	$\lambda_k^{\gamma}\to\lambda$ and $v_k^{\gamma}\to v$, respectively,
	and using $\alpha_k:=(k \norm{x_k - \bar x}^{2\gamma-1})^{-1}$ for all $k\in\N$
	as well as \eqref{eq:definition_of_surrogate_sequences},
	we find $v_k^{\gamma} = \alpha_k \lambda_k^{\gamma}$ for all $k\in\N$,
	$\alpha_k\to\norm{v}/\norm{\lambda}=:\alpha$, and $v=\alpha\lambda$.
	Criticality of $u$ for \eqref{eq:nonsmooth_problem} at $\bar x$ 
	has been shown above.
	Thus, situation~\ref{item:pseudo_dir_MSt} has been verified.

	If $\{\lambda_k^\gamma\}_{k\in\N}$ is not bounded, 
	we pass to a subsequence (without relabeling) such that $\nnorm{\lambda_k^{\gamma}} \to \infty$
	and so we also get $v_k^{\gamma} \to 0$ along this subsequence by boundedness of $\{y_k^*\}_{k\in\N}$.
	This means that $u$ is actually critical
	of order $(\gamma_0,\gamma)$ for \eqref{eq:nonsmooth_problem} at $\bar x$
	and so all conditions stated in~\ref{item:irregular_case_pseudo_MSt} have been verified
	since \eqref{eq:asymptotic_stationarity_regular_tools_gamma} follows from
	\eqref{eq:asymptotic_necessary_condition_of_order_gamma}.
	
	Finally, let us argue that option~\ref{item:irregular_case_pseudo_MSt} can be avoided, 
	i.e., that the sequence
	$\{\lambda_k^\gamma\}_{k\in\N}$ from above remains bounded
	if we assume that $\Phi$ is metrically pseudo-subregular
	of order $\gamma$ in direction $u$ at $(\xb,\yb)$.
	By boundedness of $\{y_k^*\}_{k\in\N}$, we immediately obtain the boundedness of $\{\lambda_k^\gamma\}_{k\in\N}$
	unless we have $v_k^\gamma\to 0$. Thus, let us assume the latter.
	By metric pseudo-subregularity of $\Phi$, 
	there is a constant $\kappa>0$ such that, for sufficiently large $k\in\N$, 
	we get the existence of $\tilde x_k \in \Phi^{-1}(\yb)$ with
	\begin{equation}\label{eq:estimate_pseudo_subregularity}
		\norm{x_k - \tilde x_k} 
		\leq
		\kappa \frac{\dist(\bar y,\Phi(x_k))}{\norm{ x_k-\bar x}^{\gamma - 1}} 
		\leq
		\kappa \frac{\norm{ y_k - \bar y}}{\norm{ x_k-\bar x}^{\gamma - 1}} 
		=
		\kappa \norm{x_k-\bar x} \nnorm{v_k^\gamma}.
	\end{equation}
	Particularly, we find $\nnorm{x_k-\tilde x_k}\to 0$ from $v_k^\gamma\to 0$,
	and $\tilde x_k\to\bar x$ follows.
	Since $(x_k,y_k)$ 
	is a global minimizer of \eqref{eq:penalized_nonsmooth_problem}, we get
	\[
		\frac{\varphi(x_k)-\varphi(\tilde x_k)}{\norm{x_k-\bar x}}
		+
		\frac k2\frac{\norm{y_k-\bar y}^2}{\norm{x_k-\bar x}}
		+
		\frac 12\frac{\norm{x_k-\bar x}^2 - \norm{\tilde x_k-\bar x}^2}{\norm{x_k-\bar x}}
		\leq 
		0
	\]
	for all sufficiently large $k\in\N$.
	Due to $v_k^\gamma \neq 0$ for all $k \in \N$,
	rearranging the above estimate and using
	\eqref{eq:definition_of_surrogate_sequences} as well as \eqref{eq:estimate_pseudo_subregularity} yield
	\[   
		\nnorm{\lambda_k^{\gamma}}
		=
		\frac{k \norm{y_k-\bar y}^2}{\norm{x_k-\bar x}\nnorm{v_k^{\gamma}}}
		\leq
		\frac{2\kappa \vert \varphi(x_k)-\varphi(\tilde x_k) \vert}{\norm{x_k - \tilde x_k}}
		+
		\frac{\kappa \big\vert \norm{x_k-\bar x}^2 - \norm{\tilde x_k-\bar x}^2\big\vert}{\norm{x_k - \tilde x_k}}.
	\]
	Boundedness of $\{\lambda_k^\gamma\}_{k\in\N}$
	thus follows from Lipschitzianity of $\varphi$
	and the estimate
	\[   
		\big\vert \norm{x_k-\bar x}^2 - \norm{\tilde x_k-\bar x}^2\big\vert
		\leq
		\norm{x_k - \tilde x_k}
		\big(\norm{\tilde x_k-\bar x} + \norm{x_k-\bar x}\big).
	\]
	This completes the proof.
\end{proof}

	Let us note that for the price of some more technicalities in the proof,
	involving the fuzzy sum rule for the regular subdifferential, see e.g.\ \cite[Exercise~2.26]{Mordukhovich2018},
	it is possible to formulate statement~\ref{item:irregular_case_pseudo_MSt} 	 
	in terms of the regular tools of variational analysis, 
	see \cite[Theorem~4.3]{BenkoMehlitz2022a} which is a preprint version of this paper.
	This more involved approach then allows for an easier comparison to available results
	in the literature which are partially stated in infinite dimensions, see e.g.\ 
	\cite{Gfrerer2014a}, where the limiting tools are of limited use and sequential 
	characterizations in terms of the regular tools are, thus, preferred.
	However, for our purposes, the way \cref{thm:higher_order_directional_asymptotic_stationarity} 
	has been formulated will be enough to proceed.

	In the rest of this subsection, we discuss some applications of
	\cref{thm:higher_order_directional_asymptotic_stationarity},
	which are then further developed in the rest of the paper.
	First, we focus on mixed-order stationarity conditions,
	involving first-order generalized derivatives of the objective function
	and pseudo-coderivatives of order $\gamma$,
	and enhance the result of \cref{cor: F-J} as follows.
	
\begin{corollary}\label{Cor:M_stat_order_gamma}
	Let $\bar x\in\mathcal F$ be a local minimizer of \eqref{eq:nonsmooth_problem} and consider $\gamma \geq 1$.
	Then the following assertions hold.
	\begin{enumerate}
	\item
	If $\Phi$ is metrically pseudo-subregular of order $\gamma$ at $(\bar x,\bar y)$ in each unit direction,
	then one of the alternatives~\ref{item:trivial_MSt},~\ref{item:pseudo_MSt_zero_mult}, or~\ref{item:pseudo_dir_MSt}
	of \cref{thm:higher_order_directional_asymptotic_stationarity} holds.
	\item\label{item:refined_conditions}
	If there are no critical directions of order $(1,\gamma)$ for \eqref{eq:nonsmooth_problem} at $\bar x$,
	then one of the alternatives~\ref{item:trivial_MSt},~\ref{item:pseudo_MSt_zero_mult}, or~\ref{item:pseudo_dir_MSt} 
    from \cref{thm:higher_order_directional_asymptotic_stationarity} is valid.
	If there exists a critical directions $u\in \mathbb S_{\mathbb X}$
	of order $(1,\gamma)$ for \eqref{eq:nonsmooth_problem} at $\bar x$
	satisfying \eqref{eq:FOSCMS_gamma}, then even
	\begin{equation}\label{eq:M_stat_order_gamma_dir_u}
		  0 \in \partial\varphi(\bar x;u) + D^\ast_{\gamma} \Phi((\xb,\yb);(u,0))(\lambda)
	\end{equation}
	holds for some $\lambda \in \mathbb Y$.
	\end{enumerate}
\end{corollary}
\begin{proof}
	The first assertion follows directly from
	\cref{thm:higher_order_directional_asymptotic_stationarity}.
	Let us now prove the second assertion.
	\cref{thm:higher_order_directional_asymptotic_stationarity}
	says that either one of the alternatives~\ref{item:trivial_MSt},~\ref{item:pseudo_MSt_zero_mult}, or~\ref{item:pseudo_dir_MSt} 
	holds, or there exists a critical direction $u\in \mathbb S_{\mathbb X}$
	of order $(1,\gamma)$ for \eqref{eq:nonsmooth_problem} at $\bar x$ (with certain properties).
	If among these critical directions, there is one that satisfies
	\eqref{eq:FOSCMS_gamma}, \cref{cor: F-J} yields
	\eqref{eq:M_stat_order_gamma_dir_u}.
\end{proof}

We conjecture that the sufficient condition \eqref{eq:FOSCMS_gamma}
can be weakened to just pseudo-subregularity of $\Phi$
in \cref{Cor:M_stat_order_gamma}\,\ref{item:refined_conditions}.
However, it would require a different proof to show this, so we will not explore this option.
For $\gamma:=1$, such a result is known to hold, see \cite[Theorem~7]{Gfrerer2013}.

	Note that \eqref{eq:M_stat_order_gamma_dir_u} is covered
	by the alternative~\ref{item:pseudo_MSt_zero_mult}
	(if $\lambda = 0$, see \eqref{eq:trivial_upper_estimate_pseudo_coderivative}) 
	or~\ref{item:pseudo_dir_MSt}
	(if $\lambda \neq 0$)
	of \cref{thm:higher_order_directional_asymptotic_stationarity}.
	Hence, the optimality conditions from \cref{Cor:M_stat_order_gamma} give 
	either M-stationarity of the underlying local minimizer or
	validity of alternative~\ref{item:pseudo_MSt_zero_mult} or~\ref{item:pseudo_dir_MSt}
	of \cref{thm:higher_order_directional_asymptotic_stationarity}
	for some critical direction (of order $(1,1)$ or $(1,\gamma)$).

\begin{remark}\label{Rem: 2_paths_CQ_vs_Crit_dir}
\cref{Cor:M_stat_order_gamma} offers two distinct paths to optimality condition 
of type ``M-stationarity or \eqref{eq:M_stat_order_gamma_dir_u}'',
both with some advantages and disadvantages.
\begin{enumerate}
\item\label{item:all_directions_general}
Assuming pseudo-subregularity in \emph{each} unit direction
yields this type of condition by
ruling out the alternative~\ref{item:irregular_case_pseudo_MSt} of \cref{thm:higher_order_directional_asymptotic_stationarity}.
However, this can sometimes be an undesirable type of assumption as pointed out in \cref{rem:D=0}.
\item\label{item:crit_directions_general}
The refined assumptions in \cref{Cor:M_stat_order_gamma}\,\ref{item:refined_conditions}
are clearly milder, but they depend on a critical direction (of order $(1,\gamma)$), 
which in turn depends also on the objective function, not just on $\Phi$.
These assumptions do not rule out the alternative~\ref{item:irregular_case_pseudo_MSt}.
Instead, they just secure that~\ref{item:trivial_MSt},~\ref{item:pseudo_MSt_zero_mult}, or~\ref{item:pseudo_dir_MSt} 
from \cref{thm:higher_order_directional_asymptotic_stationarity} holds.
\end{enumerate}
These two types of assumptions will be prevalent throughout this section.
\end{remark}

	Recall that all the assumptions in \cref{Cor:M_stat_order_gamma} become less restrictive as
	$\gamma$ increases, see \cref{sec:sufficient_conditions_pseudo_regularity} as well.
	On the contrary, with increasing $\gamma$, the involved pseudo-coderivatives become more difficult 
	to handle which, exemplary, can be seen for constraint mappings when comparing the cases 
	$\gamma:=1$ and $\gamma:=2$ from \cref{sec:variational_analysis_constraint_mapping}.
	In this regard, in \cref{Cor:M_stat_order_gamma}, $\gamma$ should be chosen as small as possible such that
	the exploited qualification condition is valid.

In \cref{sec:constraint_mappings}, we work out the conditions from \cref{Cor:M_stat_order_gamma} for $\gamma:=2$ 
in the setting where $\Phi$ is a constraint mapping as the appearing pseudo-coderivatives actually can be computed, 
see \cref{sec:variational_analysis_constraint_mapping}, and, hence, we obtain conditions in terms of initial problem data.
In \cref{sec:applications}, 
we further apply these results to two specific problem classes
and compare them with similar results based on 2-regularity.
	
\cref{thm:higher_order_directional_asymptotic_stationarity} also opens a way to the
identification of new conditions which guarantee that local minimizers of
\eqref{eq:nonsmooth_problem} are M-stationary.
One of the most prominent conditions that implies this is metric subregularity,
and the corresponding result, which we state next, can be obtained simply by setting $\gamma := 1$
in \cref{Cor:M_stat_order_gamma}, taking also into account \cite[Theorem~7]{Gfrerer2013}.
For us, this result serves as a basis for comparison.
Later on, we will derive new conditions, which are independent of (directional) metric subregularity,
but which are milder than various known sufficient conditions for metric subregularity.

\begin{corollary}\label{lem:directional_M_stationarity_via_metric_subregularity}
	A local minimizer $\bar x\in\mathcal F$ of \eqref{eq:nonsmooth_problem} is M-stationary
	if one of the following conditions holds.
\begin{enumerate}
\item\label{item:all_directions_subreg} The mapping $\Phi$ is metrically subregular at $(\bar x,\bar y)$ in each unit direction.
\item\label{item:crit_directions_subreg} There are no critical directions for \eqref{eq:nonsmooth_problem} at $\bar x$,
or there is a critical direction $u\in\mathbb S_{\mathbb X}$ for \eqref{eq:nonsmooth_problem} at $\bar x$
and $\Phi$ is metrically subregular at $(\bar x,\bar y)$ in direction $u$,
in which case there is a multiplier $\lambda\in\mathbb Y$ such that
	\[
		0\in\partial\varphi(\bar x;u)+D^*\Phi((\bar x,\bar y);(u,0))(\lambda).
	\]
\end{enumerate}
\end{corollary}

Let us now discuss two novel approaches to M-stationarity.
The first approach corresponds to using \cref{Cor:M_stat_order_gamma} with $\gamma > 1$ and then making sure that the
derived mixed-order conditions in terms of pseudo-coderivatives actually yield M-stationarity.
To formalize the idea, we introduce the following assumption.
\begin{assumption}\label{ass:A_gamma_u}
	Given $u \in \mathbb S_{\mathbb X}$ and $\gamma \geq 1$, we say that $A^\gamma(u)$ holds if \eqref{eq:FOSCMS_gamma} is satisfied and
	\begin{equation}\label{eq:sufficient_condition_asym_reg_pseudo}
		\widetilde{D}^*_\gamma\Phi((\bar x,\bar y);(u,0))(0)
		\cup
		\bigcup_{w\in\mathbb S_{\mathbb Y}}
		D^*_\gamma\Phi((\bar x,\bar y);(u,\alpha w))(\beta w)
		\subset
		\Im D^*\Phi(\bar x,\bar y)
	\end{equation}
	is valid for all $\alpha,\beta \geq 0$.
\end{assumption}
Let us mention that
\begin{equation}\label{eq:sufficient_condition_asym_reg_pseudo_rough}
		\Im\widetilde D^*_\gamma\Phi((\bar x,\bar y);(u,0))
		\subset
		\Im D^*\Phi(\bar x,\bar y)
\end{equation}
is a sufficient condition for \eqref{eq:sufficient_condition_asym_reg_pseudo} due to 
\eqref{eq:trivial_upper_estimate_pseudo_coderivative}.
Assumption $A^\gamma(u)$ leads to the problem of how to compute or estimate the appearing pseudo-coderivatives.
As mentioned above, for $\gamma:=2$ and in the setting where $\Phi$ is a constraint mapping, these objects can be
computed and assumption $A^\gamma(u)$ can be made explicit.
We discuss this case in detail in \cref{Sec:5_constr_mappings}, where we show
that $A^2(u)$ is (strictly) weaker than FOSCMS$(u)$ as well as its refinement SOSCMS$(u)$ in the polyhedral case.
Here, we just explain how $A^\gamma(u)$ can be used to secure M-stationarity and how to compare it
with sufficient conditions for metric subregularity.

To proceed, let $\bar x\in\mathcal F$ be a local minimizer of \eqref{eq:nonsmooth_problem} and consider $\gamma\geq 1$.
Assuming that $A^\gamma(u)$ holds in every unit direction $u$ implies
that options~\ref{item:pseudo_MSt_zero_mult} or~\ref{item:pseudo_dir_MSt} 
from \cref{thm:higher_order_directional_asymptotic_stationarity} yield M-stationarity of $\bar x$,
and that option~\ref{item:irregular_case_pseudo_MSt} cannot occur.
Thus, we end up with $\bar x$ being M-stationary.
Similarly, if $A^\gamma(u)$ holds in a critical direction $u\in\mathbb S_{\mathbb X}$ 
of order $(1,\gamma)$ for \eqref{eq:nonsmooth_problem} at $\bar x$,
\eqref{eq:M_stat_order_gamma_dir_u} is satisfied.
Due to \eqref{eq:sufficient_condition_asym_reg_pseudo}, this also shows M-stationarity of $\bar x$.
Thus, we obtain the following from \cref{Cor:M_stat_order_gamma}.
\begin{corollary}\label{lem:M_stationarity_via_A_gamma}
	Let $\bar x\in\mathcal F$ be a local minimizer of \eqref{eq:nonsmooth_problem} and consider $\gamma\geq 1$.
	Then each of the following conditions implies that $\bar x$ is M-stationary.
	\begin{enumerate}
		\item\label{item:all_directions_gamma} 
			Condition $A^\gamma(u)$ holds in each unit direction $u$.	
		\item\label{item:crit_directions_gamma} 
			There are no critical directions for \eqref{eq:nonsmooth_problem} at $\bar x$,
			or there is a critical direction $u\in\mathbb S_{\mathbb X}$ of order $(1,\gamma)$ for \eqref{eq:nonsmooth_problem} at $\bar x$
			such that $A^\gamma(u)$ holds.	
	\end{enumerate}
\end{corollary}

In the following remark, we compare our approach from
\cref{lem:M_stationarity_via_A_gamma} with the results from
\cref{lem:directional_M_stationarity_via_metric_subregularity} in the presence
of any sufficient condition for directional metric subregularity.
\begin{remark}
	Due to \cref{lem:directional_M_stationarity_via_metric_subregularity},
	directional metric subregularity serves as a constraint qualification guaranteeing M-stationarity
	of local minimizers.
	However, given $u\in\mathbb S_{\mathbb X}$, metric subregularity in direction $u$ is difficult to verify,
	so it is often replaced by some stronger condition which is easier to check - exemplary, FOSCMS$(u)$.
	Let us label such a sufficient condition as SCMS$(u)$.
	Clearly, \cref{lem:directional_M_stationarity_via_metric_subregularity}
	can be restated in terms of SCMS$(u)$.
	Suppose that we can show that $A^\gamma(u)$ is milder than SCMS$(u)$
	for every $u\in\mathbb S_{\mathbb X}$ (even strictly milder for some $u$).
	Naturally, option~\ref{item:all_directions_gamma} from \cref{lem:M_stationarity_via_A_gamma}
	then provides a (strictly) milder assumption than requiring
	SCMS$(u)$ to hold for all unit directions.
	However, does an analogous relationship hold for the more complicated 
	option~\ref{item:crit_directions_gamma} from \cref{lem:M_stationarity_via_A_gamma}?
	Both approaches yield M-stationarity if there are no critical directions.
	If there is a critical direction $u\in\mathbb S_{\mathbb X}$ such that SCMS$(u)$
	holds, then \cref{lem: crit_dir+pseudo_subreg}
	yields that $u$ is actually critical of order $(1,\gamma)$ and, thus, the
	milder assumption $A^\gamma(u)$ from the case \ref{item:crit_directions_gamma} 
	of \cref{lem:M_stationarity_via_A_gamma} can be applied.
	This means that our approach via \cref{lem:M_stationarity_via_A_gamma}
	is indeed better than an approach via any
	sufficient condition for metric subregularity in direction $u$
	which is stronger that $A^\gamma(u)$.
\end{remark}

	The second approach to M-stationarity can be called ``asymptotic'' 
	and is based on the following result,
	a generalization of \cite[Theorem~3.9]{Mehlitz2020a}, which
	reinspects \cref{thm:higher_order_directional_asymptotic_stationarity}
	in the situation $\gamma:=1$.
	Particularly, we exploit that, in this case,
	both notions of a pseudo-coderivative from \cref{def:coderivatives}
	coincide with the directional limiting coderivative.

\begin{corollary}\label{thm:directional_asymptotic_stationarity}
	Let $\bar x\in\mathcal F$ be a local minimizer of \eqref{eq:nonsmooth_problem}.
	Then $\bar x$ is M-stationary or there exist a critical direction 
	$u\in\mathbb S_{\mathbb X}$ for \eqref{eq:nonsmooth_problem} at $\bar x$,
	some $y^* \in \mathbb Y$, and sequences
	$\{x_k\}_{k\in\N},\{\eta_k\}_{k\in\N}\subset\mathbb X$ as well as
	$\{y_k\}_{k\in\N}\subset\mathbb Y$ such that
	$x_k\notin\Phi^{-1}(\bar y)$ and $y_k\neq\bar y$ for all $k\in\N$,
	satisfying the convergence properties
	\begin{subequations}\label{eq:convergences_gamma=1}
		\begin{align}
			\label{eq:convergences_gamma=1_basic}
				x_k&\to\bar x,&		\qquad	y_k&\to\bar y,&	\qquad	\eta_k&\to 0,&
				\\
			\label{eq:convergences_gamma=1_directional}
				\frac{x_k-\bar x}{\nnorm{x_k-\bar x}}&\to u,&	\qquad
				\frac{y_k-\bar y}{\nnorm{x_k-\bar x}}&\to 0,&	\qquad
				&&
				\\
			\label{eq:convergences_gamma=1_multiplier}
			k\frac{\nnorm{y_k-\bar y}}{\nnorm{x_k-\bar x}}(y_k-\bar y)&\to y^*,&\qquad
			k\nnorm{y_k-\bar y}&\to \infty,&\qquad
			&&
		\end{align}
	\end{subequations}	
	and
	\begin{equation}\label{eq:asymptotic_stationarity_gamma=1}
		\forall k\in\N\colon\quad
		\eta_k\in\partial\varphi(x_k)+D^*\Phi(x_k,y_k)\left(k(y_k-\bar y)\right).
	\end{equation}
\end{corollary}

The above result shows that each local minimizer of \eqref{eq:nonsmooth_problem}
either is M-stationary or satisfies asymptotic
stationarity conditions w.r.t.\ a certain critical direction and an unbounded
sequence of multiplier estimates $\{\lambda_k\}_{k\in\N}$ given by
\begin{equation}\label{eq:multiplier_estimate}
	\forall k\in\N\colon\quad
	\lambda_k:=k(y_k-\bar y).
\end{equation}
Note that in the case where $\{\lambda_k\}_{k\in\N}$ would be bounded, one could simply
take the limit in \eqref{eq:asymptotic_stationarity_gamma=1} along a suitable
subsequence and, respecting the
convergences from \eqref{eq:convergences_gamma=1_basic}, would end up with M-stationarity
again taking into account robustness of the limiting subdifferential and coderivative. 
Thus, divergence of the multiplier estimates is natural since not all local minimizers of
\eqref{eq:nonsmooth_problem} are M-stationary in general, see \cite[Lemma~3.4]{Mehlitz2020a} as well.
Related results in nondirectional form can be found
in \cite{KrugerMehlitz2021,Mehlitz2020a}. The story of asymptotic stationarity
conditions in variational analysis, however, can be traced back to
\cite{Kruger1985,KrugerMordukhovich1980}. This concept has been rediscovered as a valuable
tool for the analysis of convergence properties for solution algorithms associated
with standard nonlinear optimization problems in 
\cite{AndreaniMartinezSvaiter2010,AndreaniHaeserMartinez2011}, and extensions were
made to disjunctive, conic, and even infinite-dimensional optimization,
see e.g.\ \cite{AndreaniHaeserSecchinSilva2019,AndreaniGomezHaeserMitoRamos2021,BoergensKanzowMehlitzWachsmuth2019,Ramos2019}
and the references therein.

The sequential information from \eqref{eq:convergences_gamma=1} describes in great detail 
what must ``go wrong'' if M-stationarity fails.
We will refer to \eqref{eq:convergences_gamma=1_basic}, \eqref{eq:convergences_gamma=1_directional}, 
and \eqref{eq:convergences_gamma=1_multiplier} as
basic, directional, and multiplier (sequential) information, respectively.
Clearly, one can secure M-stationarity of a local minimizer 
by ruling out the second alternative in \cref{thm:directional_asymptotic_stationarity} and,
as we will show, various known constraint qualifications for M-stationarity indeed do precisely that.
Let us mention here two such conditions.
Rescaling \eqref{eq:asymptotic_stationarity_gamma=1} by
$\norm{\lambda_k}$, for $\{\lambda_k\}_{k\in\N}$ as given in \eqref{eq:multiplier_estimate}, 
and taking the limit $k\to\infty$ leads to a contradiction 
with the Mordukhovich criterion \eqref{eq:Mordukhovich_criterion}, i.e.,  
metric regularity of $\Phi$ at $(\bar x,\bar y)$.
Respecting also the directional information \eqref{eq:convergences_gamma=1_directional} 
yields a contradiction with FOSCMS$(u)$ at $(\bar x,\bar y)$.

In both cases, we have essentially discarded the multiplier information \eqref{eq:convergences_gamma=1_multiplier}
which deserves some remarks.
We have used $\norm{\lambda_k} \to \infty$, but this information is not really very important 
since, as we already explained,
if the multipliers remain bounded, we end up with M-stationarity anyway.
The fact that $\{\lambda_k \norm{y_k-\bar y}/\norm{x_k-\bar x}\}_{k\in\N}$ converges tells us how fast the multipliers 
$\{\lambda_k\}_{k\in\N}$ blow up.
We note that the concept of super-coderivatives from \cref{def:super_coderivative} collects this information,
and we will come back to it in \cref{sec:asymptotic_regularity_via_super_coderivative}, 
where it is used to design constraint qualifications for M-stationarity.
	As we will show in \cref{sec:asymptotic_regularity_via_super_coderivative}, 
	this approach is closely related to the hypothesis $A^\gamma(u)$ which we formulated in \cref{ass:A_gamma_u},
	and its role as a constraint qualification has already been illustrated in \cref{lem:M_stationarity_via_A_gamma}.

Finally, note that $(y_k-\bar y)/\nnorm{y_k-\bar y} = \lambda_k/\nnorm{\lambda_k}$
means that the multipliers
precisely capture the direction from which $\{y_k\}_{k\in\N}$ converges to $\bar y$.
Particularly, we find $\dual{\lambda_k/\nnorm{\lambda_k}}{(y_k-\bar y)/\nnorm{y_k-\bar y}} = 1$,
which is clearly more restrictive than the condition $\dual{\lambda_k/\nnorm{\lambda_k}}{(y_k-\bar y)/\nnorm{y_k-\bar y}} \to 1$.
The latter convergence,
which is used in the sufficient condition for metric subregularity in \cite[Corollary~1]{Gfrerer2014a},
can be recast as $(y_k-\bar y)/\nnorm{y_k-\bar y} - \lambda_k/\nnorm{\lambda_k}\to 0$.
This information is respected by the new constraint qualifications which we are going to suggest
in \cref{sec:directional_asymptotic_regularity}.

\subsection{Mixed-order necessary optimality conditions for optimization problems with geometric constraints in the case $\gamma:=2$}
\label{sec:constraint_mappings}

In this part, we apply \cref{Cor:M_stat_order_gamma} with $\gamma := 2$ to the case where 
$\Phi\colon\mathbb X\tto\mathbb Y$ is given in
the form of a constraint mapping, i.e., $\Phi(x) := g(x) - D$, $x\in\mathbb X$, holds where $g\colon\mathbb X\to\mathbb Y$ is twice continuously differentiable 
and $D \subset \mathbb Y$ is a closed set.
Since, in \cref{sec:variational_analysis_constraint_mapping},
we computed the pseudo-coderivative and the graphical pseudo-derivative of order $2$ of $\Phi$,
we are able to derive explicit conditions in terms of initial problem data.
For that purpose, we assume $\bar y:=0$ in \eqref{eq:nonsmooth_problem} 
throughout the section which can be done
without loss of generality.

We start with a description of critical directions of order $(1,2)$ and $(2,2)$.
\begin{lemma}\label{Lem:Crit_dir_constraints}
	Fix $\bar x\in\mathcal F$ and let $u \in \mathbb S_{\mathbb X}$
	be a critical direction of order $(1,2)$ of \eqref{eq:nonsmooth_problem} at $\bar x$.
	Suppose that $\mathbb Y:=\R^m$ and $D$ is locally polyhedral around $g(\xb)$.
	Then
	\begin{equation*}
		u
		\in
		\mathcal C^{1,2}(\bar x)
		:=
		\{
			u \in \mathbb X
			\,|\,
			\dr \varphi(\bar x)(u) \leq 0,\,
			\exists s \in \mathbb X\colon\, w_s(u,0) \in \mathbf T(u) 
		\},
	\end{equation*}
	where $w_s(u,0)$ and $\mathbf T(u)$ are defined in \eqref{eq:Tu_and_ws}.
	If $\varphi$ is continuously differentiable, $\mathcal C^{1,2}(\bar x)$
	corresponds precisely to the set of critical directions of order $(1,2)$
	of \eqref{eq:nonsmooth_problem} at $\bar x$.
	Moreover, if $\varphi$ is even twice continuously differentiable at $\bar x$,
	the set of all critical directions of order $(2,2)$ 
	of \eqref{eq:nonsmooth_problem} at $\bar x$
	equals
	\begin{equation}\label{eq:Crit_dir_2_2_constraint}
		\{
			u \in \mathbb X
			\,|\,
			\exists s \in \mathbb X\colon\,
			\nabla \varphi(\bar x) s + 1/2 \nabla^2 \varphi(\bar x)[u,u] \in \mathcal T_{\R_-}(\nabla \varphi(\bar x) u),
			w_s(u,0) \in \mathbf T(u) 
		\}.
	\end{equation}
\end{lemma}
\begin{proof}
	A critical direction $u$ of order $(1,2)$ of \eqref{eq:nonsmooth_problem} at $\bar x$ satisfies
	$\dr \varphi(\bar x)(u) \leq 0$ and $0 \in D_2 \Phi(\bar x,0)(u)$,
	with equivalence being valid if $\varphi$ is continuously differentiable at $\bar x$.
	Hence, the first statement follows from \cref{The : NCgen_2}.
	\\
	By \cref{Pro : critical_direction_interpretation},
	a direction $u$ is critical of order $(2,2)$ of \eqref{eq:nonsmooth_problem} at $\bar x$
	if and only if
	$u \in \ker D_2 M(\bar x,(\varphi(\bar x),0))$
	for $M\colon\mathbb X\tto\R\times\R^m$ given by
	$M(x):=(\varphi(x),g(x)) - (\R_- \times D)$, $x\in\mathbb X$.
	Hence, \cref{The : NCgen_2} can be applied again,
	yielding the second statement.
\end{proof}

\begin{remark}
 Note that \cref{Lem:Crit_dir_constraints}
 shows that the set of directions $C_2(\bar x)$ from
 \cite[Theorem~3]{AvakovArutunovIzmailov2007}
 and its extension labeled \emph{second-order tightened critical cone} in
 \cite[Theorem~3]{ArutyunovAvakovIzmailov2008}
 actually correspond to $\mathcal C^{1,2}(\bar x)$,
 while the set of directions used in
 \cite[Theorem~3(2.)]{Gfrerer2014a}
 corresponds to the one in \eqref{eq:Crit_dir_2_2_constraint}.
 We believe that interpreting these directions
 as critical (of some order) is very natural.
 Moreover, our approach justifies the name.
 Indeed, as already mentioned, our definition of criticality is an extension of 
 the one stated in \cite[Definition~5]{Gfrerer2013}.
 More importantly, we have shown in \cref{Cor:M_stat_order_gamma}
 that in the absence of nonzero critical directions 
 (of order $(1,\gamma)$ for some $\gamma\geq 1$),
 the corresponding mixed-order optimality conditions 
 (involving a pseudo-coderivative of order $\gamma$) are satisfied
 without any additional assumptions.
 \end{remark}

Based on \cref{The : NCgen,The : NCgen_2} as well as \cref{Cor:M_stat_order_gamma}, we obtain the following result.

\begin{proposition}\label{Pro:M-stat_via_second_order}
 Let $\bar x\in\mathcal F$ be a local minimizer of \eqref{eq:nonsmooth_problem}.
 	\begin{enumerate}
  	\item\label{item:pseudo_stationarity_constraint_maps}
		If \eqref{eq:CQ_pseudo_subregularity_II},
  		as well as \eqref{eq:CQ_pseudo_subregularity_Ia} or, in the case $\nabla g(\bar x)u\neq 0$,
  		\eqref{eq:CQ_pseudo_subregularity_Ib} hold for every unit direction, then $\xb$ is M-stationary or there exist
  		a critical direction $u \in \mathbb S_{\mathbb X}$
  		and
  		\begin{equation}\label{eq:Props_of_y*_z*}
  			y^* \in \mathcal N_{D}(g(\xb);\nabla g(\xb)u)  \cap \ker\nabla g(\xb)^*,
  			\qquad 	
  			z^* \in D\mathcal N_{D}(g(\xb),y^*)(\nabla g(\xb)u)
  		\end{equation}
  		such that
  		\begin{equation}\label{eq:Sec_order_optim_cond}
   				0 
   				\in 
   				\partial\varphi(\bar x;u) 
   				+ 
   				\nabla^2\innerprod{y^*}{g}(\xb)(u) + \nabla g(\xb)^* z^*.
  		\end{equation}
		If there exists a critical direction $u \in \mathbb S_{\mathbb X}$
		of order $(1,2)$ of \eqref{eq:nonsmooth_problem} at $\bar x$
		satisfying \eqref{eq:CQ_pseudo_subregularity_II},
  		as well as \eqref{eq:CQ_pseudo_subregularity_Ia} or, in the case $\nabla g(\bar x)u\neq 0$,
  		\eqref{eq:CQ_pseudo_subregularity_Ib},
   		then there exist $y^*,z^*\in\mathbb Y$ satisfying
   		\eqref{eq:Props_of_y*_z*} and \eqref{eq:Sec_order_optim_cond} for this $u$.
   		\item\label{item:pseudo_stationarity_constraint_maps_polyhedral_refined} 
  		Let $\mathbb Y:=\R^m$ and $D$ be locally polyhedral around $g(\xb)$.
  		If either $\mathcal C^{1,2}(\bar x) = \{0\}$
  		or if \eqref{eq:CQ_pseudo_subregularity_polyhedral_II}
  		holds for every unit direction,
  		then $\xb$ is M-stationary or there exist
  		a critical direction $u \in \mathbb S_{\mathbb X}$,
  		$s\in\mathbb X$, $y^*, z_i^* \in \R^m$ for $i=1,2$, and $\alpha \geq 0$,
  		satisfying $\nabla g(\xb)^* y^* =0$,
  		\[
  			y^*, z_1^* \in \mathcal N_{\mathbf T(u)}(w_s(u,v)),
  	 		\qquad
  	 		z_2^* \in \mathcal T_{\mathcal N_{\mathbf T(u)}(w_s(u,v))}(y^*),
  		\]
  		and \eqref{eq:Sec_order_optim_cond} (with $z^* = z_i^*, i=1,2$),
  		where $v:= \alpha y^*$, and $w_s(u,v)$ and $\mathbf T(u)$ have been defined in \eqref{eq:Tu_and_ws}.
  		If there exists 
		$u \in \mathcal C^{1,2}(\bar x) \cap \mathbb S_{\mathbb X}$ satisfying
  		\eqref{eq:CQ_pseudo_subregularity_polyhedral_II},
  		then there exist $s\in\mathbb X$ and
  		\[
  			y^*, z_1^* \in \mathcal N_{\mathbf T(u)}(w_s(u,0)),
  	 		\qquad
  	 		z_2^* \in \mathcal T_{\mathcal N_{\mathbf T(u)}(w_s(u,0))}(y^*)
  		\]
  		satisfying \eqref{eq:Sec_order_optim_cond} (with $z^* = z_i^*, i=1,2$)
  		as well as $\nabla g(\xb)^* y^* =0$.
 	\end{enumerate}
\end{proposition}
\begin{proof}
		For the proof of~\ref{item:pseudo_stationarity_constraint_maps},
		in the first alternative,
		we apply \cref{cor:sufficient_condition_pseudo_subregularity}
		in order to verify that \eqref{eq:FOSCMS_gamma} holds for every unit direction.
		\cref{Cor:M_stat_order_gamma} in turn yields
		that $\bar x$ is M-stationary or
		one of the cases~\ref{item:pseudo_MSt_zero_mult} and~\ref{item:pseudo_dir_MSt}
		from \cref{thm:higher_order_directional_asymptotic_stationarity} holds.
		In the case of 
		\cref{thm:higher_order_directional_asymptotic_stationarity}\,\ref{item:pseudo_MSt_zero_mult}, 
		however, from \cref{The : NCgen}\,\ref{item:general_estimate_+CQ} we get
		$0 \in \partial \varphi(\xb;u) + \nabla g(\xb)^* z^*$ 
		with
		\[      
		 z^*
		 \in
		 D\mathcal N_{D}(g(\xb),0)(\nabla g(\xb)u)
		 \subset
		 \mathcal N_D(g(\bar x);\nabla g(\bar x) u)
		 \subset
		 \mathcal N_D(g(\bar x)),
		\]
		see \cref{lem:technical_property_D_normal_cone_map},
		and M-stationarity of $\bar x$ follows.
		In the case of \cref{thm:higher_order_directional_asymptotic_stationarity}\,\ref{item:pseudo_dir_MSt}, 
		from
		\eqref{eq:trivial_upper_estimate_pseudo_coderivative}
		and
		\cref{The : NCgen}\,\ref{item:general_estimate_+CQ}
		we precisely obtain $y^*$ and $z^*$ as stated.
		Similarly, the second alternative follows from
		successively applying
		\cref{cor:sufficient_condition_pseudo_subregularity,Cor:M_stat_order_gamma},
		\eqref{eq:trivial_upper_estimate_pseudo_coderivative}, and
		\cref{The : NCgen}\,\ref{item:general_estimate_+CQ}.

		For the proof of~\ref{item:pseudo_stationarity_constraint_maps_polyhedral_refined}, we first would like to hint to \cref{Lem:Crit_dir_constraints}.
		In the first alternative, taking into account
		\cref{cor:sufficient_condition_pseudo_subregularity},
		\cref{Cor:M_stat_order_gamma} yields
		that $\bar x$ is M-stationary or
		one of the cases~\ref{item:pseudo_MSt_zero_mult} and~\ref{item:pseudo_dir_MSt}
		from \cref{thm:higher_order_directional_asymptotic_stationarity} holds.
		As before, in the case of 
		\cref{thm:higher_order_directional_asymptotic_stationarity}\,\ref{item:pseudo_MSt_zero_mult}, 
		from \cref{The : NCgen}\,\ref{item:polyhedral_estimate} we get
		$0 \in \partial \varphi(\xb;u) + \nabla g(\xb)^* z^*$ 
		with $z^* \in \mathcal N_{\mathcal T_D(g(\bar x))}(\nabla g(\bar x) u) \subset \mathcal N_D(g(\bar x))$,
		see \cref{lem:some_properties_of_polyhedral_sets}, 
		and M-stationarity of $\bar x$ follows.
		In the case of \cref{thm:higher_order_directional_asymptotic_stationarity}\,\ref{item:pseudo_dir_MSt}, 
		from \cref{The : NCgen_2}
		we precisely obtain $y^*, z_1^*$, and $z_2^*$ as stated.
		The second alternative follows from
		\cref{cor:sufficient_condition_pseudo_subregularity,Cor:M_stat_order_gamma}
		as well as \cref{The : NCgen_2}.
\end{proof}

Similar optimality conditions involving a mixture of
first- and second-order derivatives were proposed e.g.\ in
\cite{ArutyunovAvakovIzmailov2008,Avakov1985,Avakov1989,AvakovArutunovIzmailov2007,Gfrerer2014a}.
Let us now explain that in the convex polyhedral case,
where $\mathbb Y:=\R^m$ holds while $D$ is convex and polyhedral,
all these optimality conditions are the same and can be stated simply as follows:
If there exists
$u \in \mathcal C^{1,2}(\bar x) \cap \mathbb S_{\mathbb X}$
satisfying \eqref{eq:2regularity_dual},
then there are $y^*, z^* \in \R^m$ satisfying 
\begin{equation}\label{eq:Sec_order_optim_cond_smooth}
	\nabla \varphi(\bar x) + \nabla^2\langle y^*,g\rangle(\bar x)(u)+\nabla g(\bar x)^*z^*=0,\,
	\nabla g(\bar x)^* y^*=0,\,
	y^*, z^* \in \mathcal N_D(g(\bar x))
\end{equation}
(for a fair comparison, we assume that $\varphi$ is continuously differentiable).

In \cref{ex:some_trivial_example,ex:2regularity_stronger},
we have shown that the 2-regularity assumption \eqref{eq:2regularity_dual}
used in \cite{ArutyunovAvakovIzmailov2008} is, in general,
strictly stronger than our condition
\eqref{eq:some_sufficient_condition_for_pseudo_regularity},
which is, in turn, strictly stronger than
the mutually equivalent conditions
\eqref{eq:Gfrerers_sufficient_condition_implicit}
from \cite{Gfrerer2014a} and
\eqref{eq:CQ_pseudo_subregularity_polyhedral_II}
from \cref{cor:sufficient_condition_pseudo_subregularity}.
However, as shown in \cref{cor:4_conditions_convex_polyhedral_case},
all these assumptions are equivalent if applied
to a critical direction $u$ of order $(1,2)$, i.e, 
$u \in \mathcal C^{1,2}(\bar x)$, as this yields the existence of $s \in \mathbb X$
with $w_s(u,0) \in \mathbf T(u)$.

Clearly, although the aforementioned qualification conditions are equivalent, 
the optimality conditions may differ due to the additional information regarding
the multipliers.
However, this is also not the case, and it can be shown following the proof of
\cref{prop:4_conditions_convex_polyhedral_case}\,\ref{item:suff_cond_pseudo_subreg_Gfrerer}.
First, as mentioned above, we automatically have $s \in \mathbb X$
with $w_s(u,0) \in \mathbf T(u)$ from $u \in \mathcal C^{1,2}(\bar x)$,
which can be added to \eqref{eq:Sec_order_optim_cond_smooth}.
Now, we are in the same situation as when proving
\cref{prop:4_conditions_convex_polyhedral_case}\,\ref{item:suff_cond_pseudo_subreg_Gfrerer},
but we have to work with 
\eqref{eq:Sec_order_optim_cond_smooth}
instead of \eqref{eq:2-reg_assumptions}.
From $u \in \mathcal C^{1,2}(\bar x)$ we also get $\nabla \varphi(\bar x) u \leq 0$, 
while $w_s(u,0) \in \mathbf T(u)$ and \cref{lem:w_s_versus_w_tilde_s}
yield $\nabla^2 \dual{y^*}{g}(\bar x)[u,u]\leq 0$,
and $\dual{z^*}{\nabla g(\bar x) u} \leq 0$ follows from
 $z^* \in \mathcal N_D(g(\bar x))$ and $\nabla g(\bar x) u \in \mathcal T_D(g(\bar x))$,
which is implicitly required due to $w_s(u,0)\in\mathbf T(u)$.
Thus, multiplying the essential equation of 
\eqref{eq:Sec_order_optim_cond_smooth}
by $u$, the three nonpositive terms sum up to zero, so they all must vanish.
Hence, the arguments which we used to prove 
\cref{prop:4_conditions_convex_polyhedral_case}\,\ref{item:suff_cond_pseudo_subreg_Gfrerer}
also work with \eqref{eq:2-reg_assumptions} replaced by \eqref{eq:Sec_order_optim_cond_smooth}.

\subsection{Applications}\label{sec:applications}

In this subsection, 
we highlight some aspects of our results from \cref{sec:constraint_mappings} in 
two popular settings of optimization theory. 
More precisely, we focus on the feasible regions of 
complementarity-constrained and nonlinear semidefinite problems.
As mentioned at the end of \cref{sec:constraint_mappings}, 
we do not obtain any new insights for standard nonlinear programs
as these can be reformulated with the aid of a constraint mapping
where the involved set is convex and polyhedral.
Hence, we do not specify our findings for this elementary setting for brevity of presentation
but refer the interested reader to \cite{Avakov1989,AvakovArutunovIzmailov2007}
where the associated mixed-order optimality conditions and constraint qualifications
are worked out.

\subsubsection{Mathematical programs with complementarity constraints}\label{sec:MPCCs}

Let us introduce
\[
	\mathcal C:=(\R_+\times\{0\})\cup(\{0\}\times\R_+),
\]
the so-called complementarity angle.
For twice continuously differentiable data functions $G,H\colon\mathbb X\to\R^m$ 
with components $G_1,\ldots,G_m\colon\mathbb X\to\R$ and $H_1,\ldots,H_m\colon\mathbb X\to\R$, 
we address the constraint region given by
\begin{equation}\label{eq:MPCC}\tag{MPCC}
	(G_i(x),H_i(x))\in\mathcal C\quad i\in I
\end{equation}
where $I:=\{1,\ldots,m\}$.
The latter is distinctive for so called 
\emph{mathematical programs with complementarity constraints}
which have been studied intensively throughout the last decades,
see e.g.\ \cite{LuoPangRalph1996,OutrataKocvaraZowe1998} for some classical references.
We observe that \eqref{eq:MPCC} can be formulated via a constraint map using 
$D:=\mathcal C^m$.
Note that standard inequality and equality constraints 
can be added without any difficulties
due to \cref{lem:product_rule_tangents_polyhedral_sets,lem:product_rule_graphical_derivative} 
when taking the findings from \cite{Avakov1989,AvakovArutunovIzmailov2007} into account.
Here, we omit them for brevity of presentation. 

Fix some feasible point $\bar x\in\mathbb X$ of \eqref{eq:MPCC}.
A critical direction $u\in\mathbb S_{\mathbb X}$ 
of the associated problem \eqref{eq:nonsmooth_problem} necessarily needs to satisfy
\begin{equation}\label{eq:linearization_cone_MPCC}
	\begin{aligned}
		\nabla G_i(\bar x) u&=0\quad& &i\in I^{0+}(\bar x),&\\
		\nabla H_i(\bar x) u&=0\quad&	&i\in I^{+0}(\bar x),&\\
		(\nabla G_i(\bar x) u,\nabla H_i(\bar x) u)&\in\mathcal C&&i\in I^{00}(\bar x),&
	\end{aligned}
\end{equation}
where we used the well-known index sets
\begin{align*}
	I^{0+}(\bar x)&:=\{i\in I\,|\,G_i(\bar x)=0,\,H_i(\bar x)>0\},\\
	I^{+0}(\bar x)&:=\{i\in I\,|\,G_i(\bar x)>0,\,H_i(\bar x)=0\},\\
	I^{00}(\bar x)&:=\{i\in I\,|\,G_i(\bar x)=0,\,H_i(\bar x)=0\}.
\end{align*}

	We start with an illustration of 
	\cref{Pro:M-stat_via_second_order}\,\ref{item:pseudo_stationarity_constraint_maps}.
	Thanks to \cref{rem:polyhedral_case_in_general_framework}, 
	we need to check the constraint qualifications
	\eqref{eq:CQ_pseudo_subregularity_II} and \eqref{eq:CQ_pseudo_subregularity_I_polyhedral_new},
	and these can be specified to the present setting with the aid of
	\cref{lem:some_properties_of_polyhedral_sets,lem:product_rule_tangents_polyhedral_sets,lem:product_rule_graphical_derivative}.
	For brevity of presentation, we abstain from a discussion of the
	case where critical directions of order $(1,2)$ are involved.
Based on the representation
\[
		\gph\mathcal N_{\mathcal C}
		=
		(\R_+\times\{0\}\times\{0\}\times\R)\cup(\{0\}\times\R_+\times\R\times\{0\})\cup(\{0\}\times\{0\}\times\R_-\times\R_-),
	\]
some elementary calculations show
\begin{equation}\label{eq:graphical_derivative_normals_compl}
		D\mathcal N_{\mathcal C}((a,b),(\mu,\nu))(v)
		=
		\begin{cases}
			\{0\}\times\R	&	a>0,\,b=\mu=0,\,v_2=0,\\
			\R\times\{0\}	&	a=\nu=0,\,b>0,\,v_1=0,\\
			\R^2			&	a=b=0,\,\mu,\nu<0,\,v=0,\\
			\{0\}\times\R	&	a=b=\mu=0,\,\nu<0,\,v_1>0,\,v_2=0,\\
			\R_-\times\R	&	a=b=\mu=0,\,\nu<0,\,v=0,\\
			\R\times\{0\}	&	a=b=\nu=0,\,\mu<0,\,v_1=0,\,v_2>0,\\
			\R\times\R_-	&	a=b=\nu=0,\,\mu<0,\,v=0,\\
			\{0\}\times\R	&	a=b=\mu=0,\,\nu>0,\,v_1\geq 0,\,v_2=0,\\
			\R\times\{0\}	&	a=b=\nu=0,\,\mu>0,\,v_1=0,\,v_2\geq 0,\\
			\{0\}\times\R	&	a=b=\mu=\nu=0,\,v_1>0,\,v_2=0,\\
			\R\times\{0\}	&	a=b=\mu=\nu=0,\,v_1=0,\,v_2>0,\\
			\mathcal N_{\mathcal C}(0)	&	a=b=\mu=\nu=0,\,v=0,\\
			\emptyset		&	\text{otherwise}
		\end{cases}
\end{equation}
for arbitrary $((a,b),(\mu,\nu))\in\gph\mathcal N_{\mathcal C}$ and $v\in\R^2$.
Consequently, for $u\in\mathbb S_{\mathbb X}$ 
satisfying \eqref{eq:linearization_cone_MPCC}, 
\eqref{eq:CQ_pseudo_subregularity_II} reduces to
\begin{equation}\label{eq:CQ_for_mixed_order_stationarity_MPCC}
	\left.
		\begin{aligned}
		&\nabla G(\bar x)^*\mu+\nabla H(\bar x)^*\nu=0,\\
		&\sum\nolimits_{i=1}^m\bigl(\mu_i\nabla^2G_i(\bar x)+\nu_i\nabla^2H_i(\bar x)\bigr)u
		\\
		&\qquad\qquad
		+\nabla G(\bar x)^*\tilde\mu+\nabla H(\bar x)^*\tilde\nu=0,\\
		&\forall i\in I^{+0}(\bar x)\cup I^{00}_{+0}(\bar x,u)\colon\,\mu_i=0,\\
		&\forall i\in I^{0+}(\bar x)\cup I^{00}_{0+}(\bar x,u)\colon\,\nu_i=0,\\
		&\forall i\in I^{00}_{00}(\bar x,u)\colon\,\mu_i,\nu_i\leq 0\,\text{ or }\,\mu_i\nu_i=0,\\
		&\forall i\in I\colon\,(\tilde\mu_i,\tilde\nu_i)
		\in D\mathcal N_{\mathcal C}
		((\bar G_i,\bar H_i),(\mu_i,\nu_i))(\nabla\bar G_iu,\nabla \bar H_iu)
		\end{aligned}
	\right\}
	\quad
	\Longrightarrow
	\quad
	\mu=\nu=0,
\end{equation}
while \eqref{eq:CQ_pseudo_subregularity_I_polyhedral_new} reads as
\begin{equation}\label{eq:CQ_for_mixed_order_stationarity_MPCC_new}
	\left.
		\begin{aligned}
		&\nabla G(\bar x)^*\mu+\nabla H(\bar x)^*\nu=0,\,
		\nabla G(\bar x)^*\tilde\mu+\nabla H(\bar x)^*\tilde\nu=0,\\
		&\forall i\in I^{+0}(\bar x)\cup I^{00}_{+0}(\bar x,u)\colon\,\mu_i=0,\\
		&\forall i\in I^{0+}(\bar x)\cup I^{00}_{0+}(\bar x,u)\colon\,\nu_i=0,\\
		&\forall i\in I^{00}_{00}(\bar x,u)\colon\,\mu_i,\nu_i\leq 0\,\text{ or }\,\mu_i\nu_i=0,\\
		&\forall i\in I\colon\,(\tilde\mu_i,\tilde\nu_i)
		\in D\mathcal N_{\mathcal C}
		((\bar G_i,\bar H_i),(\mu_i,\nu_i))(\nabla\bar G_iu,\nabla \bar H_iu)
		\end{aligned}
	\right\}
	\quad
	\Longrightarrow
	\quad
	\tilde\mu=\tilde\nu=0.
\end{equation}
Above, for each $i\in I$, we used $\bar G_i:=G_i(\bar x)$, $\bar H_i:=H_i(\bar x)$, 
$\nabla\bar G_iu:=\nabla G_i(\bar x) u$, and
$\nabla\bar H_iu:=\nabla H_i(\bar x) u$ for brevity as well as the index sets
\begin{align*}
	I^{00}_{0+}(\bar x,u)&:=\{i\in I^{00}(\bar x)\,|\,\nabla \bar G_iu=0,\,\nabla \bar H_iu>0\},\\
	I^{00}_{+0}(\bar x,u)&:=\{i\in I^{00}(\bar x)\,|\,\nabla \bar G_iu>0,\,\nabla \bar H_iu=0\},\\
	I^{00}_{00}(\bar x,u)&:=\{i\in I^{00}(\bar x)\,|\,\nabla \bar G_iu=0,\,\nabla \bar H_iu=0\}.
\end{align*}
The first assertion of \cref{Pro:M-stat_via_second_order}\,\ref{item:pseudo_stationarity_constraint_maps} now yields that whenever 
$\bar x$ is a local minimizer for the associated problem \eqref{eq:nonsmooth_problem} and for each
$u\in\mathbb S_{\mathbb X}$ satisfying \eqref{eq:linearization_cone_MPCC}, 
\eqref{eq:CQ_for_mixed_order_stationarity_MPCC} and \eqref{eq:CQ_for_mixed_order_stationarity_MPCC_new} hold,
then $\bar x$ is either M-stationary, i.e., there are multipliers $\mu,\nu\in\R^m$ satisfying
\begin{align*}
	&0\in\partial \varphi(\bar x)+\nabla G(\bar x)^*\mu+\nabla H(\bar x)^*\nu,\\
	&\forall i\in I^{+0}(\bar x)\colon\,\mu_i=0,\\
	&\forall i\in I^{0+}(\bar x)\colon\,\nu_i=0,\\
	&\forall i\in I^{00}(\bar x)\colon\,\mu_i,\nu_i\leq 0\,\text{ or }\,\mu_i\nu_i=0,
\end{align*}
or we find $u\in\mathbb S_{\mathbb X}$ satisfying \eqref{eq:linearization_cone_MPCC} 
and $\mathrm d\varphi(\bar x)(u)\leq 0$ 
as well as multipliers $\mu,\nu,\tilde\mu,\tilde\nu\in\R^m$ such that
\begin{equation}\label{eq:mixed_order_stationarity_MPCC}
	\begin{aligned}
	&0\in\partial \varphi(\bar x;u)
		+\sum\nolimits_{i=1}^m\bigl(\mu_i\nabla^2G_i(\bar x)+\nu_i\nabla^2H_i(\bar x)\bigr)u
		+\nabla G(\bar x)^*\tilde\mu+\nabla H(\bar x)^*\tilde\nu,\\
	&0=\nabla G(\bar x)^*\mu+\nabla H(\bar x)^*\nu,\\
	&\forall i\in I^{+0}(\bar x)\cup I^{00}_{+0}(\bar x,u)\colon\,\mu_i=0,\\
	&\forall i\in I^{0+}(\bar x)\cup I^{00}_{0+}(\bar x,u)\colon\,\nu_i=0,\\
	&\forall i\in I^{00}_{00}(\bar x,u)\colon\,\mu_i,\nu_i\leq 0\,\text{ or }\,\mu_i\nu_i=0,\\
	&\forall i\in I\colon\,(\tilde\mu_i,\tilde\nu_i)
	\in D\mathcal N_{\mathcal C}
	((\bar G_i,\bar H_i),(\mu_i,\nu_i))(\nabla\bar G_iu,\nabla \bar H_iu).
	\end{aligned}
\end{equation}

For brevity, we present the results from 
\cref{Pro:M-stat_via_second_order}\,\ref{item:pseudo_stationarity_constraint_maps_polyhedral_refined}
only in simplified from, where $w_s(u,v)$ is replaced by $0$, see \cref{rem:refined_polyhedral_situation_normal_cone_relation} as well,
and we do not comment on the cases where critical directions of order $(1,2)$ are involved,
but this would clearly yield further refinements.

In order to characterize \eqref{eq:CQ_pseudo_subregularity_polyhedral_II},
we observe that
\begin{align*}
	\mathcal N_{\mathcal T_{\mathcal C}(\bar G_i,\bar H_i)}(\nabla\bar G_iu,\nabla\bar H_iu)
	=
	\begin{cases}
		\{0\}\times\R	&i\in I^{+0}(\bar x)\cup I^{00}_{+0}(\bar x,u),\\
		\R\times\{0\}	&i\in I^{0+}(\bar x)\cup I^{00}_{0+}(\bar x,u),\\
		\mathcal N_{\mathcal C}(0)	&i\in I^{00}_{00}(\bar x,u)
	\end{cases}
\end{align*}
is valid for each $i\in I$. 
For each pair $(\mu_i,\nu_i)\in\mathcal N_{\mathcal T_{\mathcal C}(\bar G_i,\bar H_i)}(\nabla\bar G_iu,\nabla\bar H_iu)$,
elementary calculations 
and a comparison with \eqref{eq:graphical_derivative_normals_compl} show
\begin{align*}
	\mathcal T_{\mathcal N_{\mathcal T_{\mathcal C}(\bar G_i,\bar H_i)}(\nabla\bar G_iu,\nabla\bar H_iu)}(\mu_i,\nu_i)
	&=
	\begin{cases}
		\{0\}\times\R	&i\in I^{+0}(\bar x)\cup I^{00}_{+0}(\bar x,u),\\
		\R\times\{0\}	&i\in I^{0+}(\bar x)\cup I^{00}_{0+}(\bar x,u),\\
		\R^2			&i\in I^{00}_{00}(\bar x,u),\,\mu_i<0,\,\nu_i<0,\\
		\R_-\times\R	&i\in I^{00}_{00}(\bar x,u),\,\mu_i=0,\,\nu_i<0,\\
		\R\times\R_-	&i\in I^{00}_{00}(\bar x,u),\,\mu_i<0,\,\nu_i=0,\\
		\{0\}\times\R	&i\in I^{00}_{00}(\bar x,u),\,\mu_i=0,\,\nu_i>0,\\
		\R\times\{0\}	&i\in I^{00}_{00}(\bar x,u),\,\mu_i>0,\,\nu_i=0,\\
		\mathcal N_{\mathcal C}(0)	&i\in I^{00}_{00}(\bar x,u),\,\mu_i=\nu_i=0
	\end{cases}
	\\
	&=
	D\mathcal N_{\mathcal C}((\bar G_i,\bar H_i),(\mu_i,\nu_i))(\nabla \bar G_i,\nabla \bar H_i).
\end{align*}
Thus, validity of \eqref{eq:CQ_for_mixed_order_stationarity_MPCC}
for each $u\in\mathbb S_{\mathbb X}$ satisfying \eqref{eq:linearization_cone_MPCC} 
is already enough to infer that whenever $\bar x$ is a local minimizer,
then it is either M-stationary or there are $u\in\mathbb S_{\mathbb X}$ satisfying \eqref{eq:linearization_cone_MPCC} 
as well as $\mathrm d\varphi(\bar x)(u)\leq 0$
and multipliers $\mu,\nu,\tilde\mu,\tilde\nu\in\R^m$ 
solving the stationarity conditions \eqref{eq:mixed_order_stationarity_MPCC}.

Let us further note that \cref{Pro:M-stat_via_second_order}\,\ref{item:pseudo_stationarity_constraint_maps_polyhedral_refined}
also allows for the consideration of a qualification and stationarity condition where simply
$(\tilde\mu_i,\tilde\nu_i)\in\mathcal N_{\mathcal T_{\mathcal C}(\bar G_i,\bar H_i)}(\nabla\bar G_iu,\nabla\bar H_iu)$
has to hold for all $i\in I$, see \cref{rem:refined_polyhedral_situation_normal_cone_relation} again.
One can easily check that there is no general inclusion between
$\mathcal N_{\mathcal T_{\mathcal C}(\bar G_i,\bar H_i)}(\nabla\bar G_iu,\nabla\bar H_iu)$
and
$\mathcal T_{\mathcal N_{\mathcal T_{\mathcal C}(\bar G_i,\bar H_i)}(\nabla\bar G_iu,\nabla\bar H_iu)}(\mu_i,\nu_i)$,
i.e., this procedure leads to conditions not related to 
\eqref{eq:CQ_for_mixed_order_stationarity_MPCC} and \eqref{eq:mixed_order_stationarity_MPCC}
which are, however, easier to evaluate.
	
	The following example illustrates a situation where \eqref{eq:CQ_for_mixed_order_stationarity_MPCC} is valid while
	\eqref{eq:CQ_for_mixed_order_stationarity_MPCC_new} is violated, i.e., where
	\cref{Pro:M-stat_via_second_order}\,\ref{item:pseudo_stationarity_constraint_maps_polyhedral_refined} is applicable while
	\cref{Pro:M-stat_via_second_order}\,\ref{item:pseudo_stationarity_constraint_maps} is not.
	This provides yet another justification of a separate consideration of the polyhedral situation.

	\begin{example}\label{ex:MPCC}
		Let us consider \eqref{eq:MPCC} with $\mathbb X:=\R$, $m:=1$, 
		and $G(x):=x$ as well as $H(x):=x^2$ for all $x\in\R$.
		We are interested in the unique feasible point $\bar x:=0$ of this system.
		The only direction from the unit sphere that satisfies \eqref{eq:linearization_cone_MPCC}
		is $u:=1$. Hence, \eqref{eq:CQ_for_mixed_order_stationarity_MPCC} reduces to
		\[
		\left.
		\begin{aligned}
		&\mu=0,\,2\nu+\tilde\mu=0,\\
		&(\tilde\mu,\tilde\nu)
		\in D\mathcal N_{\mathcal C}
		((0,0),(\mu,\nu))(1,0)
		\end{aligned}
		\right\}
		\quad
		\Longrightarrow
		\quad
		\mu=\nu=0.
		\]
		Let the premise be valid and assume $\nu\neq 0$.
		This gives $D\mathcal N_{\mathcal C}((0,0),(0,\nu))(1,0)=\{0\}\times\R$ due to
		\eqref{eq:graphical_derivative_normals_compl},
		i.e., $\tilde\mu=0$, and, thus, $\nu=0$ which yields a contradiction.
		Hence, this constraint qualification holds.
		However, \eqref{eq:CQ_for_mixed_order_stationarity_MPCC_new} is given by
		\[
		\left.
		\begin{aligned}
		&\mu=0,\,\tilde\mu=0,\\
		&(\tilde\mu,\tilde\nu)
		\in D\mathcal N_{\mathcal C}
		((0,0),(\mu,\nu))(1,0)
		\end{aligned}
		\right\}
		\quad
		\Longrightarrow
		\quad
		\tilde\mu=\tilde\nu=0,
		\]
		and one can easily check with the aid of
		\eqref{eq:graphical_derivative_normals_compl}
		that the premise holds for $(\mu,\nu):=(\tilde\mu,\tilde\nu):=(0,1)$,
		i.e., this condition is violated.
	\end{example}

	Finally, we would like to refer the interested reader to \cite[Section~6]{IzmailovSolodov2002} 
	and \cite{IzmailovSolodov2002b} where the theory of
	$2$--regularity is first extended to mappings which are once but not twice differentiable and then
	applied to a suitable reformulation of complementarity constraints 
		as a system of once but not twice differentiable equations.
		We abstain from a detailed comparison of our findings with the ones 
		from \cite{IzmailovSolodov2002,IzmailovSolodov2002b}
		for the following reasons.
		First, in these papers, a different way of stating the system of complementarity constraints is used, 
		and it would be laborious to transfer the results to the formulation \eqref{eq:MPCC}.
		Second, at least in \cite{IzmailovSolodov2002}, some additional assumptions are used to simplify the
		calculations while we do not need to assume anything artificial to make the calculus accessible.
		Third, the final characterization of 2-regularity obtained in these papers does
		not comprise any second-order derivatives of the involved data functions and, thus,
		is anyhow clearly different from \eqref{eq:CQ_for_mixed_order_stationarity_MPCC}.
		Let us, however, close with the remark that the system of necessary optimality conditions derived
		in \cite[Theorem~4.2]{IzmailovSolodov2002b} is closely related to \eqref{eq:mixed_order_stationarity_MPCC}. 

\subsubsection{Semidefinite programming}

Let us consider the Hilbert space $\SSS_m$ of all real symmetric matrices equipped with the
standard (Frobenius) inner product. We denote by $\SSS_m^+$ and $\SSS_m^-$ the cone of all
positive and negative semidefinite matrices in $\SSS_m$, respectively.
The foundations of variational analysis in this space can be found, 
e.g., in \cite[Section~5.3]{BonnansShapiro2000}.
For some twice continuously differentiable mapping $g\colon\mathbb X\to\SSS_m$, 
we investigate the
constraint system
\begin{equation}\label{eq:semidefinite_NLP}\tag{NLSD}
	g(x)\in\SSS_m^+.
\end{equation} 
It is well known that the closed, convex cone $\SSS_m^+$ is not polyhedral.
Nevertheless, the constraint \eqref{eq:semidefinite_NLP}, 
associated with so-called \emph{nonlinear semidefinite programming}, 
can be encoded via a constraint map.
	Subsequently, we merely illustrate the first assertion of
	\cref{Pro:M-stat_via_second_order}\,\ref{item:pseudo_stationarity_constraint_maps}.
	As $\SSS_m^+$ is not polyhedral, \cref{Lem:Crit_dir_constraints} cannot
	be used for a characterization of critical directions of order $(1,2)$.

Let $\bar x\in\mathbb X$ be feasible to \eqref{eq:semidefinite_NLP} and, 
for some $u\in\mathbb S_{\mathbb X}$, fix
$\Omega\in\mathcal N_{\SSS_m^+}(g(\bar x);\nabla g(\bar x)u)$. 
For later use, fix an orthogonal matrix 
$\mathbf P\in\R^{m\times m}$ and a diagonal matrix $\mathbf\Lambda\in\R^{m\times m}$
whose diagonal elements $\lambda_1,\ldots,\lambda_m$ are ordered nonincreasingly 
such that $g(\bar x)+\Omega=\mathbf P\mathbf\Lambda\mathbf P^\top$. The index sets
corresponding to the positive, zero, and negative entries on the main diagonal 
of $\mathbf\Lambda$ are denoted by 
$\alpha$, $\beta$, and $\gamma$, respectively.
We emphasize that, here and throughout the subsection,
$\alpha$ is a constant index set while
$\beta$ and $\gamma$ depend on the precise choice of $\Omega$.
Subsequently, we use 
$\mathbf Q^{\mathbf P}:=\mathbf P^\top \mathbf Q\mathbf P$ and 
$\mathbf Q^{\mathbf P}_{IJ}:=(\mathbf Q^{\mathbf P})_{IJ}$
for each matrix $\mathbf Q\in\mathcal S_m$ 
and index sets $I,J\subset\{1,\ldots,m\}$
where $\mathbf M_{IJ}$ is the submatrix of $\mathbf M\in\mathcal S_m$
which possesses only those rows and columns of $\mathbf M$ whose
indices can be found in $I$ and $J$, respectively. 

The above constructions yield
\[
	g(\bar x)=\mathbf P\max(\mathbf\Lambda,\mathbf O)\mathbf P^\top,\qquad
	\Omega=\mathbf P\min(\mathbf\Lambda,\mathbf O)\mathbf P^\top
\]
where $\max$ and $\min$ have to be understood in entrywise fashion and $\OOO$ is an all-zero matrix of appropriate dimensions.
Due to
\[
	\nabla g(\bar x)u\in\mathcal T_{\mathcal S_m^+}(g(\bar x))
	=
	\left\{
		\QQQ\in\mathcal S_m\,\middle|\,
		\QQQ^{\PPP}_{\beta\cup\gamma,\beta\cup\gamma}\in\mathcal S_{|\beta\cup\gamma|}^+
	\right\},
\]
we find
\begin{align*}
	0
	&=
	\innerprod{\Omega}{\nabla g(\bar x)u}
	=
	\trace(\Omega\,\nabla g(\bar x)u)
	=
	\trace(\PPP\min(\LLL,\mathbf O)\PPP^\top \PPP[\nabla g(\bar x)u]^{\PPP}\PPP^\top)
	\\
	&=
	\trace(\min(\LLL,\mathbf O)[\nabla g(\bar x)u]^{\PPP})
	=
	\sum_{i\in\gamma}\underbrace{\lambda_i}_{<0}\,\underbrace{[\nabla g(\bar x)u]^{\PPP}_{i,i}}_{\geq 0}
\end{align*}
which directly gives us $[\nabla g(\bar x)u]^\PPP_{\beta\gamma}=\OOO$, $[\nabla g(\bar x)u]^\PPP_{\gamma\gamma}=\OOO$, and $[\nabla g(\bar x)u]^\PPP_{\beta\beta}\in\SSS_{|\beta|}^+$.
Furthermore, we note
\[
	\mathcal N_{\SSS^+_m}(g(\bar x))
	=
	\left\{
		\tilde\Omega\in\SSS_m\,\middle|\,
		\tilde\Omega^\PPP_{\alpha\alpha}=\OOO,\,\tilde\Omega^\PPP_{\alpha\beta}=\OOO,\,\tilde\Omega^{\PPP}_{\alpha\gamma}=\OOO,\,
		\tilde\Omega^\PPP_{\beta\cup\gamma,\beta\cup\gamma}\in\SSS^-_{|\beta\cup\gamma|}
	\right\}.
\]
Finally, let $\Xi_{\alpha\gamma}\in\R^{|\alpha|\times|\gamma|}$ be the matrix given by
\[
	\forall i\in \alpha\,\forall j\in\gamma\colon\quad
	[\Xi_{\alpha\gamma}]_{ij}:=-\frac{\lambda_j}{\lambda_i}.
\]
It is well known that the projection onto $\SSS_m^+$ is directionally
differentiable.
With the aid of \cref{lem:graphical_derivatives_of_normal_cone_map} and \cite[Corollary~3.1]{WuZhangZhang2014}, we find
\begin{align*}
	D\mathcal N_{\SSS_m^+}(g(\bar x),\Omega)(\nabla g(\bar x)u)
	&=
	\left\{\tilde\Omega\in\SSS_m\,\middle|\,
		\begin{aligned}
			&\tilde\Omega^\PPP_{\alpha\alpha}=\OOO,\,\tilde\Omega^\PPP_{\alpha\beta}=\OOO,\,
			\tilde\Omega^\PPP_{\alpha\gamma}=\Xi_{\alpha\gamma}\bullet[\nabla g(\bar x)u]^\PPP_{\alpha\gamma},\\
			&\tilde\Omega^\PPP_{\beta\beta}\in\SSS_{|\beta|}^-,\,\innerprod{\tilde\Omega^\PPP_{\beta\beta}}{[\nabla g(\bar x)u]^\PPP_{\beta\beta}}=0
		\end{aligned}
	\right\},
\end{align*}
and if $\nabla g(\bar x)u\neq\OOO$, we obtain
\begin{align*}
	D_\textup{sub}\mathcal N_{\SSS_m^+}(g(\bar x),\Omega)\left(\frac{\nabla g(\bar x)u}{\norm{\nabla g(\bar x)u}}\right)
	&\subset
	\left\{\tilde\Omega\in\SSS_m\,\middle|\,
		\begin{aligned}
			&\tilde\Omega^\PPP_{\alpha\alpha}=\OOO,\,\tilde\Omega^\PPP_{\alpha\beta}=\OOO,\,
			\tilde\Omega^\PPP_{\alpha\gamma}=\OOO,\\
			&\tilde\Omega^\PPP_{\beta\beta}\in\SSS_{|\beta|}^-,\,
			\innerprod{\tilde\Omega^\PPP_{\beta\beta}}{[\nabla g(\bar x)u]^\PPP_{\beta\beta}}=0
		\end{aligned}
	\right\}.
\end{align*}
Above, $\bullet$ represents the \emph{Hadamard}, i.e., entrywise product.
Note that validity of the final orthogonality condition in the estimate 
for the graphical subderivative follows from
\cref{lem:graphical_derivatives_of_normal_cone_map} since 
$\tilde\Omega\in D_\textup{sub}\mathcal N_{\SSS_m^+}(g(\bar x),\Omega)(\nabla g(\bar x)u/\norm{\nabla g(\bar x)u})$ 
and $\norm{\nabla g(\bar x)u}>0$ yield
\begin{align*}
	0
	&\leq
	\ninnerprod{\tilde\Omega}{\nabla g(\bar x)u}
	=
	\trace\bigl(\tilde\Omega\,\nabla g(\bar x)u\bigr)
	=
	\trace\bigl(\PPP\tilde\Omega^{\PPP}\PPP^\top\PPP[\nabla g(\bar x)u]^\PPP\PPP^\top\bigr)
	\\
	&=
	\trace\bigl(\tilde\Omega^\PPP[\nabla g(\bar x)u]^\PPP\bigr)
	=
	\trace\bigl(\tilde\Omega^\PPP_{\beta\beta}[\nabla g(\bar x)u]^\PPP_{\beta\beta}\bigr)
	\leq 
	0
\end{align*}
due to $\tilde\Omega^\PPP_{\alpha\alpha}=\OOO$, $\tilde\Omega^\PPP_{\alpha\beta}=\OOO$, $\tilde\Omega^\PPP_{\alpha\gamma}=\OOO$,
$\tilde\Omega^\PPP_{\beta\beta}\in\SSS_{|\beta|}^-$, $[\nabla g(\bar x)u]^\PPP_{\beta\gamma}=\OOO$, $[\nabla g(\bar x)u]^\PPP_{\gamma\gamma}=\OOO$,
and $[\nabla g(\bar x)u]^\PPP_{\beta\beta}\in\SSS_{|\beta|}^+$.
Thus, for each $u\in\mathbb S_{\mathbb X}$, \eqref{eq:CQ_pseudo_subregularity_II} takes the form
\[
	\left.
		\begin{aligned}
			\nabla g(\bar x)^*\Omega=0,\,
			\nabla^2\innerprod{\Omega}{g}(\bar x)(u)+\nabla g(\bar x)^*\tilde\Omega=0,\\
			\Omega^\PPP_{\alpha\alpha}=\OOO,\,\Omega^\PPP_{\alpha\beta}=\OOO,\,
			\Omega^{\PPP}_{\alpha\gamma}=\OOO,\,
			\Omega^\PPP_{\beta\cup\gamma,\beta\cup\gamma}\in\SSS_{|\beta\cup\gamma|}^-,\\
			[\nabla g(\bar x)u]^\PPP_{\beta\gamma}=\OOO,\,
			[\nabla g(\bar x)u]^\PPP_{\gamma\gamma}=\OOO,\,
			[\nabla g(\bar x)u]^\PPP_{\beta\beta}\in\SSS_{|\beta|}^+,\\
			\tilde\Omega^\PPP_{\alpha\alpha}=\OOO,\,\tilde\Omega^\PPP_{\alpha\beta}=\OOO,\,
			\tilde\Omega^\PPP_{\alpha\gamma}
			=
			\Xi_{\alpha\gamma}\bullet[\nabla g(\bar x)u]^\PPP_{\alpha\gamma},\\
			\tilde\Omega^\PPP_{\beta\beta}\in\SSS_{|\beta|}^-,\,
			\innerprod{\tilde\Omega^\PPP_{\beta\beta}}{[\nabla g(\bar x)u]^\PPP_{\beta\beta}}=0
		\end{aligned}
	\right\}
	\quad
	\Longrightarrow\quad
	\Omega=\OOO,
\]
while \eqref{eq:CQ_pseudo_subregularity_Ia} and \eqref{eq:CQ_pseudo_subregularity_Ib} 
(the latter in the case $\nabla g(\bar x)u\neq\OOO$) are both implied by
\[
	\left.
		\begin{aligned}
			\nabla g(\bar x)^*\Omega=0,\,\nabla g(\bar x)^*\tilde\Omega=0,\\
			\Omega^\PPP_{\alpha\alpha}=\OOO,\,\Omega^\PPP_{\alpha\beta}=\OOO,\,
			\Omega^{\PPP}_{\alpha\gamma}=\OOO,\,
			\Omega^\PPP_{\beta\cup\gamma,\beta\cup\gamma}\in\SSS_{|\beta\cup\gamma|}^-,\\
			[\nabla g(\bar x)u]^\PPP_{\beta\gamma}=\OOO,\,
			[\nabla g(\bar x)u]^\PPP_{\gamma\gamma}=\OOO,\,
			[\nabla g(\bar x)u]^\PPP_{\beta\beta}\in\SSS_{|\beta|}^+,\\
			\tilde\Omega^\PPP_{\alpha\alpha}=\OOO,\,\tilde\Omega^\PPP_{\alpha\beta}=\OOO,\,
			\tilde\Omega^\PPP_{\alpha\gamma}=\OOO,\\
			\tilde\Omega^\PPP_{\beta\beta}\in\SSS_{|\beta|}^-,\,
			\innerprod{\tilde\Omega^\PPP_{\beta\beta}}{[\nabla g(\bar x)u]^\PPP_{\beta\beta}}=0
		\end{aligned}
	\right\}
	\quad
	\Longrightarrow\quad
	\tilde\Omega=\OOO.
\]
In the case where $\bar x$ is a local minimizer of the associated problem \eqref{eq:nonsmooth_problem}, 
validity of these conditions for each $u\in\mathbb S_{\mathbb X}$ guarantees that
$\bar x$ is either M-stationary (we omit stating this well-known system here) 
or we find $u\in\mathbb S_{\mathbb X}$ 
and $\Omega,\tilde\Omega\in\SSS_m$ such that
\begin{align*}
			&0\in\partial \varphi(\bar x;u)+\nabla^2\innerprod{\Omega}{g}(\bar x)(u)
			+\nabla g(\bar x)^*\tilde\Omega,\,\nabla g(\bar x)^*\Omega=0,\\
			&\Omega^\PPP_{\alpha\alpha}=\OOO,\,\Omega^\PPP_{\alpha\beta}=\OOO,\,
			\Omega^{\PPP}_{\alpha\gamma}=\OOO,\,
			\Omega^\PPP_{\beta\cup\gamma,\beta\cup\gamma}\in\SSS_{|\beta\cup\gamma|}^-,\\
			&\mathrm d\varphi(\bar x)(u)\leq 0,\,
			[\nabla g(\bar x)u]^\PPP_{\beta\gamma}=\OOO,\,
			[\nabla g(\bar x)u]^\PPP_{\gamma\gamma}=\OOO,\,
			[\nabla g(\bar x)u]^\PPP_{\beta\beta}\in\SSS_{|\beta|}^+,\\
			&\tilde\Omega^\PPP_{\alpha\alpha}=\OOO,\,
			\tilde\Omega^\PPP_{\alpha\beta}=\OOO,\,
			\tilde\Omega^\PPP_{\alpha\gamma}
			=\Xi_{\alpha\gamma}\bullet[\nabla g(\bar x)u]^\PPP_{\alpha\gamma},\\
			&\tilde\Omega^\PPP_{\beta\beta}\in\SSS_{|\beta|}^+,\,
			\innerprod{\tilde\Omega^\PPP_{\beta\beta}}{[\nabla g(\bar x)u]^\PPP_{\beta\beta}}=0.
\end{align*}

\section{Directional asymptotic regularity in nonsmooth optimization}\label{sec:directional_asymptotic_regularity}

In this section, we focus on (directional) asymptotic regularity conditions,
which essentially correspond to conditions ensuring that
(directional) asymptotic stationarity 
	from \cref{thm:directional_asymptotic_stationarity},
	which serves as a necessary optimality condition for \eqref{eq:nonsmooth_problem}
	even in the absence of constraint qualifications,
translates into M-stationarity.
We provide a comprehensive comparison of (directional) asymptotic regularity with various
known constraint qualifications.
Throughout the section, we consider a set-valued mapping 
$\Phi\colon\mathbb X\tto\mathbb Y$ with a closed graph.

\subsection{On the concept of directional asymptotic regularity}\label{sec:directional_asymptotic_regularity_basics}

Based on \cref{thm:directional_asymptotic_stationarity}, the following definition introduces
concepts which may serve as (directional) qualification conditions for the model problem \eqref{eq:nonsmooth_problem}.
\begin{definition}\label{def:asymptotic_regularity}
	Let $(\bar x,\bar y)\in\gph\Phi$ be fixed.
	\begin{enumerate}
		\item\label{item:def_asymp_reg}
			The map $\Phi$ is said to be \emph{asymptotically regular at $(\bar x,\bar y)$}
			whenever the following condition holds:
			for every sequences $\{(x_k,y_k)\}_{k\in\N}\subset\gph\Phi$,
			$\{x_k^*\}_{k\in\N}\subset\mathbb X$, and $\{\lambda_k\}_{k\in\N}\subset\mathbb Y$
			as well as $x^*\in\mathbb X$ satisfying $x_k\to\bar x$, $y_k\to\bar y$,
			$x_k^*\to x^*$, and $x_k^*\in\widehat{D}^*\Phi(x_k,y_k)(\lambda_k)$ for all
			$k\in\N$, we find $x^*\in\Im D^*\Phi(\bar x,\bar y)$.
		\item\label{item:def_dir_asymp_reg} 
			For the fixed direction $u\in\mathbb S_{\mathbb X}$, 
			$\Phi$ is said to be \emph{asymptotically regular 
			at $(\bar x,\bar y)$ in direction $u$} whenever
			the following condition holds:
			for every sequences $\{(x_k,y_k)\}_{k\in\N}\subset\gph\Phi$,
			$\{x_k^*\}_{k\in\N}\subset\mathbb X$, and $\{\lambda_k\}_{k\in\N}\subset\mathbb Y$ 
			as well as $x^*\in\mathbb X$ and $y^*\in\mathbb Y$ satisfying 
			$x_k\notin\Phi^{-1}(\bar y)$, $y_k\neq\bar y$, and 
			$x_k^*\in \widehat{D}^*\Phi(x_k,y_k)(\lambda_k)$
			for each $k\in\N$ as well as the convergences
			\begin{equation}\label{eq:convergences_directional_asymptotic_regularity}
				\begin{aligned}
				x_k&\to\bar x,&
				\qquad
				y_k&\to\bar y,&
				\qquad
				x_k^*&\to x^*,&
				\\
				\frac{x_k-\bar x}{\norm{x_k-\bar x}}&\to u,&
				\qquad
				\frac{y_k-\bar y}{\norm{x_k-\bar x}}&\to 0,&
				\qquad  
				&&
				\\
				\frac{\norm{y_k-\bar y}}{\norm{x_k-\bar x}}\lambda_k&\to y^*,& 
				\qquad
				\norm{\lambda_k}&\to\infty,&
				\qquad
				\frac{y_k-\bar y}{\norm{y_k-\bar y}}-\frac{\lambda_k}{\norm{\lambda_k}}
				&\to 0,&
				\end{aligned}
			\end{equation}
			we find $x^*\in\Im D^*\Phi(\bar x,\bar y)$.
		\item\label{item:def_dir_strong_asymp_reg}
			For the fixed direction $u\in\mathbb S_{\mathbb X}$, 
			$\Phi$ is said to be \emph{strongly asymptotically regular 
			at $(\bar x,\bar y)$ in direction $u$} whenever
			the following condition holds:
			for every sequences $\{(x_k,y_k)\}_{k\in\N}\subset\gph\Phi$,
			$\{x_k^*\}_{k\in\N}\subset\mathbb X$, and $\{\lambda_k\}_{k\in\N}\subset\mathbb Y$ 
			as well as $x^*\in\mathbb X$ and $y^*\in\mathbb Y$ satisfying 
			$x_k\notin\Phi^{-1}(\bar y)$, $y_k\neq\bar y$, 
			and $x_k^*\in \widehat{D}^*\Phi(x_k,y_k)(\lambda_k)$
			for each $k\in\N$ as well as the convergences
			\eqref{eq:convergences_directional_asymptotic_regularity}, 
			we have $x^*\in \Im D^*\Phi((\bar x,\bar y);(u,0))$.
	\end{enumerate}
\end{definition}

Before commenting in detail on these conditions, we would like to emphasize that they can be 
equivalently formulated in terms of limiting coderivatives completely.
The mainly technical proof of this result can be found in \cref{sec:appendix}.

\begin{proposition}\label{prop:asymptotic_regularity_via_limiting_tools}
	\Cref{def:asymptotic_regularity} can equivalently be formulated in terms of limiting normals.
\end{proposition}

Having \cref{prop:asymptotic_regularity_via_limiting_tools} available,
let us briefly note that asymptotic regularity of a set-valued mapping 
$\Phi\colon\mathbb X\tto\mathbb Y$ at some point $(\bar x,0)\in\gph\Phi$ in the sense
of \cref{def:asymptotic_regularity} equals AM-regularity of the set $\Phi^{-1}(0)$
at $\bar x$ mentioned in \cite[Remark~3.17]{Mehlitz2020a}.
The concepts of directional asymptotic regularity from 
\cref{def:asymptotic_regularity}\,\ref{item:def_dir_asymp_reg} and~\ref{item:def_dir_strong_asymp_reg} 
are new.

In the subsequent remark, we summarize some obvious relations between the different concepts
from \cref{def:asymptotic_regularity}.
\begin{remark}\label{rem:relations_asymptotic_regularity}
	Let $(\bar x,\bar y)\in\gph\Phi$ be fixed.
	Then the following assertions hold.
	\begin{enumerate}
		\item Let $u\in\mathbb S_{\mathbb X}$ be arbitrarily chosen.
			If $\Phi$ is strongly asymptotically regular at $(\bar x,\bar y)$ in direction
			$u$, it is asymptotically regular at $(\bar x,\bar y)$ in direction $u$.
		\item\label{item:dir_asymp_reg_equals_asymp_reg} 
			If $\Phi$ is asymptotically regular at $(\bar x,\bar y)$, then
			it is asymptotically regular at $(\bar x,\bar y)$ 
			in each direction from $\mathbb S_{\mathbb X}$.
	\end{enumerate}
\end{remark}

We note that strong asymptotic regularity in each unit direction is indeed not related to asymptotic regularity. 
On the one hand, the subsequently stated example, taken from \cite[Example~3.15]{Mehlitz2020a},
shows that asymptotic regularity does not imply strong asymptotic regularity in each unit direction.
On the other hand, \cref{ex:FOSCMS_but_not_asymptotically_regular} from below
illustrates that strong asymptotic regularity in each unit direction does not
yield asymptotic regularity.
\begin{example}\label{ex:asymptotic_regularity_but_not_strong}
	We consider $\Phi\colon\R\tto\R$ given by
	\[
		\forall x\in\R\colon\quad
		\Phi(x)
		:=
		\begin{cases}
			\R		&\text{if }x\leq 0,\\
			[x^2,\infty)	&\text{if }x>0
		\end{cases}
	\]
	at $(\bar x,\bar y):=(0,0)$. 
	It is demonstrated in \cite[Example~3.15]{Mehlitz2020a} that $\Phi$ is 
	asymptotically regular at $(\bar x,\bar y)$.
	We find $\mathcal T_{\gph\Phi}(\bar x,\bar y)=\{(u,v)\in\R^2\,|\,u\leq 0\,\text{ or }\,v\geq 0\}$
	so $(\pm 1,0)\in\mathcal T_{\gph\Phi}(\bar x,\bar y)$.
	Let us consider $u:=1$.
	Then we find $\Im D^*\Phi((\bar x,\bar y);(u,0))=\{0\}$.
	Taking $x^*:=1$, $y^*:=1/2$, as well as
	\[
		\forall k\in\N\colon\quad
		x_k:=\frac1k,
		\qquad
		y_k:=\frac1{k^2},
		\qquad
		x_k^*:=1,
		\qquad
		\lambda_k:=\frac k2,
	\]
	we have $x_k^*\in\widehat D^*\Phi(x_k,y_k)(\lambda_k)$ for all $k\in\N$
	as well as the convergences \eqref{eq:convergences_directional_asymptotic_regularity}.
	However, due to $x_k^*\to x^*\notin \Im D^*\Phi((\bar x,\bar y);(u,0))$, $\Phi$ is
	not strongly asymptotically regular at $(\bar x,\bar y)$ in direction $u$.
\end{example}

Combining \cref{thm:directional_asymptotic_stationarity} with the concepts from \cref{def:asymptotic_regularity}, 
we immediately obtain the following result which motivates our interest in directional asymptotic regularity.
\begin{corollary}\label{cor:M_stationarity_via_directional_asymptotic_regularity}
 	Let $\bar x\in\mathcal F$ be a local minimizer of \eqref{eq:nonsmooth_problem} 
 	such that, for each critical direction $u\in\mathbb S_{\mathbb X}$ 
 	for \eqref{eq:nonsmooth_problem} at $\bar x$, $\Phi$ is
	asymptotically regular at $(\bar x,\bar y)$ in direction $u$.
	Then $\bar x$ is M-stationary.
\end{corollary}
\begin{proof}
	Due to \cref{thm:directional_asymptotic_stationarity}, it suffices to consider the
	situation where there are a critical direction $u\in\mathbb S_{\mathbb X}$ for \eqref{eq:nonsmooth_problem}
	at $\bar x$ and $y^*\in\mathbb Y$ as well as sequences
	$\{x_k\}_{k\in\N},\{\eta_k\}_{k\in\N},\{x_k^*\}_{k\in\N}\subset\mathbb X$
	and  $\{y_k\}_{k\in\N}\subset\mathbb Y$ such that
	$x_k\notin\Phi^{-1}(\bar y)$, $y_k\neq\bar y$, $x_k^*\in\partial\varphi(x_k)$, and
	\begin{equation*}
		\eta_k-x_k^*\in D^*\Phi(x_k,y_k)\left(k(y_k-\bar y)\right)
	\end{equation*}
	for all $k\in\N$ as well as the convergences \eqref{eq:convergences_gamma=1} are valid.

	Since $\varphi$ is a locally Lipschitz continuous function, $\{x_k^*\}_{k\in\N}$ is bounded,
	see e.g.\ \cite[Theorem~1.22]{Mordukhovich2018}, and, thus, converges (along a subsequence),
	to some point $x^*\in\mathbb X$ which belongs to $\partial\varphi(\bar x)$
	by robustness of the limiting subdifferential.
	
	We can set $\lambda_k:=k(y_k-\bar y)$ for each $k\in\N$
	and obtain $\lambda_k\nnorm{y_k-\bar y}/\nnorm{x_k-\bar x}\to y^*$ and $\nnorm{\lambda_k}\to\infty$ 
	from \eqref{eq:convergences_gamma=1_multiplier} as well as
	$(y_k-\bar y)/\nnorm{y_k-\bar y}=\lambda_k/\nnorm{\lambda_k}$ for each $k\in\N$ by construction.
	Additionally, $\eta_k-x_k^*\in D^*\Phi(x_k,y_k)(\lambda_k)$ is valid for each $k\in\N$.
	
	Now, asymptotic regularity of $\Phi$ at $(\bar x,\bar y)$ in direction $u$,
	\cref{prop:asymptotic_regularity_via_limiting_tools}, and the
	remaining convergences from \eqref{eq:convergences_gamma=1} yield
	$-x^*\in\Im D^*\Phi(\bar x,\bar y)$, i.e., there exists $\lambda\in\mathbb Y$ such that
	$-x^*\in D^*\Phi(\bar x,\bar y)(\lambda)$.
	Recalling $x^*\in\partial\varphi(\bar x)$
	 shows the claim.
\end{proof}

In the light of 
\cref{rem:relations_asymptotic_regularity}\,\ref{item:dir_asymp_reg_equals_asymp_reg}, 
our result from
\cref{cor:M_stationarity_via_directional_asymptotic_regularity} improves
\cite[Theorem~3.9]{Mehlitz2020a} by a directional refinement of the constraint qualification
since it suffices to check asymptotic regularity w.r.t.\ particular directions.

We point out that, unlike typical constraint qualifications, (directional) asymptotic regularity
allows the existence of sequences 
satisfying \eqref{eq:convergences_directional_asymptotic_regularity}
as long as the limit $x^*$ is included in $\Im D^*\Phi(\bar x,\bar y)$ 
which is enough for M-stationarity.

\begin{remark}
	\cref{cor:M_stationarity_via_directional_asymptotic_regularity}
	requires asymptotic regularity in every (critical) unit direction.
	Taking into account \cref{Rem: 2_paths_CQ_vs_Crit_dir},
	we could also consider an alternative approach to secure M-stationarity,
	demanding either that there does not exist a critical direction
	together with the sequences from 
	\cref{def:asymptotic_regularity}\,\ref{item:def_dir_asymp_reg},
	or, in the case of existence, 
	that $\Phi$ is asymptotically regular at least in one of these critical directions.
	For brevity of presentation, we abstain from developing this approach further.
\end{remark}

Since (directional) asymptotic regularity (w.r.t.\ all critical unit directions) yields
M-stationarity of a local minimizer 
by \cref{cor:M_stationarity_via_directional_asymptotic_regularity},
in the remaining part of the paper, we put it into context of other
common assumptions that work as a constraint qualification
for M-stationarity associated with problem \eqref{eq:nonsmooth_problem}.
Let us clarify here some rather simple or known connections.
\begin{enumerate}
 \item\label{item:polyhedrality} 
 	A polyhedral mapping is asymptotically regular at each point of its graph.
 \item\label{item:metric_regularity} 
 	Metric regularity implies asymptotic regularity.
 \item\label{item:strong_metric_subregularity} 
 	Strong metric subregularity implies asymptotic regularity.
 \item\label{item:FOSCMS_does_not_imply_asymptotic_regularity} 
 	FOSCMS does not imply asymptotic regularity,
 	but it implies strong asymptotic regularity in each unit direction.
 \item\label{item:metric_subregularity_not_sufficient_for_asymptotic_regularity}
 	Metric subregularity does not imply asymptotic regularity in each unit direction.
 	However, if the map of interest is metrically subregular 
 	at every point of its graph near the reference point with a \emph{uniform} constant,
 	then strong asymptotic regularity in each unit direction follows.
 \item\label{item:asymptotic_regularity_does_not_yield_exact_penalty}
 	Neither asymptotic regularity nor strong directional asymptotic regularity yields
 	the directional exact penalty property of \cref{lem:directional_exact_penalization}.
\end{enumerate}

Statements~\ref{item:polyhedrality} and~\ref{item:metric_regularity} 
were shown in \cite[Theorems 3.10 and 3.12]{Mehlitz2020a}.
Let us now argue that strong metric subregularity (the ``inverse'' property associated with
isolated calmness), see \cite{DontchevRockafellar2014}, 
also implies asymptotic regularity at the point.
This follows easily from the discussion above
\cite[Corollary~4.6]{BenkoMehlitz2020}, which yields
that the domain of the limiting coderivative, at the point 
where the mapping is isolatedly calm, is the whole space.
Equivalently, the range of the limiting coderivative, at the point 
where the mapping is strongly metrically subregular, is the whole space 
and asymptotic regularity thus follows trivially.
Thus, statement~\ref{item:strong_metric_subregularity} follows.

Regarding~\ref{item:FOSCMS_does_not_imply_asymptotic_regularity}, the fact that FOSCMS implies
strong asymptotic regularity in each unit direction easily follows
by similar arguments 
that show that
metric regularity implies asymptotic regularity, 
see \cite[Lemma~3.11, Theorem~3.12]{Mehlitz2020a}.
Actually, it can be proved that validity of FOSCMS$(u)$ for some unit direction $u$
implies strong asymptotic regularity in direction $u$.
For constraint mappings, this also follows from \cref{cor:FOSCMS_SOSCMS_imply_sAR}
from below.

The following example shows that FOSCMS
does not imply asymptotic regularity.
\begin{example}\label{ex:FOSCMS_but_not_asymptotically_regular}
    Let $\Phi\colon \R \tto \R$ be given by
    \begin{equation*}
        \forall x\in\R\colon\quad
        \Phi(x):= 
        \begin{cases}
             [x,\infty) & \text{if } x \leq 0, \\
             \left[\frac{1}{k} - \frac{1}{k}\left(x - \frac{1}{k}\right),\infty\right)
             	& \text{if } x \in \left(\frac{1}{k+1},\frac{1}{k}\right] \text{ for some }k\in\N,\\
             \emptyset	&\text{otherwise.}
        \end{cases}
    \end{equation*}
    Then $\{(1/k,1/k)\}_{k\in\N}\subset\gph \Phi$ converges to
    $(\bar x,\bar y):=(0,0)$ and
    \begin{equation*}
        \mathcal N_{\gph \Phi}(1/k,1/k) =
        \{(x^*,y^*)\in\R^2 \,\vert\, y^* \leq 0, y^* \leq k x^*\}
    \end{equation*}
    is valid showing that $\Im D^*\Phi(1/k,1/k) = \R$ holds for all $k\in\N$.
    On the other hand, we have
    \begin{equation*}
        \mathcal N_{\gph \Phi}(0,0) =
        \{(x^*,y^*)\in\R^2 \,\vert\, x^* \geq 0, y^* \leq 0\},
    \end{equation*}
    and, thus, $\Im D^*\Phi(0,0) = \R_+$.
    This means that $\Phi$ is not asymptotically regular at $(\bar x,\bar y)$.
    
    On the other hand, we find
    \[
    	\mathcal T_{\gph\Phi}(\bar x,\bar y)=\{(u,v)\in\R^2\,|\,u\leq v\}.
    \]
    Each pair $(u,0)\in \mathcal T_{\gph\Phi}(\bar x,\bar y)$ with $u\neq 0$ satisfies
    $u<0$, i.e., the direction $(u,0)$ points into the interior of $\gph\Phi$. 
    Thus, we have $\mathcal N_{\gph\Phi}((\bar x,\bar y),(u,0))=\{(0,0)\}$
    which shows that FOSCMS is valid.
\end{example}

Regarding~\ref{item:metric_subregularity_not_sufficient_for_asymptotic_regularity}, 
let us fix $(\bar x,\bar y)\in\gph\Phi$ and note that metric subregularity of $\Phi$ 
on a neighborhood of $(\bar x,\bar y)$ (restricted to $\gph\Phi$) with a uniform constant
$\kappa>0$ is clearly milder
than metric regularity at $(\bar x,\bar y)$ 
since it is automatically satisfied, e.g., for polyhedral mappings.
To see that it implies asymptotic regularity, consider
sequences $\{(x_k,y_k)\}_{k\in\N}\subset\gph\Phi$, $\{x_k^*\}_{k\in\N}\subset\mathbb X$,
and $\{\lambda_k\}_{k\in\N}\subset\mathbb Y$ as well as $x^*\in\mathbb X$ 
and $y^*\in\mathbb Y$ satisfying 
$x_k^*\in \widehat{D}^*\Phi(x_k,y_k)(\lambda_k)$ for each $k\in\N$ 
and the convergences \eqref{eq:convergences_directional_asymptotic_regularity}
for some unit direction $u\in\mathbb S_{\mathbb X}$.
Due to \cite[Theorem~3.2]{BenkoMehlitz2020} and 
$-x_k^*\in \dom\widehat{D}^*\Phi^{-1}(y_k,x_k)$, we find
$x_k^*\in \widehat{\mathcal N}_{\Phi^{-1}(y_k)}(x_k) \subset \mathcal N_{\Phi^{-1}(y_k)}(x_k)$
for each $k\in\N$.
Furthermore, \cite[Theorem~3.2]{BenkoMehlitz2020} also gives the existence of 
$\tilde\lambda_k\in\mathbb Y$
with $\nnorm{\tilde\lambda_k} \leq \kappa \norm{x_k^*}$ 
and $x_k^*\in D^*\Phi(x_k,y_k)(\tilde\lambda_k)$.
Noting that $\{x_k^*\}_{k\in\N}$ converges, this shows that there is 
an accumulation point $\lambda\in\mathbb Y$ 
of $\{\tilde\lambda_k\}_{k\in\N}$ which satisfies 
$x^*\in D^*\Phi((\bar x,\bar y);(u,0))(\lambda)$
by robustness of the directional limiting coderivative, 
see \cref{lem:robustness_directional_limiting_normals}.
Hence, $\Phi$ is strongly asymptotically regular at $(\bar x,\bar y)$ in direction $u$.
Note that for the above arguments to work, we only need uniform metric subregularity 
along all sequences $\{(x_k,y_k)\}_{k\in\N}\subset\gph\Phi$
converging to $(\bar x, \bar y)$ from direction $(u,0)$.

The following example shows that
metric subregularity in the neighborhood of the point of interest
does not imply asymptotic regularity in each unit direction.

\begin{example}\label{ex:metric_subregularity_vs_asymptotic_regularity}
	We consider the mapping $\Phi\colon\R\tto\R$ given by
	\[
		\forall x\in\R\colon\quad
		\Phi(x):=\{0,x^2\}.
	\]
	Due to $\Phi^{-1}(0)=\R$, $\Phi$ is metrically subregular at all points $(x,0)$ where $x\in\R$
	is arbitrary.
	Furthermore, at all points $(x,x^2)$ where $x\neq 0$ holds, the Mordukhovich criterion
	\eqref{eq:Mordukhovich_criterion} shows that $\Phi$ is metrically regular. 
	Thus, $\Phi$ is metrically subregular at each point of its graph.
	Note that the moduli of metric subregularity tend to $\infty$ along the points 
	$(x,x^2)$ as $x\downarrow 0$ or $x\uparrow 0$.
	
	Let us consider the point $(\bar x,\bar y):=(0,0)$ where we have 
	$\mathcal N_{\gph\Phi}(\bar x,\bar y)=\{0\}\times\R$
	and, thus, $\Im D^*\Phi(\bar x,\bar y)=\{0\}$. 
	Choosing $x^*:=1$, $y^*:=1/2$, as well as
	\[
		\forall k\in\N\colon\quad
		x_k:=\frac1k,\qquad y_k:=\frac{1}{k^2},\qquad x_k^*:=1,\qquad \lambda_k:=\frac{k}{2},
	\]
	we have $x_k^*\in\widehat D^*\Phi(x_k,y_k)(\lambda_k)$ for all $k\in\N$ as well as
	the convergences \eqref{eq:convergences_directional_asymptotic_regularity} for $u:=1$.
	Due to $x_k^*\to x^*\notin\Im D^*\Phi(\bar x,\bar y)$, $\Phi$ is not asymptotically regular
	at $(\bar x,\bar y)$ in direction $u$.
\end{example}

Finally, let us address item~\ref{item:asymptotic_regularity_does_not_yield_exact_penalty}
with the aid of an example.
\begin{example}\label{ex:asymptotic_regularity_and_no_exact_penalty}
	Let us define $\varphi\colon\R\to\R$ and $\Phi\colon\R\tto\R$ by means of
	\[
		\forall x\in\R\colon\quad
		\varphi(x):=-x,
		\qquad
		\Phi(x)
		:=
		\begin{cases}
			\R	&\text{if }x\leq 0,\\
			[x^2,\infty)	&\text{if }x=\frac1k\text{ for some }k\in\N,\\
			\emptyset	&\text{otherwise.}
		\end{cases}
	\]
	Furthermore, we fix $\bar y:=0$.
	One can easily check that $\bar x:=0$ is the uniquely determined global
	minimizer of the associated problem \eqref{eq:nonsmooth_problem}.
	Furthermore, we have 
	$\Im D^*\Phi(\bar x,\bar y)=\Im D^*\Phi((\bar x,\bar y);(1,0))=\R$ which
	shows that $\Phi$ is asymptotically regular at $(\bar x,\bar y)$ as well
	as strongly asymptotically regular at $(\bar x,\bar y)$ in direction $1$.
	Furthermore, it is obvious that $\Phi$ is strongly asymptotically regular
	at $(\bar x,\bar y)$ in direction $-1$.
	Finally, let us mention that $\Phi$ fails to be metrically subregular at
	$(\bar x,\bar y)$ in direction $1$.
	
	Now, define $x_k:=1/k$ for each $k\in\N$ and observe that for each constant
	$C>0$ and sufficiently large $k\in\N$, we have 
	$\varphi(x_k)+C\,\dist(\bar y,\Phi(x_k))=-1/k+C/k^2<0=\varphi(\bar x)$, i.e., 
	$\bar x$ is not a minimizer 
	of \eqref{eq:directionally_penalized_problem}
	for any choice of $C>0$, $\varepsilon>0$, $\delta>0$, and $u:=1$.	
\end{example}

\subsection{Directional pseudo- and quasi-normality}\label{sec:pseudo_quasi_normality}

In this section, we connect asymptotic regularity with the notions of pseudo- and quasi-normality.
Note that the latter concepts have been introduced for
standard nonlinear programs in 
\cite{BertsekasOzdaglar2002,Hestenes1975},
and extensions to more general
geometric constraints have been established in \cite{GuoYeZhang2013}.
Furthermore, problem-tailored notions of these conditions have
been coined e.g.\ for so-called cardinality-, complementarity-, and switching-constrained
optimization problems, see \cite{KanzowRaharjaSchwartz2021b,KanzowSchwartz2010,LiangYe2021}.
Let us point out that these conditions are comparatively
mild constraint qualifications and sufficient for the
presence of metric subregularity of the associated constraint mapping,
see e.g.\ \cite[Theorem~5.2]{GuoYeZhang2013}.
Here, we extend pseudo- and quasi-normality
from the common setting of geometric constraint
systems to arbitrary set-valued mappings and comment
on the qualitative properties of these conditions.
Naturally, we aim for directional versions of these concepts, which,
in the setting of geometric constraints,
were recently introduced in \cite{BaiYeZhang2019} and further explored
in \cite{BenkoCervinkaHoheisel2019}.

\subsubsection{On the general concept of directional pseudo- and quasi-normality}

The definition below introduces the notions of our interest.

\begin{definition}\label{def:direction_quasi_pseudo_normality}
	Fix $(\bar x,\bar y)\in\gph\Phi$ and a direction $u\in\mathbb S_{\mathbb X}$.
	\begin{enumerate}
		\item We say that \emph{pseudo-normality in direction $u$} holds at $(\bar x,\bar y)$
			if there does not exist a nonzero vector 
			$\lambda\in\ker D^*\Phi((\bar x,\bar y);(u,0))$
			satisfying the following condition:
			there are sequences $\{(x_k,y_k)\}_{k\in\N}\subset\gph\Phi$ 
			with $x_k\neq\bar x$ for all $k\in\N$ and
			$\{\lambda_k\}_{k\in\N}\subset\mathbb Y$, $\{\eta_k\}_{k\in\N}\subset\mathbb X$,
			such that
			\begin{equation}\label{eq:convergences_definition_quasi_normality}
				\begin{aligned}
				x_k&\to\bar x,&
				\qquad
				y_k&\to\bar y,&
				\qquad
				\lambda_k&\to\lambda,&
				\\
				\eta_k&\to 0,&
				\qquad
				\frac{x_k-\bar x}{\norm{x_k-\bar x}}&\to u,&
				\qquad
				\frac{y_k-\bar y}{\norm{x_k-\bar x}}&\to 0,&
				\end{aligned}
			\end{equation}
			and $\eta_k\in \widehat{D}^*\Phi(x_k,y_k)(\lambda_k)$ 
			as well as $\dual{\lambda}{y_k-\bar y}>0$
			for all $k\in\N$.
		\item Let $\mathcal E:=\{e_1,\ldots,e_m\}\subset\mathbb Y$ be an orthonormal basis of
			$\mathbb Y$. We say that \emph{quasi-normality in direction $u$} holds at $(\bar x,\bar y)$
			w.r.t.\ $\mathcal E$ if there does not exist a nonzero vector 
			$\lambda\in\ker D^*\Phi((\bar x,\bar y);(u,0))$
			satisfying the following condition:
			there are sequences $\{(x_k,y_k)\}_{k\in\N}\subset\gph\Phi$ 
			with $x_k\neq\bar x$ for all $k\in\N$ and
			$\{\lambda_k\}_{k\in\N}\subset\mathbb Y$, $\{\eta_k\}_{k\in\N}\subset\mathbb X$,
			such that we have the convergences from
			\eqref{eq:convergences_definition_quasi_normality}
			and, for all $k\in\N$ and $i\in\{1,\ldots,m\}$, 
			$\eta_k\in \widehat{D}^*\Phi(x_k,y_k)(\lambda_k)$ 
			as well as $\dual{\lambda}{e_i}\dual{y_k-\bar y}{e_i}>0$
			if $\dual{\lambda}{e_i}\neq 0$.
	\end{enumerate}
\end{definition}

In the case where the canonical basis is chosen in $\mathbb Y:=\R^m$, the
above concept of quasi-normality is a direct generalization of the original
notion from \cite{BertsekasOzdaglar2002} which was coined for standard
nonlinear problems and neglected directional information.
Let us just mention that a reasonable, basis-independent definition of quasi-normality 
would require that there exists some basis w.r.t.\ which the mapping of interest is quasi-normal,
see also \cref{thm:quasi_normality_yields_asymptotic_regularity}.

Note that the sequence $\{y_k\}_{k\in\N}$ in the definition of directional pseudo- and
quasi-normality needs to satisfy $y_k\neq\bar y$ for all $k\in\N$. In the definition
of directional pseudo-normality, this is clear from $\dual{\lambda}{y_k-\bar y}>0$ for all
$k\in\N$. Furthermore, in the definition of directional quasi-normality, observe that
$\lambda\neq 0$ implies the existence of $j\in\{1,\ldots,m\}$ such that 
$\dual{\lambda}{e_j}\neq 0$ holds, so that $\dual{y_k-\bar y}{e_j}\neq 0$ is necessary
for each $k\in\N$.

In the following lemma, we show the precise relation between directional pseudo- and quasi-normality.

\begin{lemma}\label{lem:pseudo_vs_quasi_normality}
	Fix $(\bar x,\bar y)\in\gph\Phi$ and some direction $u\in\mathbb S_{\mathbb X}$.
	Then $\Phi$ is pseudo-normal at $(\bar x,\bar y)$ in direction $u$ if and only
	if $\Phi$ is quasi-normal at $(\bar x,\bar y)$ in direction $u$ w.r.t.\
	each orthonormal basis of $\mathbb Y$.
\end{lemma}
\begin{proof}
	$[\Longrightarrow]$ Let $\Phi$ be pseudo-normal at $(\bar x,\bar y)$ in direction $u$,
		let $\mathcal E:=\{e_1,\ldots,e_m\}\subset\mathbb Y$ be an orthonormal basis of
		$\mathbb Y$, and pick $\lambda\in\ker D^*\Phi((\bar x,\bar y);(u,0))$
		as well as sequences $\{(x_k,y_k)\}_{k\in\N}\subset\gph\Phi$ 
		with $x_k\neq\bar x$ for all $k\in\N$ and
		$\{\lambda_k\}_{k\in\N}\subset\mathbb Y$, $\{\eta_k\}_{k\in\N}\subset\mathbb X$,
		satisfying the convergences \eqref{eq:convergences_definition_quasi_normality}
		and, for all $k\in\N$ and $i\in\{1,\ldots,m\}$, 
		$\eta_k\in \widehat{D}^*\Phi(x_k,y_k)(\lambda_k)$ 
		as well as $\dual{\lambda}{e_i}\dual{y_k-\bar y}{e_i}>0$
		if $\dual{\lambda}{e_i}\neq 0$.
		Observing that we have
		\begin{align*}
			\dual{\lambda}{y_k-\bar y}
			&=
			\dual{\sum\nolimits_{i=1}^m\dual{\lambda}{e_i}e_i}{\sum\nolimits_{j=1}^m\dual{y_k-\bar y}{e_j}e_j}
			\\
			&=
			\sum\nolimits_{i=1}^m\sum\nolimits_{j=1}^m\dual{\lambda}{e_i}\dual{y_k-\bar y}{e_j}\dual{e_i}{e_j}
			\\
			&=
			\sum\nolimits_{i=1}^m\dual{\lambda}{e_i}\dual{y_k-\bar y}{e_i},
		\end{align*}
		validity of pseudo-normality at $(\bar x,\bar y)$ in direction $u$ gives $\lambda=0$,
		i.e., $\Phi$ is quasi-normal at $(\bar x,\bar y)$ in direction $u$ w.r.t.\ $\mathcal E$.\\
	$[\Longleftarrow]$ Assume that $\Phi$ is quasi-normal at $(\bar x,\bar y)$ in direction $u$ w.r.t.\
		each orthonormal basis of $\mathbb Y$.
		Suppose that $\Phi$ is not pseudo-normal at $(\bar x,\bar y)$ in direction $u$.
		Then we find some nonzero $\lambda\in\ker D^*\Phi((\bar x,\bar y);(u,0))$
		as well as sequences $\{(x_k,y_k)\}_{k\in\N}\subset\gph\Phi$ 
		with $x_k\neq\bar x$ for all $k\in\N$ and
		$\{\lambda_k\}_{k\in\N}\subset\mathbb Y$, $\{\eta_k\}_{k\in\N}\subset\mathbb X$,
		satisfying the convergences \eqref{eq:convergences_definition_quasi_normality}
		and	$\eta_k\in \widehat{D}^*\Phi(x_k,y_k)(\lambda_k)$ 
		as well as $\dual{\lambda}{y_k-\bar y}>0$ for all $k\in\N$.
		Noting that $\lambda$ does not vanish, we can construct an orthonormal basis
		$\mathcal E_\lambda:=\{e_1^\lambda,\ldots,e_m^\lambda\}$ of $\mathbb Y$ with
		$e_1^\lambda:=\lambda/\norm{\lambda}$. 
		Note that, for $i\in\{1,\ldots,m\}$, 
		we have $\ninnerprod{\lambda}{e_i^\lambda}\neq 0$ if and only
		if $i=1$ by construction of $\mathcal E_\lambda$. Furthermore, we find
		\begin{align*}
			\ninnerprod{\lambda}{e_1^\lambda}\ninnerprod{y_k-\bar y}{e_1^\lambda}
			=
			\norm{\lambda}\ninnerprod{\lambda/\norm{\lambda}}{y_k-\bar y}
			=
			\ninnerprod{\lambda}{y_k-\bar y}
			>
			0.
		\end{align*}
		This, however, contradicts quasi-normality of $\Phi$ at $(\bar x,\bar y)$ in direction $u$
		w.r.t.\ $\mathcal E_{\lambda}$.
\end{proof}

Let us note that \cite[Example~1]{BertsekasOzdaglar2002} shows in the nondirectional
situation of standard nonlinear programming that pseudo-normality might be more 
restrictive than quasi-normality w.r.t.\ the canonical basis in $\R^m$.
On the other hand, due to \cref{lem:pseudo_vs_quasi_normality}, there must exist another 
basis such that quasi-normality w.r.t.\ this basis fails since pseudo-normality fails.
This depicts that validity of quasi-normality indeed may depend on the chosen basis.
In \cite{BaiYeZhang2019}, the authors define directional quasi-normality 
for geometric constraints in Euclidean spaces in componentwise fashion although this 
is somehow unclear in situations where the image space is different from $\R^m$. 
Exemplary, in the $\tfrac12m(m+1)$-dimensional space $\mathcal S_m$
of all real symmetric $m\times m$-matrices, the canonical basis, 
which seems to be associated with a componentwise calculus, comprises 
$\tfrac12(m-1)m$ matrices with precisely two nonzero entries. 
Our definition of quasi-normality from \cref{def:direction_quasi_pseudo_normality} 
gives some more freedom since the choice of the underlying basis allows to \emph{rotate}
the coordinate system.

Following the arguments in \cite[Section~3.2]{BenkoCervinkaHoheisel2019},
it also might be reasonable to define intermediate conditions bridging
pseudo- and quasi-normality. In the light of this paper, however, the concepts
from \cref{def:direction_quasi_pseudo_normality} are sufficient for
our purposes.

As the following theorem shows, directional quasi- and, thus, pseudo-normality
also serve as sufficient conditions for strong directional asymptotic regularity
and directional metric subregularity which explains our interest in these conditions.
Both statements follow once we clarify that pseudo- and quasi-normality
are in fact specifications of the multiplier sequential information in  \eqref{eq:convergences_directional_asymptotic_regularity},
namely $(y_k-\bar y)/\nnorm{y_k-\bar y} - \lambda_k/\nnorm{\lambda_k} \to 0$.

\begin{theorem}\label{thm:quasi_normality_yields_asymptotic_regularity}
  If $\Phi\colon\mathbb X \tto \mathbb Y$ is quasi-normal in direction
  $u \in \mathbb S_{\mathbb X}$ at $(\bar x, \bar y) \in \gph \Phi$
  w.r.t.\ some orthonormal basis $\mathcal E:=\{e_1,\ldots,e_m\}\subset\mathbb Y$ 
  of $\mathbb Y$, then it is also
  strongly asymptotically regular as well as metrically subregular in direction $u$ at $(\bar x, \bar y)$.
\end{theorem}
\begin{proof}
	Fix arbitrary sequences $\{(x_k,y_k)\}_{k\in\N}\subset\gph\Phi$,
	$\{x_k^*\}_{k\in\N}\subset\mathbb X$, and $\{\lambda_k\}_{k\in\N}\subset\mathbb Y$ 
	as well as $x^*\in\mathbb X$ and $y^*\in\mathbb Y$ satisfying 
	$x_k\notin\Phi^{-1}(\bar y)$, $y_k\neq\bar y$, and 
	$x_k^*\in \widehat{D}^*\Phi(x_k,y_k)(\lambda_k)$
	for each $k\in\N$ as well as the convergences
	\eqref{eq:convergences_directional_asymptotic_regularity}.
	Let us define $w_k:=(y_k-\bar y)/\norm{y_k-\bar y}$ 
	and $\tilde\lambda_k:=\lambda_k/\norm{\lambda_k}$ for each $k\in\N$.
	The requirements from \eqref{eq:convergences_directional_asymptotic_regularity} 
	imply that
	$\{w_k\}_{k\in\N}$ and  $\{\tilde\lambda_k\}_{k\in\N}$ converge, 
	along a subsequence (without relabeling),
	to the same nonvanishing limit which we will call $\lambda\in\mathbb S_{\mathbb Y}$.
	Moreover, given $i\in\{1,\ldots,m\}$ with $\dual{\lambda}{e_i} \neq 0$,
	for sufficiently large $k\in\N$, we get $\dual{w_{k}}{e_i} \neq 0 $ and
	\[
		0 
		< 
		\dual{\lambda}{e_i}\dual{ w_{k}}{e_i} 
		= 
		\dual{\lambda}{e_i}\dual{y_k - \bar y}{e_i} /\norm{y_k - \bar y}.
	\]
	Observing that we have $x_k^*/\norm{\lambda_k}\to 0$ from
	\eqref{eq:convergences_directional_asymptotic_regularity}, we find
	$\lambda\in\ker D^*\Phi((\bar x,\bar y);(u,0))$ by definition of the 
	directional limiting coderivative.
	This contradicts validity of quasi-normality of $\Phi$ at $(\bar x,\bar y)$
	in direction $u$ w.r.t.\ $\mathcal E$.
	Particularly, such sequences $\{(x_k,y_k)\}_{k\in\N}$, $\{x_k^*\}_{k\in\N}$, 
	and $\{\lambda_k\}_{k\in\N}$ cannot exist which means
	that $\Phi$ is strongly asymptotically regular in direction $u$ at $(\bar x,\bar y)$.

	The claim about metric subregularity now follows from \cite[Corollary~1]{Gfrerer2014a},
	since the only difference from quasi-normality is the requirement
	\[
	\dual{\lambda_k/\nnorm{\lambda_k}}{(y_k-\bar y)/\nnorm{y_k-\bar y}} \to 1
	\]
	which is the same as $(y_k-\bar y)/\nnorm{y_k-\bar y} - \lambda_k/\nnorm{\lambda_k} \to 0$
	as mentioned in the comments after \cref{thm:directional_asymptotic_stationarity}.
\end{proof}

Relying on this result, \cite[Theorem~7]{Gfrerer2013}
yields that directional pseudo- and quasi-normality provide constraint qualifications
for \eqref{eq:nonsmooth_problem} which ensure validity of directional M-stationarity
at local minimizers.

We would like to point the reader's attention to the fact that nondirectional versions of
pseudo- and quasi-normality are not comparable with the nondirectional version of
asymptotic regularity. This has been observed in the context of standard nonlinear
programming, see \cite[Sections~4.3, 4.4]{AndreaniMartinezRamosSilva2016}.
The reason is that the standard version of asymptotic regularity makes no use of the
multiplier information \eqref{eq:convergences_gamma=1_multiplier}.

	In \cite[Section~4.2]{BenkoMehlitz2022}, which is a preprint version of this paper,
	our new notions of directional pseudo- and quasi-normality from
	\cref{def:direction_quasi_pseudo_normality} are worked out for so-called optimization problem with
	equilibrium constraints which cover models with variational inequality constraints,
	see e.g.\ \cite{FacchneiPang2003,LuoPangRalph1996,OutrataKocvaraZowe1998},
	or bilevel optimization problems, see e.g.\ \cite{Dempe2002,DempeKalashnikovPerezValdesKalashnykova2015}.

\subsubsection{Directional pseudo- and quasi-normality for geometric constraint systems}

Let us now also justify the terminology by showing that the new notions from \cref{def:direction_quasi_pseudo_normality}
coincide with directional pseudo- and quasi-normality in the case of standard constraint mappings 
as studied in \cite{BenkoCervinkaHoheisel2019}.

We start with a general result relying on calmness of the constraint function.
Note that we consider $\bar y:=0$ for simplicity of notation.
Furthermore, we only focus on the concept of directional quasi-normality in
our subsequently stated analysis. Analogous results hold for
directional pseudo-normality.

\begin{proposition}\label{pro:quasi_normality_for_constraint_maps_calm_case}
  A constraint mapping  $\Phi\colon\mathbb X\tto\mathbb Y$ given by $\Phi(x) := g(x) - D$, 
  $x\in\mathbb X$,  where $g\colon\mathbb X \to \mathbb Y$ is 
  a continuous function which is
  calm in direction $u \in \mathbb S_{\mathbb X} $
  at $\bar x\in\mathbb X$ such that $(\bar x, 0) \in \gph \Phi$ 
  and $D \subset \mathbb Y$ is closed,
  is quasi-normal in direction $u$ at $(\bar x, 0)$
  w.r.t.\ some orthonormal basis $\mathcal E:=\{e_1,\ldots,e_m\}\subset\mathbb Y$ of $\mathbb Y$
  provided there do not exist a direction $v \in \mathbb Y$ and 
  a nonzero vector $\lambda\in \mathcal N_D(g(\bar x);v)$ 
  with $0 \in D^*g(\bar x;(u,v))(\lambda)$ satisfying the following condition: 
  there are sequences $\{x_k\}_{k\in\N}\subset\mathbb X$
  with $x_k\neq\bar x$ for all $k\in\N$, $\{z_k\}_{k\in\N}\subset D$,
  $\{\lambda_k\}_{k\in\N}\subset\mathbb Y$, and $\{\eta_k\}_{k\in\N}\subset\mathbb X$
  satisfying $x_k\to\bar x$, $z_k\to g(\bar x)$, $\lambda_k\to\lambda$, $\eta_k\to 0$,
	\begin{equation}\label{eq:directional_convergences_quasi_normality_2}
		\frac{x_k-\bar x}{\norm{x_k-\bar x}}\to u,
		\qquad
		\frac{z_k-g(\bar x)}{\norm{x_k-\bar x}}\to v,
                \qquad
                \frac{g(x_k)-g(\bar x)}{\norm{x_k-\bar x}}\to v,
	\end{equation}
  and, for all $k\in\N$ and $i\in\{1,\ldots,m\}$, $\eta_k \in \widehat{D}^*g(x_k)(\lambda_k)$,
  $\lambda_k\in\widehat{\mathcal N}_D(z_k)$, 
  as well as $\dual{\lambda}{e_i}\dual{g(x_k)-z_k}{e_i}>0$
  if $\dual{\lambda}{e_i}\neq 0$.
  
  Moreover, if $g$ is even calm near $\bar x$, the two conditions are equivalent.
\end{proposition}
\begin{proof}
    $[\Longleftarrow]$ Choose $\lambda\in\ker D^*\Phi((\bar x,0);(u,0))$ 
    and sequences $\{(x_k,y_k)\}_{k\in\N}\subset\gph\Phi$ with $x_k\neq\bar x$ 
    for all $k\in\N$ and
	$\{\lambda_k\}_{k\in\N}\subset\mathbb Y$, $\{\eta_k\}_{k\in\N}\subset\mathbb X$ satisfying
	\eqref{eq:convergences_definition_quasi_normality} with $\bar y:=0$ and, 
	for all $k\in\N$ and $i\in\{1,\ldots,m\}$,
	$\eta_k\in \widehat{D}^*\Phi(x_k,y_k)(\lambda_k)$ 
	as well as $\dual{\lambda}{e_i}\dual{y_{k}}{e_i}>0$ if $\dual{\lambda}{e_i}\neq 0$.
	Applying 
	\cref{lem:coderivatives_constraint_maps}\,\ref{item:constraint_maps_regular_coderivative}
	yields $\eta_k \in \widehat{D}^*g(x_k)(\lambda_k)$ and
	$\lambda_k\in\widehat{\mathcal N}_D(g(x_k) - y_k)$ for each $k\in\N$.
	The assumed calmness of $g$ at $\bar x$ in direction $u$ yields boundedness
	of the sequence $\{(g(x_k) - g(\bar x))/\norm{x_k - \bar x}\}_{k\in\N}$,
	i.e., along a subsequence (without relabeling) it converges to some $v\in\mathbb Y$.
	Note also that
	$(u,v)\in\mathcal T_{\gph g}(\bar x,g(\bar x))$, i.e., $v\in Dg(\bar x)(u)$,
	and that $\{(x_k,g(x_k))\}_{k\in\N}$ converges to $(\bar x,g(\bar x))$ 
	from direction $(u,v)$.
	Setting $z_k := g(x_k) - y_k$ for each $k\in\N$, 
	we get $z_k \to g(\bar x)$ by continuity of $g$
	as well as $\lambda_k\in\widehat{\mathcal N}_D(z_k)$ 
	and $\dual{\lambda}{e_i}\dual{g(x_k)-z_k}{e_i}>0$ if $\dual{\lambda}{e_i}\neq 0$
	for each $k\in\N$ and $i\in\{1,\ldots,m\}$.
	Moreover, we have
	\begin{align*}
  		\frac{z_k - g(\bar x)}{\norm{x_k-\bar x}}
  		=
  		\frac{g(x_k)- g(\bar x)}{\norm{x_k-\bar x}}
  		-
  		\frac{y_k}{\norm{x_k-\bar x}}
  		\to 
  		v - 0 
  		= 
  		v
  	\end{align*}
  	and $v\in \mathcal T_D(g(\bar x))$ follows as well.
    Finally, taking the limit yields
    $\lambda\in\mathcal N_D(g(\bar x);v)$ and $0 \in D^*g(\bar x;(u,v))(\lambda)$,
    so that the assumptions of the proposition imply $\lambda=0$.
    Consequently, $\Phi$ is quasi-normal in direction $u$ at $(\bar x,0)$ w.r.t.\ $\mathcal E$.
    \\
	$[\Longrightarrow]$ Assume that quasi-normality in direction $u$ holds at 
	$(\bar x,0)$ w.r.t.\ $\mathcal E$ and that $g$ is calm around $\bar x$.
	Suppose that there are some $v \in \mathbb Y$,
	$\lambda\in\mathcal N_D(g(\bar x);v)$ with $0 \in D^*g(\bar x;(u,v))(\lambda)$,
	and sequences $\{x_k\}_{k\in\N}\subset\mathbb X$ with $x_k\neq\bar x$ 
	for all $k\in\N$ and $\{z_k\}_{k\in\N}\subset D$, 
	$\{\lambda_k\}_{k\in\N}\subset\mathbb Y$, $\{\eta_k\}_{k\in\N}\subset\mathbb X$ with
	$x_k\to\bar x$, $z_k\to g(\bar x)$, $\lambda_k\to\lambda$, $\eta_k \to 0$,
	\eqref{eq:directional_convergences_quasi_normality_2}, and, for 
	all $k\in\N$ and $i\in\{1,\ldots,m\}$, $\eta_k \in \widehat{D}^*g(x_k)(\lambda_k)$,
	$\lambda_k\in\widehat{\mathcal N}_D(z_k)$, as well as 
	$\dual{\lambda}{e_i}\dual{g(x_k)-z_k}{e_i}>0$
  	as soon as $\dual{\lambda}{e_i}\neq 0$.
  	Set $y_k:=g(x_k)-z_k$ for each $k\in\N$. Then we have $y_k\to 0$,
  	\begin{align*}
  		\frac{y_k}{\norm{x_k-\bar x}}
  		=
  		\frac{g(x_k)-z_k}{\norm{x_k-\bar x}}
  		=
  		\frac{g(x_k)-g(\bar x)}{\norm{x_k-\bar x}}
  		-\frac{z_k-g(\bar x)}{\norm{x_k-\bar x}}
  		\to v - v = 0,
  	\end{align*}
  	and, for all $k\in\N$ and $i\in\{1,\ldots,m\}$, 
  	$\lambda_k\in\widehat{\mathcal N}_D(g(x_k)-y_k)$ 
  	as well as $\dual{\lambda}{e_i}\dual{y_{k}}{e_i}>0$ if $\dual{\lambda}{e_i}\neq 0$.
  	Since $\eta_k \in \widehat{D}^*g(x_k)(\lambda_k)$,
  	calmness of $g$ at $x_k$ implies $\eta_k\in \widehat{D}^*\Phi(x_k,y_k)(\lambda_k)$
  	due to
  	\cref{lem:coderivatives_constraint_maps}\,\ref{item:constraint_maps_regular_coderivative},
  	and taking the limit yields $\lambda\in\ker D^*\Phi((\bar x,\bar y);(u,0))$.
  	Thus, the assumed quasi-normality of $\Phi$ at $(\bar x,0)$ 
  	in direction $u$ w.r.t.\ $\mathcal E$ yields $\lambda=0$ and the claim follows.
\end{proof}

If $g$ is continuously differentiable, the situation becomes a bit simpler
and we precisely recover the notion of directional quasi-normality for geometric constraint
systems as discussed in \cite[Definition~3.4]{BenkoCervinkaHoheisel2019}.

\begin{corollary}\label{cor:quasi_normality_for_constraint_maps_smooth_case}
  A constraint mapping $\Phi\colon\mathbb X\tto\mathbb Y$ given by $\Phi(x) = g(x) - D$, 
  $x\in\mathbb X$, where $g\colon\mathbb X \to \mathbb Y$ is continuously differentiable and 
  $D \subset \mathbb Y$ is closed,
  is quasi-normal in direction $u \in \mathbb S_{\mathbb X}$ at $(\bar x, 0) \in \gph \Phi$ 
  w.r.t.\ some orthonormal basis $\{e_1,\ldots,e_m\}\subset\mathbb Y$ of $\mathbb Y$
  if and only if there does not exist a nonzero vector 
  $\lambda\in\mathcal N_D(g(\bar x);\nabla g(\bar x) u)$ 
  with $ \nabla g(\bar x)^*\lambda=0$
  satisfying the following condition: 
  there are sequences $\{x_k\}_{k\in\N}\subset\mathbb X$ with $x_k\neq\bar x$ 
  for all $k\in\N$, $\{z_k\}_{k\in\N}\subset D$, and $\{\lambda_k\}_{k\in\N}\subset\mathbb Y$
  satisfying $x_k\to\bar x$, $z_k\to g(\bar x)$, $\lambda_k\to\lambda$, 
	\begin{equation}\label{eq:directional_convergences_quasi_normality_3}
		\frac{x_k-\bar x}{\norm{x_k-\bar x}}\to u,
		\qquad
		\frac{z_k-g(\bar x)}{\norm{x_k-\bar x}}\to \nabla g(\bar x)u,
	\end{equation}
  and, for all $k\in\N$ and $i\in\{1,\ldots,m\}$, 
  $\lambda_k\in\widehat{\mathcal N}_D(z_k)$ as well as 
  $\dual{\lambda}{e_i}\dual{g(x_k)-z_k}{e_i}>0$
  if $\dual{\lambda}{e_i}\neq 0$.
\end{corollary}

In \cite[Section~3.3]{BenkoCervinkaHoheisel2019}, it has been reported that under additional
conditions on the set $D$, we can drop the sequences $\{z_k\}_{k\in\N}$ and
$\{\lambda_k\}_{k\in\N}$ from the characterization of directional
quasi-normality in \cref{cor:quasi_normality_for_constraint_maps_smooth_case}.
Particularly, this can be done for so-called \emph{ortho-disjunctive}
programs which cover, e.g., standard nonlinear, complementarity-, cardinality-, or
switching-constrained optimization problems. 
In this regard, \cref{cor:quasi_normality_for_constraint_maps_smooth_case}
reveals that some results from 
\cite{BertsekasOzdaglar2002,Hestenes1975,KanzowRaharjaSchwartz2021b,KanzowSchwartz2010,LiangYe2021}
are covered by our general concept from \cref{def:direction_quasi_pseudo_normality}.

Let us briefly compare our results with the approach from \cite{BaiYeZhang2019}.

\begin{remark}
	Let us consider the setting discussed in \cref{pro:quasi_normality_for_constraint_maps_calm_case}.
    The directional versions of pseudo- and quasi-normality from \cite{BaiYeZhang2019} operate with all nonzero pairs
    of directions $(u,v)$, rather than just a fixed $u$.
    The advantage is that calmness of $g$ plays no role.
    The reason is, however, that the authors in \cite{BaiYeZhang2019} only derive statements
    regarding metric subregularity, but not metric
    subregularity in some fixed direction.
    Calmness of $g$ is needed precisely for preservation
    of directional information.
    We believe that it is useful to know how to verify if a mapping
    is metrically subregular in a specific direction since only
    some directions play a role in many situations.
    We could drop the calmness assumption from \cref{pro:quasi_normality_for_constraint_maps_calm_case}, but, 
    similarly as in \cite[Theorem~3.1]{BenkoGfrererOutrata2019}, additional directions of the type $(0,v)$ for a nonzero $v$ would appear.
    Clearly, such directions are included among all nonzero pairs $(u,v)$, but the connection to the original direction $u$ 
    would have been lost.
\end{remark}

\subsection{Sufficient conditions for asymptotic regularity via pseudo-coderivatives}\label{sec:asymptotic_regularity_via_super_coderivative}

\subsubsection{The role of super-coderivatives}

We start this section by interrelating the concept of super-coderivatives from \cref{def:super_coderivative} 
and asymptotic regularity.
Fix $(\bar x,\bar y)\in\gph\Phi$ and choose
$\{(x_k,y_k)\}_{k\in\N}\subset\gph\Phi$, 
$\{x_k^*\}_{k\in\N}\subset\mathbb X$, 	and $\{\lambda_k\}_{k\in\N}\subset\mathbb Y$ 
as well as $x^*\in\mathbb X$ and $y^*\in\mathbb Y$ satisfying
$x_k\notin\Phi^{-1}(\bar y)$, $y_k\neq\bar y$, 
and $x_k^*\in\widehat{D}^*\Phi(x_k,y_k)(\lambda_k)$
for all $k\in\N$ as well as the convergences
\eqref{eq:convergences_directional_asymptotic_regularity}.
For each $k\in\N$, we set $t_k:=\norm{x_k-\bar x}$, $\tau_k:=\norm{y_k-\bar y}$,
\[
	u_k:=\frac{x_k-\bar x}{\norm{x_k-\bar x}},\qquad 
	v_k:=\frac{y_k-\bar y}{\norm{y_k-\bar y}},\qquad
	y_k^*:=\frac{\norm{y_k-\bar y}}{\norm{x_k-\bar x}}\lambda_k,
\]
and find $\tau_k/t_k\to 0$ as well as
\[
	\forall k\in\N\colon\quad
		x_k^*\in \widehat{D}^*\Phi(\bar x+t_ku_k,\bar y+\tau_kv_k)((t_k/\tau_k)y_k^*).
\]
Along a subsequence (without relabeling), $v_k\to v$ holds for some 
$v\in\mathbb S_{\mathbb Y}$.
Thus, taking the limit $k\to\infty$, 
we have $x^*\in D^*_\textup{sup}\Phi((\bar x,\bar y);(u,v))(y^*)$
by definition of the super-coderivative.
Moreover, from \eqref{eq:convergences_directional_asymptotic_regularity} 
we also know that $y^* = \norm{y^*} v$.
Consequently, we come up with the following lemma.

\begin{lemma}\label{lem:super_coderivative_vs_asymptotic_regularity}
	Let $(\bar x,\bar y)\in\gph\Phi$ and $u\in\mathbb S_{\mathbb X}$ be fixed.
	If 
	\[
		\bigcup_{v\in\mathbb S_{\mathbb Y}}
		D^*_\textup{sup}\Phi((\bar x,\bar y);(u,v))(\beta v)
		\subset
		\Im D^*\Phi(\bar x,\bar y) 
	\]
	holds for all $\beta \geq 0$, then $\Phi$ is asymptotically regular at $(\bar x,\bar y)$ in direction $u$.
	If the above estimate holds for all $\beta\geq 0$
	with $\Im D^*\Phi(\bar x,\bar y)$ replaced by $\Im D^*\Phi((\bar x,\bar y);(u,0))$,
	then $\Phi$ is strongly asymptotically regular at $(\bar x,\bar y)$ in direction $u$.
\end{lemma}

The next result, 
which is based on hypothesis $A^\gamma(u)$, see \cref{ass:A_gamma_u},
follows as a corollary of 
\cref{lem:super_coderivative_vs_asymptotic_regularity,lem:super_coderivative_vs_pseudo_coderivative},
and gives new sufficient conditions for directional asymptotic regularity.
Note that strong directional asymptotic regularity can be handled analogously
	by employing an adjusted version of $A^\gamma(u)$ where
	$\Im D^*\Phi(\bar x,\bar y)$ in the right-hand side of
	\eqref{eq:sufficient_condition_asym_reg_pseudo} is replaced
	by $\Im D^*\Phi((\bar x,\bar y);(u,0))$.
\begin{theorem}\label{thm:asymptotic_regularity_via_pseudo_coderivatives}
	Let $(\bar x,\bar y)\in\gph\Phi$, $u\in\mathbb S_{\mathbb X}$, and $\gamma>1$ be fixed.
	If $A^\gamma(u)$ holds,
	then $\Phi$ is asymptotically regular at $(\bar x,\bar y)$ in direction $u$.
\end{theorem}

In the case where the pseudo-coderivatives involved in 
the construction of $A^\gamma(u)$
can be computed or estimated from above, new applicable sufficient conditions for (strong)
directional asymptotic regularity are 
provided by \cref{thm:asymptotic_regularity_via_pseudo_coderivatives}.
Particularly, in situations where $\Phi$ is given in form of a constraint
mapping and $\gamma:=2$ is fixed, we can rely on the results obtained
in \cref{sec:variational_analysis_constraint_mapping} in order to make the findings
of \cref{thm:asymptotic_regularity_via_pseudo_coderivatives} more specific.
This will be done in the next subsection.

\subsubsection{The case of constraint mappings}\label{Sec:5_constr_mappings}

Throughout the section, we assume that $\Phi\colon\mathbb X\tto\mathbb Y$
is given by $\Phi(x):=g(x)-D$, $x\in\mathbb X$, where $g\colon\mathbb X\to\mathbb Y$ is a
twice continuously differentiable function and $D\subset\mathbb Y$
is a closed set. Furthermore, for simplicity of notation, we fix $\bar y:=0$ which
is not restrictive as already mentioned earlier.

We start with a general result which does not rely on any additional structure
of the set $D$.
\begin{theorem}\label{thm:asymptotic_regularity_pseudo_coderivative_non_polyhedral}
	Let $(\bar x,0)\in\gph\Phi$ as well as $u\in\mathbb S_{\mathbb X}$
	be fixed. Assume that \eqref{eq:CQ_pseudo_subregularity_II} holds,
  		as well as \eqref{eq:CQ_pseudo_subregularity_Ia} or, in the case $\nabla g(\bar x)u\neq 0$,
  		\eqref{eq:CQ_pseudo_subregularity_Ib}.
	If, for each $x^*\in\mathbb X$ and $y^*,z^*\in\mathbb Y$ satisfying
			\begin{subequations}\label{eq:system_pseudo_coderivative_nonpolyhedral}
				\begin{align}
					\label{eq:pseudo_coderivative_nonpolyhedral_x}
						x^*
						&=
						\nabla^2\dual{y^*}{g}(\bar x)(u)
						+
						\nabla g(\bar x)^*z^*,
						\\
					\label{eq:pseudo_coderivative_nonpolyhedral_y}
						y^*
						&\in 
						\mathcal N_D(g(\bar x);\nabla g(\bar x)u)
						\cap
						\ker\nabla g(\bar x)^*,
						\\
					\label{eq:pseudo_coderivative_nonpolyhedral_z}
						z^*
						&\in 
						D\mathcal N_D(g(\bar x),y^*)(\nabla g(\bar x)u),
				\end{align}
			\end{subequations}
			there is some $\lambda\in\mathcal N_D(g(\bar x))$
			such that $x^*=\nabla g(\bar x)^*\lambda$,
			then $\Phi$ is asymptotically regular at $(\bar x,0)$ in direction $u$.
			Moreover, $\Phi$ is even	strongly asymptotically regular at $(\bar x,0)$ in direction $u$
			if $\lambda$ can be chosen from $\mathcal N_D(g(\bar x);\nabla g(\bar x)u)$.
\end{theorem}
\begin{proof}
	\cref{The : NCgen}\,\ref{item:general_estimate_+CQ} implies
	$\ker\widetilde D^*_2\Phi((\bar x,0);(u,0))=\{0\}$
	(and so, due to \eqref{eq:trivial_upper_estimate_pseudo_coderivative}, 
	also $\ker D^*_2\Phi((\bar x,0);(u,0))=\{0\}$)
	as well as that for $x^*\in\Im \widetilde D^*_2\Phi((\bar x,0);(u,0))$,
	we find $y^*,z^*\in\mathbb Y$ satisfying 
	\eqref{eq:system_pseudo_coderivative_nonpolyhedral}.
	The assumptions guarantee that we can find $\lambda\in\mathcal N_D(g(\bar x))$
	such that $x^*=\nabla g(\bar x)^*\lambda\in \Im D^*\Phi(\bar x,0)$
	where we used 
	\cref{lem:coderivatives_constraint_maps}\,\ref{item:constraint_maps_limiting_coderivative}.
	It follows $\Im \widetilde D^*_2\Phi((\bar x,0);(u,0))\subset\Im D^*\Phi(\bar x,0)$.
	Thus, \cref{thm:asymptotic_regularity_via_pseudo_coderivatives}
	shows that $\Phi$ is asymptotically regular at $(\bar x,0)$ in direction $u$.
	The statement regarding strong asymptotic regularity follows in analogous way while
	respecting 
	\cref{lem:coderivatives_constraint_maps}\,\ref{item:constraint_maps_directional_limiting_coderivative}.
\end{proof}

We note that \eqref{eq:CQ_pseudo_subregularity_Ia} is stronger than \eqref{eq:CQ_pseudo_subregularity_Ib}
when $\nabla g(\bar x)u\neq 0$ holds, see \eqref{eq:trivial_upper_estimate_graphical_subderivative}.
Naturally, this means that it is sufficient to check \eqref{eq:CQ_pseudo_subregularity_Ia}
regardless whether $\nabla g(\bar x)u$ vanishes or not. In the case $\nabla g(\bar x)u\neq 0$, however,
it is already sufficient to check the milder condition \eqref{eq:CQ_pseudo_subregularity_Ib}.
This will be important later on, see \cref{Pro:Milder_than_FOSCMS} and \cref{rem:forgetting_directional_information_in_CQ} below.

Note also that we implicitly relied on condition \eqref{eq:sufficient_condition_asym_reg_pseudo_rough} 
(with $\bar y:=0$ and $\gamma:=2$)
in the proof of \cref{thm:asymptotic_regularity_pseudo_coderivative_non_polyhedral}, 
and not on the milder refined condition \eqref{eq:sufficient_condition_asym_reg_pseudo} 
(again with $\bar y:=0$ and $\gamma:=2$) which appears in the
statement of $A^\gamma(u)$. This happened due to the generality of the setting 
in \cref{thm:asymptotic_regularity_pseudo_coderivative_non_polyhedral}.
In the polyhedral situation, \eqref{eq:sufficient_condition_asym_reg_pseudo} can be employed to obtain the
following improved result.

\begin{theorem}\label{thm:asymptotic_regularity_pseudo_coderivative_polyhedral_II}
	Let $(\bar x,0)\in\gph\Phi$ as well as $u\in\mathbb S_{\mathbb X}$
	be fixed. 
	Let $\mathbb Y:=\R^m$ and let $D$ be polyhedral
	locally around $g(\bar x)$.
	Assume that condition \eqref{eq:CQ_pseudo_subregularity_polyhedral_II} holds for each $s\in\mathbb X$.
	If, for each 
	$x^*,s\in\mathbb X$, $y^*,z^*\in\R^m$, and $\alpha\geq 0$ 
	satisfying \eqref{eq:pseudo_coderivative_nonpolyhedral_x} and
			\begin{equation}\label{eq:system_pseudo_coderivative_polyhedral_II}
				\begin{aligned}
						y^*
						&\in 
						\mathcal N_{\mathbf T(u)}(w_s(u,v))
						\cap
						\ker\nabla g(\bar x)^*,
						\\
						z^*
						&\in 
							\mathcal N_{\mathbf T(u)}(w_s(u,v))
						\quad
						\big(
						\textrm{or } \ z^*
						\, \in 
						\mathcal T_{
							\mathcal N_{\mathbf T(u)}(w_s(u,v))
							}
							(y^*)
						\big),
				\end{aligned}
			\end{equation}
		where $v:=\alpha y^*$, 
		and $\mathbf T(u)$ as well as $w_s(u,v)$ have been defined in \eqref{eq:Tu_and_ws},
		there is some $\lambda\in\mathcal N_D(g(\bar x))$
		such that $x^*=\nabla g(\bar x)^*\lambda$,
		then $\Phi$ is asymptotically regular at $(\bar x,0)$ in direction $u$.
		Moreover, $\Phi$ is even strongly asymptotically regular at $(\bar x,0)$ in direction $u$
		if $\lambda$ can be chosen from $\mathcal N_{\mathcal T_D(g(\bar x))}(\nabla g(\bar x)u)$.
\end{theorem}
\begin{proof}
	Due to \cref{The : NCgen_2}, \eqref{eq:CQ_pseudo_subregularity_polyhedral_II} yields 
	$\ker D^*_2\Phi((\bar x,0);(u,0))\subset\{0\}$ in the present situation.
	Now, fix $x^*\in \widetilde D^*_2\Phi((\bar x,0);(u,0))(0)$.
	Then \cref{The : NCgen}\,\ref{item:polyhedral_estimate} shows the existence
	of $z^*\in\mathcal N_{\mathcal T_D(g(\bar x))}(\nabla g(\bar x)u)$ such that
	$x^*=\nabla g(\bar x)^*z^*$.
		Let us now consider the case $x^*\in D^*_2\Phi((\bar x,0);(u,\bar \alpha w))(\bar \beta w)$ for 
		some $w\in\mathbb S_{\R^m}$ and $\bar \alpha,\bar \beta \geq 0$.
		If $\bar\beta=0$ holds, we can employ \eqref{eq:trivial_upper_estimate_pseudo_coderivative}
		to find $x^*\in \widetilde D^*_2\Phi((\bar x,0);(u,0))(0)$ and, thus,
		the above argumentation applies.
		Thus, let us consider $\bar\beta>0$ and set $\alpha:=\bar\alpha/\bar\beta$.
		Then we have $x^*\in D^*_2\Phi((\bar x,0);(u,v))(y^*)$ for
		$v=\alpha y^*$ with $v:=\bar\alpha w$ and $y^*:=\bar\beta w$.
	\cref{The : NCgen_2} implies the existence of $s\in\mathbb X$ such that
	\eqref{eq:pseudo_coderivative_nonpolyhedral_x} and
	\eqref{eq:system_pseudo_coderivative_polyhedral_II} hold with $v=\alpha y^*$.
	Now, the postulated assumptions guarantee the existence of
	$\lambda\in\mathcal N_{D}(g(\bar x))$
	such that $x^*=\nabla g(\bar x)^*\lambda$.
	Respecting 
	\cref{lem:coderivatives_constraint_maps}\,\ref{item:constraint_maps_limiting_coderivative},
	this shows
	\eqref{eq:sufficient_condition_asym_reg_pseudo} with $\bar y:=0$ and $\gamma:=2$.
	Thus, \cref{thm:asymptotic_regularity_via_pseudo_coderivatives} yields
	that $\Phi$ is asymptotically regular at $(\bar x,0)$ in direction $u$.
	The statement regarding strong asymptotic regularity follows analogously.
\end{proof}

	Due to \cref{cor:M_stationarity_via_directional_asymptotic_regularity}, 
	\cref{thm:asymptotic_regularity_pseudo_coderivative_non_polyhedral,thm:asymptotic_regularity_pseudo_coderivative_polyhedral_II}
	provide constraint qualifications for M-stationarity.
	Interestingly, one can easily check that the same conditions can also be
	obtained from \cref{Pro:M-stat_via_second_order} by demanding that
	any mixed-order stationary point is already M-stationary.	

In the remaining part of the section, we prove that the assumptions of
\cref{thm:asymptotic_regularity_pseudo_coderivative_non_polyhedral}
are not stronger than 
FOSCMS$(u)$ while the assumptions of
\cref{thm:asymptotic_regularity_pseudo_coderivative_polyhedral_II}
are strictly weaker than the so-called \emph{Second-Order Sufficient Condition for Metric Subregularity} (SOSCMS) 
in direction $u$.

Given a point $\bar x\in\mathbb X$ with $(\bar x,0)\in\gph\Phi$,
\cref{lem:coderivatives_constraint_maps}\,\ref{item:constraint_maps_directional_limiting_coderivative}
shows that the condition
\[
	u\in\mathbb S_{\mathbb X},\,
	\nabla g(\bar x)u\in\mathcal T_D(g(\bar x)),\,
	\nabla g(\bar x)^*y^*=0,\,
	y^*\in\mathcal N_D(g(\bar x);\nabla g(\bar x)u)
	\quad
	\Longrightarrow
	\quad
	y^*=0
\]
equals FOSCMS in the current setting.
In the case where $D$ is locally polyhedral around $g(\bar x)$, the refined condition
\[
	\left.
	\begin{aligned}
	&u\in\mathbb S_{\mathbb X},\,
	\nabla g(\bar x)u\in\mathcal T_D(g(\bar x)),\,
	\nabla g(\bar x)^*y^*=0,\\
	&
	\nabla^2\innerprod{y^*}{g}(\bar x)[u,u]\geq 0,\,
	y^*\in\mathcal N_D(g(\bar x);\nabla g(\bar x)u)
	\end{aligned}
	\right\}
	\quad
	\Longrightarrow
	\quad
	y^*=0,
\]
is referred to as SOSCMS in the literature.
As these names suggest, both conditions are sufficient for metric subregularity of
$\Phi$ at $(\bar x,0)$, see \cite[Corollary~1]{GfrererKlatte2016}.
Particularly, they provide constraint qualifications for M-stationarity of local minimizers.
Again, with the aid of \cref{lem:coderivatives_constraint_maps}\,\ref{item:constraint_maps_directional_limiting_coderivative}, 
one can easily check that
\[
	\nabla g(\bar x)^*y^*=0,\,
	y^*\in\mathcal N_D(g(\bar x);\nabla g(\bar x)u)
	\quad
	\Longrightarrow
	\quad
	y^*=0
\]
equals FOSCMS$(u)$ in the present setting, and
\[
	\nabla g(\bar x)^*y^*=0,\,
	\nabla^2\innerprod{y^*}{g}(\bar x)[u,u]\geq 0,\,
	y^*\in\mathcal N_D(g(\bar x);\nabla g(\bar x)u)
	\quad
	\Longrightarrow
	\quad
	y^*=0
\]
will be denoted by SOSCMS$(u)$.
Each of the conditions FOSCMS$(u)$ and SOSCMS$(u)$
is sufficient for metric subregularity of $\Phi$ at $(\bar x,0)$ in direction $u$.

\begin{proposition}\label{Pro:Milder_than_FOSCMS}
	Consider $(\xb,0) \in \gph \Phi$ and $u\in\mathbb S_{\mathbb X}$.
	Under FOSCMS$(u)$ all assumptions of \cref{thm:asymptotic_regularity_pseudo_coderivative_non_polyhedral} are satisfied.
\end{proposition}
\begin{proof}
Let $y^* \in \mathcal N_{D}(g(\xb);\nabla g(\xb) u)$
be such that $\nabla g(\xb)^* y^* = 0$.
Then FOSCMS$(u)$ yields $y^* = 0$ 
and so \eqref{eq:CQ_pseudo_subregularity_II} is satisfied.
Moreover, we only need to show the remaining assertions for $y^* = 0$.

Assume that $\nabla g(\xb)u \neq 0$ holds.
Suppose now that \eqref{eq:CQ_pseudo_subregularity_Ib} is violated, i.e., 
there exists
$\hat z^* \in D_{\textrm{sub}}\mathcal N_{D}(g(\xb),0)(q)$ for $q:=\nabla g(\xb)u/\norm{\nabla g(\xb)u}$ with $\nabla g(\bar x)^*\hat z^*=0$. 
By \cref{lem:technical_property_D_normal_cone_map} and FOSCMS$(u)$,
we thus get $\hat z^* = 0$
which is a contradition since $\hat z^* \in  \mathbb S_{\mathbb Y}$ by \cref{def:graphical_derivative}.
Similarly, in the case $\nabla g(\bar x)u=0$, we can verify \eqref{eq:CQ_pseudo_subregularity_Ia} which reduces to
\[
	\nabla g(\bar x)^*\hat z^*=0,\quad
	\hat{z}^* \in D\mathcal N_{D}(g(\xb),0)(0)
	\quad
	\Longrightarrow
	\quad
	\hat z^*=0.
\] 
Applying \cref{lem:technical_property_D_normal_cone_map} again,
we get $\hat{z}^* \in \mathcal N_D(g(\xb))$ which implies $\hat{z}^* = 0$
since FOSCMS$(u)$ corresponds to the Mordukhovich criterion 
due to $\nabla g(\xb)u = 0$.
Thus, we have shown that \eqref{eq:CQ_pseudo_subregularity_Ia} or,
in the case $\nabla g(\bar x)u\neq 0$, \eqref{eq:CQ_pseudo_subregularity_Ib}
holds.

Validity of the last assumption follows immediately
since $z^* \in \mathcal N_D(g(\xb);\nabla g(\xb)u)$
is obtained from \cref{lem:technical_property_D_normal_cone_map}, 
and so we can just take $\lambda := z^*$ due to $y^*=0$.
\end{proof}

\begin{remark}\label{rem:forgetting_directional_information_in_CQ}
	Note that for $u\in\mathbb S_{\mathbb X}$ satisfying $\nabla g(\bar x)u\neq 0$, we have the trivial
	upper estimate
	$D_{\textup{sub}}\mathcal N_{D}(g(\xb),y^*)(\nabla g(\xb)u/\norm{\nabla g(\bar x)u}) 
	\subset D\mathcal N_{D}(g(\xb),y^*)(0)$.
	Hence, in \cref{thm:asymptotic_regularity_pseudo_coderivative_non_polyhedral}, it is possible to
	replace validity of \eqref{eq:CQ_pseudo_subregularity_Ia} or, in the case $\nabla g(\bar x)u\neq 0$,
  	\eqref{eq:CQ_pseudo_subregularity_Ib} by the slightly stronger assumption that
  	\eqref{eq:CQ_pseudo_subregularity_Ia} has to hold (even in the case $\nabla g(\bar x)u\neq 0$).
  	However, we cannot show anymore that FOSCMS$(u)$ is sufficient for this stronger assumption to hold,
  	i.e., dropping directional information comes for a price.
\end{remark}

\begin{proposition}\label{Pro:Milder_than_SOSCMS}
	Let $(\bar x,0)\in\gph\Phi$ as well as $u\in\mathbb S_{\mathbb X}$
	be fixed, let $\mathbb Y:=\R^m$, and let $D$ be polyhedral
	locally around $g(\bar x)$.
	If SOSCMS$(u)$ is valid, then the assumptions of
	\cref{thm:asymptotic_regularity_pseudo_coderivative_polyhedral_II}
	are satisfied.
\end{proposition}
\begin{proof}
The key step is to realize that if $y^* \in \mathcal N_{\mathbf T(u)}(w_s(u,v))\cap\ker\nabla g(\bar x)^*$
for some $s\in\mathbb X$ and $v \in \R^m$,
then we get
\[
	\frac 12 \nabla^2\langle y^*,g\rangle(\bar x)[u,u]
	= 
	\innerprod{w_s(u,v)}{ y^*}+ \innerprod{v}{y^*}
	=
	 \innerprod{v}{ y^* }
\]
by \cref{rem:refined_polyhedral_situation_normal_cone_relation} and $\nabla g(\xb)^* y^* = 0$,
and $y^*\in\mathcal N_D(g(\bar x);\nabla g(\bar x)u)$ also holds,
again by \cref{rem:refined_polyhedral_situation_normal_cone_relation}.

Then \eqref{eq:CQ_pseudo_subregularity_polyhedral_II}
follows because for $y^* \in \mathcal N_{\mathbf T(u)}(w_s(u,0))\cap\ker\nabla g(\bar x)^*$,
$\nabla^2\langle y^*,g\rangle(\bar x)[u,u] = 0$ is obtained,
and SOSCMS$(u)$ yields $y^* = 0$.

Next, for arbitrary $y^* \in \mathcal N_{\mathbf T(u)}(w_s(u,v))\cap\ker\nabla g(\bar x)^*$ with $s\in\mathbb X$ and 
$v:=\alpha y^*$ for some $\alpha\geq 0$,
we get $\nabla^2\langle y^*,g\rangle(\bar x)[u,u]
= 2\innerprod{v}{ y^* }=2\alpha\nnorm{y^*}^2 \geq 0$,
so SOSCMS$(u)$ can still be applied to give $y^* = 0$.
Now, we can always take $\lambda := z^*$ since $z^* \in \mathcal N_{\mathbf T(u)}(w_s(u,v)) \subset \mathcal N_D(g(\bar x);\nabla g(\bar x)u)$.
\end{proof}

We immediately arrive at the following corollary.

\begin{corollary}\label{cor:FOSCMS_SOSCMS_imply_sAR}
	The constraint mapping $\Phi$ is strongly asymptotically regular at $(\xb,0) \in \gph \Phi$ in direction $u\in\mathbb S_{\mathbb X}$ if
	FOSCMS$(u)$ holds or if $\mathbb Y:=\R^m$, $D$ is locally polyhedral around $g(\xb)$, and SOSCMS$(u)$ holds.
\end{corollary}

The following example shows that our new conditions from \cref{thm:asymptotic_regularity_pseudo_coderivative_polyhedral_II} 
are in fact strictly milder than SOSCMS.

\begin{example}\label{ex:milder_than_SOSCMS}
	Let $g\colon\R \to \R^2$ and $D \subset \R^2$ be given by
	$g(x) := (x,-x^2)$, $x\in\R$, and $D:=(\R_+ \times \R) \cup (\R \times \R_+)$.
	Observe that $D$ is a polyhedral set.
	We consider the constraint map $\Phi\colon\R\tto\R^2$ given by $\Phi(x):=g(x)-D$, $x\in\R$.
	We note that $\Phi^{-1}(0)=[0,\infty)$ holds.
	Hence, fixing $\bar x:=0$, we can easily check that $\Phi$ is metrically
	subregular at $(\bar x,0)$ in direction $1$ but not in direction $-1$,
	i.e., FOSCMS and SOSCMS must be violated.
	
	First, we claim that all the assumptions from \cref{thm:asymptotic_regularity_pseudo_coderivative_polyhedral_II} are satisfied for $u=\pm 1$.
	Taking into account \cref{rem:refined_polyhedral_situation_normal_cone_relation}, it suffices to verify these assumptions
	for $w_s(u,v)$ replaced by $0$.
	Let us fix $u=\pm 1$, $y^*, z^* \in\mathcal N_{\mathcal T_D(g(\xb))}(\nabla g(\bar x)u)$
	such that $\nabla g(\bar x)^* y^* = 0$ and 
	$\nabla^2\langle y^*,g\rangle(\bar x)(u)+\nabla g(\bar x)^* z^*=x^*$ for $x^* \in \R$.
	We have $\nabla g(\bar x)u=(u,0)$, 
	$\nabla^2\langle y^*,g\rangle(\bar x)(u)=-2y_2^*u$, and
	\[
		\mathcal N_{\mathcal T_D(g(\xb))}(\nabla g(\bar x)u) 
		=
		\begin{cases} 
		\{0\} \times \R_-	& u=-1,\\
		\{(0,0)\}			& u=1.
		\end{cases}
	\]
	Thus, for $u=1$, we have $y^*=0$ regardless of $x^*$.
	Hence, condition \eqref{eq:CQ_pseudo_subregularity_polyhedral_II} holds trivially
	and we can choose $\lambda:=z^*$ to find $x^* = \nabla g(\bar x)^*\lambda$ as well as $\lambda\in\mathcal N_{\mathcal T_D(g(\xb))}(\nabla g(\bar x)u)$. 
	For $u=-1$, we get $y_1^*=z_1^*=0$ and $y^*_2\leq 0$.
	Thus, if $x^*=0$, from $-2 y_2^*u + z_1^* = 0$
	we deduce $y_2^*=0$, and \eqref{eq:CQ_pseudo_subregularity_polyhedral_II} follows.
	For arbitrary $x^*\in\R$, we get $x^* = -2 y_2^*u + z_1^* = 2 y_2^* \leq 0$
	and we can choose $\lambda := (x^*,0) \in \mathcal N_D(g(\bar x))$ to obtain $\nabla g(\bar x)^*\lambda=x^*$.
	Note, however, that $(x^*,0) \notin \mathcal N_{\mathcal T_D(g(\xb))}(\nabla g(\bar x)u) = \{0\} \times \R_-$
	unless $x^* = 0$.
	
	Regarding the assumptions of \cref{thm:asymptotic_regularity_pseudo_coderivative_non_polyhedral},
	let us just mention, without providing the details, that \eqref{eq:CQ_pseudo_subregularity_Ia} and
	\eqref{eq:CQ_pseudo_subregularity_Ib}
	fail since the graphical (sub)derivative is too large.
	Particularly, this clarifies that these assumptions are not necessary e.g.\ in the polyhedral setting,
	but not because they would be satisfied automatically.
\end{example}
	
\section{Concluding remarks}\label{sec:conclusions}

In this paper, we enriched the general concepts of asymptotic stationarity and regularity with the aid
of tools from directional limiting variational analysis. Our central result
\cref{thm:higher_order_directional_asymptotic_stationarity} states that, even in the absence of any
constraint qualification, local minimizers of a rather general optimization problem are M-stationary,
mixed-order stationary in terms of a suitable pseudo-coderivative, or asymptotically stationary
in a critical direction (of a certain order). 
By ruling out the last option, we were in position to
distill new mixed-order necessary optimality conditions. Some novel upper estimates
for the second-order directional pseudo-coderivative of constraint mappings were
successfully employed to make these results fully explicit in the presence of
geometric constraints. Our findings also gave rise to the formulation of directional
notions of asymptotic regularity for set-valued mappings.
These conditions have been shown to serve as constraint qualifications guaranteeing
M-stationarity of local minimizers in nonsmooth optimization. We embedded these new qualification
conditions into the landscape of constraint qualifications which are
already known from the literature, showing that these conditions
are comparatively mild.
Noting that directional asymptotic regularity might be difficult to check in practice,
we then focused on the derivation of applicable sufficient conditions for its validity.
First, we suggested directional notions of pseudo- and quasi-normality for that purpose
which have been shown to generalize related concepts for geometric constraint systems to
arbitrary set-valued mappings. 
Second, with the aid of so-called super- and pseudo-coderivatives, sufficient conditions
for the presence of directional asymptotic regularity for geometric constraint systems
in terms of first- and second-order derivatives of the associated mapping as well as
standard variational objects associated with the underlying set were derived. 
These sufficient conditions turned out to be not stronger than 
the First- and Second-Order Sufficient Condition for Metric Subregularity from the literature.

In this paper, we completely neglected to study the potential value of directional
asymptotic regularity in numerical optimization 
which might be a promising topic of future research.
Furthermore, it has been shown in \cite{Mehlitz2020a} that nondirectional asymptotic
regularity can be applied nicely as a qualification condition in the limiting variational
calculus. Most likely, directional asymptotic regularity may play a similar role 
in the directional limiting calculus. 
Finally, it seems desirable to further develop the calculus for pseudo-coderivatives 
for mappings which possess a more difficult structure than constraint mappings.

\subsection*{Acknowledgments}

The authors would like to thank the referees and the associated editor for valuable comments 
which helped to improve the presentation of the material.
Particularly, the authors are grateful to one of the reviewers who pointed out the close
relationship with 2-regularity and suggested \cref{ex:MPCC}.
The research of Mat\'u\v{s} Benko was supported by the Austrian Science Fund (FWF) under grant P32832-N
as well as by the infrastructure of the Institute of Computational Mathematics, Johannes Kepler University Linz, Austria.


\appendix

\section{Missing proofs}\label{sec:appendix}

\begin{proof}[Proof of \cref{lem:pseudo_coderivative_via_limiting_normals}]
	We only verify the (more technical) assertion 
	regarding \cref{def:coderivatives}\,\ref{item:new_pseudo_coderivative} 
	as the proof for the assertion which addresses
	\cref{def:coderivatives}\,\ref{item:Gfrerer_pseudo_coderivative} 
	follows in similar (but slightly easier) fashion.
	
	Thus, fix $x^*\in\mathbb X$ and $y^*\in\mathbb Y$ as well as 
	$\{u_k\}_{k\in\N},\{x_k^*\}_{k\in\N}\subset\mathbb X$,
	$\{v_k\}_{k\in\N},\{y_k^*\}_{k\in\N}\subset\mathbb Y$,
	and $\{t_k\}_{k\in\N}\subset\R_+$
	which satisfy
	$u_k\to u$, $v_k\to v$, $t_k\downarrow 0$, $x_k^*\to x^*$, $y_k^*\to y^*$, and
	\[
		\forall k\in\N\colon\quad
		\left(x_k^*,-\frac{y_k^*}{(t_k\nnorm{u_k})^{\gamma-1}}\right)
		\in 
		\mathcal N_{\gph\Phi}(\bar x+t_ku_k,\bar y+(t_k\nnorm{u_k})^\gamma v_k).
	\] 
	By definition of the limiting normal cone, for each $k\in\N$, we find
	$\{x_{k,\ell}\}_{\ell\in\N},\{x_{k,\ell}^*\}_{\ell\in\N}\subset\mathbb X$
	and
	$\{y_{k,\ell}\}_{\ell\in\N},\{y_{k,\ell}^*\}_{\ell\in\N}\subset\mathbb Y$
	such that
	$x_{k,\ell}\to\bar x+t_ku_k$, $y_{k,\ell}\to \bar y+(t_k\nnorm{u_k})^\gamma v_k$,
	$x_{k,\ell}^*\to x_k^*$, and $y_{k,\ell}^*\to y_k^*/(t_k\nnorm{u_k})^{\gamma-1}$
	as $\ell\to\infty$ and
	$(x_{k,\ell}^*,-y_{k,\ell}^*)\in\widehat{\mathcal N}_{\gph\Phi}(x_{k,\ell},y_{k,\ell})$
	for all $\ell\in\N$.
	
	For each $k\in\N$, let us define sequences $\{u_{k,\ell}\}_{\ell\in\N}\subset\mathbb X$
	and $\{v_{k,\ell}\}_{\ell\in\N},\{\hat y_{k,\ell}^*\}_{\ell\in\N}\subset\mathbb Y$ by means of
	\[
		\forall \ell\in\N\colon\quad
		u_{k,\ell}:=\frac{x_{k,\ell}-\bar x}{t_k},\qquad
		v_{k,\ell}:=\frac{y_{k,\ell}-\bar y}{(t_k\nnorm{u_{k,\ell}})^\gamma},\qquad
		\hat y_{k,\ell}^*:=(t_k\nnorm{u_{k,\ell}})^{\gamma-1}y_{k,\ell}^*.
	\]
	This gives
	\begin{equation}\label{eq:robustness_pseudo_coderivative_surrogate_sequences}
		\forall\ell\in\N\colon\quad
		\left(x_{k,\ell}^*,-\frac{\hat y_{k,\ell}^*}{(t_k\nnorm{u_{k,\ell}})^{\gamma-1}}\right)
		\in 
		\widehat{\mathcal N}_{\gph\Phi}
			\bigl(\bar x+t_ku_{k,\ell},\bar y+(t_k\nnorm{u_{k,\ell}})^\gamma v_{k,\ell}\bigr).
	\end{equation}
	Furthermore, we have the convergences $u_{k,\ell}\to u_k$, $v_{k,\ell}\to v_k$, and
	$\hat y_{k,\ell}^*\to y_k^*$ as $\ell\to\infty$ by construction.
	Thus, for each $k\in\N$, we find an index $\ell(k)\in\N$ such that
	\[
			\nnorm{u_{k,\ell(k)}-u_k}\leq\frac{1}{k},\quad 
			\nnorm{v_{k,\ell(k)}-v_k}\leq\frac{1}{k},\quad
			\nnorm{x_{k,\ell(k)}^*-x_k^*}\leq\frac{1}{k},\quad
			\nnorm{\hat y_{k,\ell(k)}^*-y_k^*}\leq\frac{1}{k}.
	\]
	Let us set $\tilde u_k:=u_{k,\ell(k)}$, $\tilde v_k:=v_{k,\ell(k)}$, $\tilde x_k^*:=x_{k,\ell(k)}^*$,
	and $\tilde y_k^*:=\hat y_{k,\ell(k)}^*$ for each $k\in\N$.
	The above estimates and $u_k\to u$, $v_k\to v$, $x_k^*\to x^*$, as well as $y_k^*\to y^*$ give
	$\tilde u_k\to u$, $\tilde v_k\to v$, $\tilde x_k^*\to x^*$, as well as $\tilde y_k^*\to y^*$.
	Additionally, \eqref{eq:robustness_pseudo_coderivative_surrogate_sequences} guarantees
	\[
		\forall k\in\N\colon\quad
		\left(\tilde x_k^*,-\frac{\tilde y_k^*}{(t_k\nnorm{\tilde u_k})^{\gamma-1}}\right)
		\in 
		\widehat{\mathcal N}_{\gph\Phi}
			\bigl(\bar x+t_k\tilde u_k,\bar y+(t_k\nnorm{\tilde u_k})^\gamma \tilde v_k\bigr).
	\]
	By definition of the directional pseudo-coderivative, we find $x^*\in D^*_\gamma\Phi((\bar x,\bar y);(u,v))(y^*)$,
	and this shows the claim.
\end{proof}

\begin{proof}[Proof of \cref{prop:asymptotic_regularity_via_limiting_tools}]
	For nondirectional asymptotic regularity,
	the proof is standard and follows from a simple diagonal sequence argument.
	The proof for strong directional asymptotic regularity parallels the one 
	for directional asymptotic regularity which is presented below.

	Since one implication is clear by definition of the
	regular and limiting coderivative, we only show the other one.
	Therefore, let $\Phi$ be asymptotically regular at $(\bar x,\bar y)$ in direction $u$.
	Let us fix sequences $\{(x_k,y_k)\}_{k\in\N}\subset\gph\Phi$,
			$\{x_k^*\}_{k\in\N}\subset\mathbb X$, and 
			$\{\lambda_k\}_{k\in\N}\subset\mathbb Y$ as well as $x^*\in\mathbb X$ 
			and $y^*\in\mathbb Y$ satisfying 
			$x_k\notin\Phi^{-1}(\bar y)$, $y_k\neq\bar y$, 
			and $x_k^*\in D^*\Phi(x_k,y_k)(\lambda_k)$
			for each $k\in\N$ as well as the convergences
			\eqref{eq:convergences_directional_asymptotic_regularity}.
			For each $k\in\N$, we find sequences 
			$\{(x_{k,\ell},y_{k,\ell})\}_{\ell\in\N}\subset\gph\Phi$,
			$\{x_{k,\ell}^*\}_{\ell\in\N}\subset\mathbb X$, 
			and $\{\lambda_{k,\ell}\}_{\ell\in\N}\subset\mathbb Y$ with
			$x_{k,\ell}\to x_k$, $x_{k,\ell}^*\to x_k^*$, $y_{k,\ell}\to y_k$, 
			and $\lambda_{k,\ell}\to \lambda_k$
			as $\ell\to\infty$ as well as 
			$x_{k,\ell}^*\in \widehat D^*\Phi(x_{k,\ell},y_{k,\ell})(\lambda_{k,\ell})$
			for each $\ell\in\N$. 
			Observing that $\Phi^{-1}(\bar y)$ is closed, 
			its complement is open so that $x_{k,\ell}\notin\Phi^{-1}(\bar y)$
			holds for sufficiently large $\ell\in\N$. 
			Furthermore, since $\norm{x_k-\bar x}>0$ and $\norm{y_k-\bar y}>0$ are valid, 
			we can choose
			an index $\ell(k)\in\N$ so large such that the estimates
			\[
				\begin{aligned}
				\nnorm{x_{k,\ell(k)}-x_k}&<\frac{1}{k}\norm{x_k-\bar x},&\quad
				\nnorm{x_{k,\ell(k)}^*-x_k^*}&<\frac{1}{k},&
				\\
				\nnorm{y_{k,\ell(k)}-y_k}&
					<\frac{1}{k}\norm{y_k-\bar y},&\quad
				\nnorm{\lambda_{k,\ell(k)}-\lambda_k}&<\frac1k&
				\end{aligned}
			\]
			and $x_{k,\ell(k)}\notin\Phi^{-1}(\bar y)$ as well as $y_{k,\ell(k)}\neq\bar y$ 
			are valid. 
			For each $k\in\N$, we set $\tilde x_k:=x_{k,\ell(k)}$, 
			$\tilde x_k^*:=x_{k,\ell(k)}^*$,
			$\tilde y_k:=y_{k,\ell(k)}$, and $\tilde\lambda_k:=\lambda_{k,\ell(k)}$.
			Clearly, we have $\tilde x_k\to\bar x$, $\tilde y_k\to\bar y$, 
			$\tilde x_k^*\to x^*$, $\nnorm{\tilde \lambda_k}\to\infty$,
			$\{(\tilde x_k,\tilde y_k)\}_{k\in\N}\subset\gph\Phi$, 
			and $\tilde x_k\notin\Phi^{-1}(\bar y)$, $\tilde y_k\neq\bar y$,
			as well as $\tilde x_k^*\in\widehat D^*\Phi(\tilde x_k,\tilde y_k)(\tilde\lambda_k)$
			for each $k\in\N$ by construction.	
			Furthermore, we find
			\[
				\norm{\tilde x_k-\bar x}
				\geq
				\norm{x_k-\bar x}-\norm{\tilde x_k-x_k}
				\geq
				\frac{k-1}{k}\norm{x_k-\bar x}
			\]
			for each $k\in\N$. 
			With the above estimates at hand, we obtain
			\begin{align*}
				\norm{
					\frac{x_k-\bar x}{\nnorm{x_k-\bar x}}
					-
					\frac{\tilde x_k-\bar x}{\norm{\tilde x_k-\bar x}}
				}
				&=
				\norm{
					\frac{x_k-\tilde x_k}{\norm{x_k-\bar x}}
					+
					(\tilde x_k-\bar x)
					\left(
						\frac{1}{\nnorm{x_k-\bar x}}-\frac{1}{\nnorm{\tilde x_k-\bar x}}
					\right)
				}
				\\
				&
				\leq 
				\frac{\norm{x_k-\tilde x_k}}{\norm{x_k-\bar x}}
				+
				\frac{\norm{\tilde x_k-\bar x}\norm{x_k-\tilde x_k}}
					{\norm{x_k-\bar x}\norm{\tilde x_k-\bar x}}		
				\leq 
				\frac{2}{k}		
			\end{align*}
			and
			\begin{equation}\label{eq:estimate_on_norm_quotients}
				\begin{aligned}
				\norm{
					\frac{y_k-\bar y}{\nnorm{x_k-\bar x}}
					-
					\frac{\tilde y_k-\bar y}{\norm{\tilde x_k-\bar x}}
				}
				&=
				\norm{
					\frac{y_k-\tilde y_k}{\norm{x_k-\bar x}}
					+
					(\tilde y_k-\bar y)
					\left(
						\frac{1}{\nnorm{x_k-\bar x}}-\frac{1}{\nnorm{\tilde x_k-\bar x}}
					\right)
				}
				\\
				&
				\leq
				\frac{\norm{y_k-\tilde y_k}}{\norm{x_k-\bar x}}
				+
				\frac{\norm{\tilde y_k-\bar y}{\norm{x_k-\tilde x_k}}}
					{\norm{x_k-\bar x}\norm{\tilde x_k-\bar x}}
				\\
				&\leq
				\frac1k\frac{\norm{y_k-\bar y}}{\norm{x_k-\bar x}}
				+
				\frac{1}{k-1}\frac{\norm{\tilde y_k-y_k}+\norm{y_k-\bar y}}{\norm{x_k-\bar x}}
				\\
				&\leq
				\left(\frac{1}{k} + \frac{1}{k(k-1)} + \frac{1}{k-1}\right)
					\frac{\norm{y_k-\bar y}}{\norm{x_k-\bar x}}
				\\
				&
				=
				\frac{2}{k-1}\frac{\norm{y_k-\bar y}}{\norm{x_k-\bar x}},
				\end{aligned}
			\end{equation}
			so that, with the aid of \eqref{eq:convergences_directional_asymptotic_regularity},
			we find
			$(\tilde x_k-\bar x)/\norm{\tilde x_k-\bar x}\to u$ 
			and $(\tilde y_k-\bar y)/\norm{\tilde x_k-\bar x}\to 0$.
			With the aid of \eqref{eq:estimate_on_norm_quotients},
			\begin{align*}
				\norm{
					\frac{\nnorm{\tilde y_k-\bar y}}{\nnorm{\tilde x_k-\bar x}}
					\tilde\lambda_k
					-
					\frac{\norm{y_k-\bar y}}{\norm{x_k-\bar x}}\lambda_k
				}
				&
				\leq
				\frac{\nnorm{\tilde y_k-\bar y}}{\nnorm{\tilde x_k-\bar x}}
					\nnorm{\tilde\lambda_k-\lambda_k}
				+
				\left|
					\frac{\nnorm{\tilde y_k-\bar y}}{\nnorm{\tilde x_k-\bar x}}
					-
					\frac{\norm{y_k-\bar y}}{\norm{x_k-\bar x}}
				\right|
					\norm{\lambda_k}
				\\
				&
				\leq
				\frac1k\frac{\nnorm{\tilde y_k-\bar y}}{\nnorm{\tilde x_k-\bar x}}
				+
				\frac{2}{k-1}\frac{\norm{y_k-\bar y}}{\norm{x_k-\bar y}}\norm{\lambda_k}
			\end{align*}
			is obtained, which gives 
			$\tilde\lambda_k\nnorm{\tilde y_k-\bar y}/\nnorm{\tilde x_k-\bar x}\to y^*$.
			Similar as above, we find
			\begin{align*}
				\norm{\frac{\tilde y_k-\bar y}{\norm{\tilde y_k-\bar y}}
					-\frac{y_k-\bar y}{\norm{y_k-\bar y}}}
				\leq\frac{2}{k}
			\end{align*}
			and
			\begin{align*}
				\bigl\Vert
					\tilde\lambda_k/\nnorm{\tilde\lambda_k}
					-
					\lambda_k/\norm{\lambda_k}
				\bigr\Vert
				\leq
				2\nnorm{\lambda_k-\tilde\lambda_k}/\norm{\lambda_k}
				\leq
				2/(k\norm{\lambda_k}),
			\end{align*}
			so that \eqref{eq:convergences_directional_asymptotic_regularity} gives us
			\begin{align*}
				\lim\limits_{k\to\infty}
				\left(\frac{\tilde y_k-\bar y}{\norm{\tilde y_k-\bar y}}
					-\frac{\tilde\lambda_k}{\nnorm{\tilde\lambda_k}}\right)
				=
				\lim\limits_{k\to\infty}
				\left(\frac{y_k-\bar y}{\norm{y_k-\bar y}}
					-\frac{\lambda_k}{\norm{\lambda_k}}\right)
				=
				0.
			\end{align*}
			Now, since $\Phi$ is asymptotically regular at $(\bar x,\bar y)$ in direction $u$, 
			we obtain $x^*\in\Im D^*\Phi(\bar x,\bar y)$.
\end{proof}

\end{document}